\documentclass{amsart}
\usepackage{diagrams}
\usepackage{amssymb}
\usepackage{amsmath}
\usepackage{amsthm}
\usepackage{epsfig}
\usepackage{mathrsfs}

\input epsf
\def\relabelbox{%
  \hbox\bgroup%
}%
\def\endrelabelbox{%
\egroup%
}%
\def\relabel #1#2 {%
  \special{ps:/a {} def}%
  \smash{\rlap{#2}}%
}%
\def\adjustrelabel <#1,#2> #3#4 {%
  \special{ps:/a {} def}%
  \smash{\rlap{\kern #1 \raise #2\hbox{#4}}}%
}%
\def\extralabel <#1,#2> #3 {\smash{\rlap{\kern #1 \raise #2\hbox{#3}}}}%

\begin{document}

\newtheorem{thm}{Theorem}[subsection]
\newtheorem{lem}[thm]{Lemma}
\newtheorem{cor}[thm]{Corollary}
\newtheorem{conj}[thm]{Conjecture}
\newtheorem{qn}[thm]{Question}

\theoremstyle{definition}
\newtheorem{defn}{Definition}[subsection]

\theoremstyle{remark}
\newtheorem{rmk}{Remark}[subsection]
\newtheorem{exa}{Example}[subsection]

\def\square{\hfill${\vcenter{\vbox{\hrule height.4pt \hbox{\vrule width.4pt
height7pt \kern7pt \vrule width.4pt} \hrule height.4pt}}}$}

\def\R{\mathbb R}
\def\Z{\mathbb Z}
\def\CP{\mathbb {CP}}
\def\H{\mathbb H}
\def\F{\mathscr F}
\def\G{\mathscr G}
\def\fG{\mathfrak G}
\def\C{\mathbb C}
\def\Q{\mathscr Q}
\def\S{\mathscr S}
\def\L{\mathscr L}
\def\T{\mathscr T}
\def\hom{\text{Homeo}}
\def\u{{\text{univ}}}
\def\stb{{\text{stab}}}
\def\til{\widetilde}
\def\I{\mathscr I}
\def\N{\mathbb N}

\newenvironment{pf}{{\it Proof:}\quad}{\square \vskip 12pt}

\title{Foliations with one--sided branching}

\author{Danny Calegari}

\address{Department of Mathematics \\ Harvard University \\ Cambridge, MA 02138}
\email{dannyc@math.harvard.edu}
\date{9/12/2002. Version 0.92}
\maketitle

\begin{abstract}
We generalize the main results from \cite{wT97} and \cite{dC99b} 
to taut foliations
with one--sided branching. First constructed in \cite{gM91},
these foliations occupy an intermediate position between $\R$--covered
foliations and arbitrary taut foliations of $3$--manifolds.

We show that for a taut foliation $\F$ with one--sided branching of
an atoroidal $3$--manifold $M$, one can
construct a pair of genuine laminations $\Lambda^\pm$ of $M$ transverse to
$\F$ with solid torus complementary regions which bind every leaf of $\F$ in a geodesic
lamination. These laminations come from a {\em universal circle},
a refinement of the universal circles proposed in \cite{wT98}, which maps
monotonely and $\pi_1(M)$--equivariantly 
to each of the circles at infinity of the leaves of $\til{\F}$,
and is minimal with respect to this property. This circle is intimately
bound up with the extrinsic geometry of the leaves of $\til{\F}$. In
particular, let $\til{\F}$ denote the pulled--back foliation
of the universal cover, and co--orient $\til{\F}$ so that the leaf
space branches in the negative direction. Then for any pair of
leaves of $\til{\F}$ with
$\mu > \lambda$, the leaf $\lambda$ is asymptotic to $\mu$ in a dense
set of directions at infinity. This is a macroscopic version of an
infinitesimal result in \cite{wT98} and gives much more drastic control
over the topology and geometry of $\F$ than is achieved in \cite{wT98}.

The pair of laminations $\Lambda^\pm$ can be used to produce a pseudo--Anosov
flow transverse to $\F$ which is {\em regulating} 
in the non--branching direction. Rigidity results for $\Lambda^\pm$ in
the $\R$--covered case are extended to the case of one--sided branching.
In particular, an $\R$--covered foliation can only be deformed to a
foliation with one--sided branching along one of the two laminations
canonically associated to the $\R$--covered foliation in \cite{dC99b},
and these laminations become exactly the laminations $\Lambda^\pm$ for
the new branched foliation. 

Other corollaries include that the ambient manifold is $\delta$--hyperbolic
in the sense of Gromov, and that a self--homeomorphism of this manifold
homotopic to the identity is isotopic to the identity.
\end{abstract}

\section{Introduction}
This paper extends the main ideas and theorems of \cite{dC99b} to
taut foliations with one--sided branching. These foliations are like
a hybrid between $\R$--covered foliations (whose leaf spaces in the
universal cover exhibit no branching) and arbitrary taut foliations, whose
leaf spaces may branch in either direction. Foliations with one--sided
branching do not exhibit the full range of phenomena to be observed in
arbitrary taut foliations, but on the other hand one has enough control
over their geometry to prove some very powerful structure theorems.

The main results of this paper are summarized in the abstract. The most
important is the fact that transverse to any taut foliation with one--sided
branching of an atoroidal $3$--manifold, one can construct 
a pair of genuine laminations $\Lambda^\pm$ 
which have solid torus guts, and which bind each leaf of the foliation.
These laminations can be ``blown down'' to a pseudo--Anosov flow transverse
to the foliation which is {\em regulating} in the nonbranching direction;
that is, flow lines in the universal cover are properly embedded in the
leaf space in the nonbranching direction. There is a history of finding
such a structure transverse or almost transverse
to taut foliations of a certain quality, which we recapitulate here:

\begin{itemize}
\item{Surface bundles over $S^1$ \cite{wT98b}}
\item{Finite depth foliations \cite{lM00} based on unpublished work of Gabai}
\item{Foliations arising from slitherings over $S^1$ \cite{wT97}}
\item{$\R$--covered foliations \cite{dC99b}, \cite{sF00}}
\end{itemize}

The case of one--sided branching is dealt with in this paper, leaving only
the case of a taut foliation with two--sided branching
still open. An approach to the two--sided case is promised in
\cite{wT98}, but this paper is largely unwritten, and details are not available.
A great deal of this program has been carried out for
foliations $\F$ with two--sided branching
{\em and every leaf dense}, which we
describe in \cite{Calprom} and its forthcoming sequel.
In particular, we construct a pseudo--Anosov flow transverse to such an $\F$
which can be split open to a pair of genuine laminations which are
transverse to each other, and intersect leaves of $\F$ in geodesic
laminations.

The existence of such a structure is important for a number of reasons.
Firstly, it sets the stage for an attempt to prove the geometrization
conjecture for the underlying manifolds, along the lines of \cite{wT98b},
and as outlined in \cite{dC99b}. Secondly, the underlying structure of
the laminations constructed is much more rigid than the foliations themselves.
This rigidity leads one to hope that the structure of a 
{\em pseudo--Anosov package} (to coin a phrase) might be described in 
finite combinatorial terms and lead to a closer relationship of the theory
of taut foliations with the algorithmic theory of $3$--manifolds; moreover,
this rigidity might manifest itself in subtler ways through
numerical invariants attached to the structure. Finally, there are immediate
corollaries ($\delta$--hyperbolicity, non--existence of exotic 
self--homeomorphisms) for the underlying manifold.
For a precise definition of the term ``pseudo--Anosov package'', see 
\cite{CalPhD} or \cite{Calprom}.

The laminations $\Lambda^\pm$ reflect the disparity between the intrinsic and the
extrinsic geometry of the leaves of $\til{\F}$ in $\til{M}$. Let $L$ denote
the leaf space of $\til{\F}$, and assume that $L$ branches in the negative, but not
the positive direction. Roughly speaking,
for each leaf $\lambda$ of $\til{\F}$, the lamination $\til{\Lambda}^+ \cap \lambda$
parameterizes the branching of the subset of $L$ consisting of leaves
on the negative side of $\lambda$.
This branching reflects the fact that leaves of $\til{\F}$ are not uniformly properly
embedded in the subspace of $\til{M}$ contained above them.
Conversely, the lamination $\til{\Lambda}^- \cap \lambda$ is a kind of
{\em eigenlamination of distortion}, and consists of geodesics in $\lambda$ which are worst
at efficiently measuring distance in the subspace of $\til{M}$ contained below $\lambda$.
This pair of laminations is parameterized by a pair of laminations 
$\Lambda^\pm_\u$ of a
{\em master circle} $S^1_\u$ on which $\pi_1(M)$ acts, and which comes equipped with
canonical monotone maps to $S^1_\infty(\lambda)$ for each leaf $\lambda$ of $\til{\F}$.

\subsection{Acknowledgments}

The material in this paper fits into Thurston's program, outlined in \cite{wT98},
to unify the theory of geometric structures and taut foliations on $3$--manifolds.
I am indebted to Bill for explaining his program to me and encouraging my participation
in it. 

I would also like to thank the referee, who made many very interesting and
insightful comments which have led to dramatic improvements in many aspects of
this paper.

\section{Extrinsic geometry of leaves}

\subsection{Definitions and conventions}

In the sequel, we will denote by $\F$ a taut foliation of a $3$--manifold $M$.
We assume $\F$ is not covered by the product foliation of $S^2 \times S^1$,
and therefore $M$ is irreducible.
We assume $M$ is orientable and $\F$ is co--orientable. Furthermore, we
assume $M$ is atoroidal. We typically assume $M$ is closed, although our
results generally extend to the case that $M$ has boundary consisting of 
a finite collection of tori. Let $\til{M}$ denote the universal cover of $M$,
and $\til{\F}$ the pullback of $\F$ to $\til{M}$. The leaf space of $\F$ will
be denoted $L$. It is well--known that $L$ is a simply--connected, possibly
non--Hausdorff $1$--manifold (see e.g. \cite{cP78}).
There is a natural partial order on $L$ coming from the
co--orientation on $\F$. 

\begin{defn}
A taut foliation $\F$ has {\em one--sided branching in the positive direction} 
if any two $\mu,\lambda \in L$ have a common lower bound in the partial order
on $L$. Conversely, $\F$ has {\em one--sided branching in the 
negative direction} if
any two $\mu,\lambda \in L$ have a common upper bound. Clearly, the direction
of branching depends on the choice of co--orientation on $\F$. If $L$ is
totally ordered, $L = \R$ and $\F$ is called {\em $\R$--covered}.
\end{defn}

Notice that for $\F$ not co--orientable, $L$ must either be $\R$--covered or
must branch in both directions.

The following easy theorem is found in \cite{dC99c}.

\begin{thm}
Suppose $M^3$ is hyperbolic, and contains a taut foliation $\F$
with one--sided branching. Then for any leaf $\lambda$ of $\til{\F}$,
the two connected components of $\til{M^3}-\lambda = \H^3 - \lambda$
cannot both contain an open half--space.
\end{thm}

\begin{cor}\label{no_compact_leaf}
If $\F$ is a taut foliation with one--sided branching
of an atoroidal $M^3$ then $\F$ contains no compact leaves.
\end{cor}

The following theorem is found in \cite{aC93}:
\begin{thm}[Candel]
Let $\Lambda$ be a lamination of a compact space $M$
with $2$--dimensional leaves.
Suppose that every invariant transverse measure supported on $\Lambda$ has
negative Euler characteristic. Then there is a metric on $M$
such that the inherited path metric makes the leaves of $\Lambda$ into
Riemann surfaces of constant curvature $-1$.
\end{thm}

The conditions of Candel's theorem are automatically satisfied for $\F$
a taut foliation of an atoroidal irreducible $3$--manifold $M$, so we assume
in the sequel that $M$ has a metric such that all the leaves of $\F$ have
constant curvature $-1$.

\begin{defn}
For any foliation, a {\em saturated set} is a union of leaves. A
{\em minimal set} is a nonempty closed saturated set, minimal with respect
to this property.
\end{defn}

Clearly, minimal sets exist. It is not at all unusual for an entire
foliated manifold $M$ to be a minimal set. The salient feature of a minimal
set $\Lambda$ is that every leaf $\lambda$ in $\Lambda$ has 
$\overline{\lambda} = \Lambda$. We will need in the sequel
the following lemma about the geometry of minimal sets.

\begin{lem}\label{big_ball_dense}
For each minimal set $\Lambda \subset \F$ and each $\epsilon$, there is a
$t$ such that each ball of radius $t$ in a leaf of $\Lambda$ is an
$\epsilon$--net for $\Lambda$.
\end{lem}
\begin{pf}
Suppose not. Then choose a sequence of real numbers $t_i \to \infty$ so that
there are a sequence of points $p_i$ such that the $t_i$--balls about
$p_i$ do not come within $\epsilon$ of some $q_i \in \Lambda$. By
compactness of $\Lambda$, we can find a convergent subsequence for which $p_i \to p$ and
$q_i \to q$ so that
the closure of the leaf through $p$ does not come within $\epsilon$ 
of $q$, violating minimality.
\end{pf}

\begin{defn}
A one dimensional foliation $X$ transverse to a taut
foliation $\F$ {\em regulates $\epsilon$--neighborhoods}
if, after lifting to the universal
cover, the restriction of $\til{X}$ (the foliation covering $X$) to the
closed $\epsilon$--neighborhood of any leaf $\lambda$ 
of $\til{\F}$ is a product foliation of $\lambda \times I$ by closed intervals
$\text{point} \times I$.
\end{defn}

For any taut $\F$, any transverse $1$--dimensional foliation
$X\pitchfork \F$ regulates $\epsilon$--neighborhoods,
for some $\epsilon$. A pair of leaves $\mu,\lambda$ of $\til{\F}$ which are sufficiently
close at some pair of points $p \in \mu,q \in \lambda$ are close on a big neighborhood.
That is, for any $T>0$ and any $\epsilon>0$
there is a $\delta>0$ such that if $d_{\til{M}}(p,q)<\delta$ then
the $T$ neighborhood of $p$ in $\mu$ is $\epsilon$--close to the $T$ neighborhood
of $q$ in $\lambda$. Moreover, by choosing $\delta$ possibly slightly smaller (but
depending only on $T,\epsilon,X,\F$),
the flow along integral leaves of $\til{X}$ from
$N_T(p)$ for time $\epsilon$ hits $\mu$ and defines a bilipschitz map to the
image, which can be taken to be contained in $N_{T+\epsilon}(q)$ and to contain
$N_{T-\epsilon}(q)$. For any $T$, for sufficiently small $\epsilon$, this flow can be
taken to be as close to an isometry as desired. For instance, if $M$ is Riemannian
$X$ could be chosen to be the integral foliation of the
normal bundle to $T\F$.

For definitions and some basic results about genuine laminations of
$3$--manifolds, see \cite{dGuO87}, \cite{dGwK93} and \cite{dGwK98}.
We recapitulate some of the main notions here for convenience.

\begin{defn}
A {\em lamination} in a $3$ manifold is a foliation of a closed subset of $M$
by $2$ dimensional leaves. The complement of this closed subset falls into
connected components, called {\em complementary regions}.
A lamination is {\em essential} if it contains
no spherical leaf or torus leaf bounding a solid torus, and furthermore if
$C$ denotes the closure (with respect to the path metric in $M$)
of a complementary region, then $C$ is irreducible and
$\partial C$ is both incompressible and {\em end incompressible} in $C$.
Here an end compressing disk is a properly embedded $(D^2 - (\text{closed
arc in } \partial D^2))$ in $C$ which is not properly isotopic rel $\partial$
in $C$ to an embedding into a leaf. Finally, an essential lamination is
{\em genuine} if it has some complementary region which is not an $I$--bundle.
\end{defn}

A complementary region to a genuine lamination falls into two pieces:
the {\em guts}, which carry the essential topology of the complementary 
region, and the {\em interstitial regions}, which are just $I$ bundles 
over non--compact surfaces, which get thinner and thinner as they go 
away from the guts. The interstitial
regions meet the guts along annuli. Ideal polygons can be properly embedded
in complementary regions, where the cusp neighborhoods of the ideal points
run up the interstitial regions as $I \times \R^+$. An end compressing disk
is a properly embedded ideal monogon which is not isotopic
rel. $\partial$ into a leaf.

\subsection{One--sided uniform properness}

\begin{defn}
An embedding $e$ of a metric space $X$ in a metric space $Y$ is {\em uniformly
proper} if there is an increasing function $f:\R^+ \to \R^+$ such that
$f(x) \to \infty$ as $x \to \infty$ with the property that for all $a,b \in X$
$$d_Y(e(a),e(b)) \ge f(d_X(a,b))$$
A subspace $X$ of a path metric space $Y$ is uniformly properly embedded if
the embedding of $X$ with its inherited path metric into $Y$ is uniformly
proper.
\end{defn}

\begin{defn}
Let $\F$ be a co--oriented foliation of $\R^3$ by separating planes. Leaves
of $\F$ are {\em uniformly proper below} if each leaf $\lambda$ is
uniformly properly embedded in the subspace $X^\lambda \subset \R^3$ consisting
of leaves below and including $\lambda$, with respect to the induced path
metrics. $\F$ is {\em uniformly proper below} if each leaf $\lambda$ 
is uniformly proper below for some function $f: \R^+ \to \R^+$ {\em independent} of
$\lambda$. Define uniformly proper above in an analogous way.
\end{defn}

In the sequel, we will adhere to the convention that our foliations have
one--sided branching {\em in the negative direction}, unless we explicitly
say otherwise.

\begin{lem}\label{uniformly_proper_below}
Let $\F$ be a taut foliation of $M$ which has one--sided branching in the
negative direction. Then $\til{\F}$ is uniformly proper below, with respect
to the path metric on $\til{M} = \R^3$ inherited from some path metric on $M$.
\end{lem}
\begin{pf}
Suppose to the contrary that there is a sequence of pairs of points
$\lbrace p_i,q_i \rbrace$ 
contained in leaves $\lambda_i$ such that there are arcs
$\alpha_i$ contained below $\lambda_i$ with length $|\alpha_i|\le t$ but
the leafwise distance from $p_i$ to $q_i$ goes to $\infty$. We can
translate the pairs $\lbrace p_i, q_i \rbrace$ by elements of $\pi_1(M)$
to get a convergent subsequence, which we also denote by $\lbrace p_i,q_i
\rbrace \to \lbrace p,q \rbrace$ joined by paths $\alpha_i$ below $\lambda_i$
which converge to a path $\alpha$ below the leaves $\lambda^\pm$ containing
$p$ and $q$. If $\lambda^+ = \lambda^-$, the leafwise distances between
$p_i$ and $q_i$ would have been uniformly bounded, contrary to assumption.
Therefore $\lambda^+ \ne \lambda^-$, and we have exhibited branching
of $L$ in the positive direction, also contrary to assumption. 

It follows that for each $t$ and each pair of points $p_i,q_i \in \lambda_i$
which can be joined by an arc below $\lambda_i$ of length at most $t$,
the leafwise distance between $p_i,q_i$ is bounded by some $f(t)$. That
is, $\til{\F}$ is uniformly proper below.
\end{pf}

For $\mu,\lambda$ leaves of $\til{\F}$ with $\lambda > \mu$, denote by
$\til{M}_\mu$ the subspace of $\til{M}$ contained above $\mu$, by
$\til{M}^\lambda$ the subspace contained below $\lambda$, and by
$\til{M}_\mu^\lambda$ the subspace contained between these two leaves.
Note that {\em in this context},
by ``above'' and ``below'' we mean the connected components
of the complement of the leaf in question. Every point $p$ in
a leaf $\lambda > \mu$ is ``above'' $\mu$, but not every point above
$\mu$ lies on a leaf comparable to $\mu$.

Notice that leaves of $\til{\F}$ are not in general
uniformly properly embedded in the subspaces contained above them. On
the other hand, since $M$ is compact, there is a uniform $\epsilon$ such
that leaves of $\til{\F}$ are quasi--geodesically embedded in their
$\epsilon$--neighborhoods.

\begin{lem}\label{straighten_transversals}
Let $\gamma$ be an arc in $\til{M}$ whose endpoints lie on comparable leaves
$\lambda,\mu$ and such that $\gamma$ is contained in $\til{M}^\lambda$.
Then there is a $g:\R^+ \to \R^+$ such that if $|\gamma|\le t$, $\gamma$
can be canonically straightened to a monotone $\gamma'$
or a transversal with the same endpoints of length $\le g(t)$.
\end{lem}
\begin{pf}
Wiggle $\gamma$ to be generic, with at most one critical point on each
leaf. Orient $\gamma$ from the point on $\mu$ to the point on $\lambda$.
Since $\gamma \subset \til{M}^\lambda$, it can be straightened to a transversal
$\gamma'$ with the same endpoints by pushing in a collection of
{\em disjoint} subarcs $\gamma_i$ whose endpoints lie on the same
leaf $\lambda_i$ and which are each contained in $\til{M}^{\lambda_i}$.

\begin{figure}[ht]
\centerline{\relabelbox\small \epsfxsize 2.0truein
\epsfbox{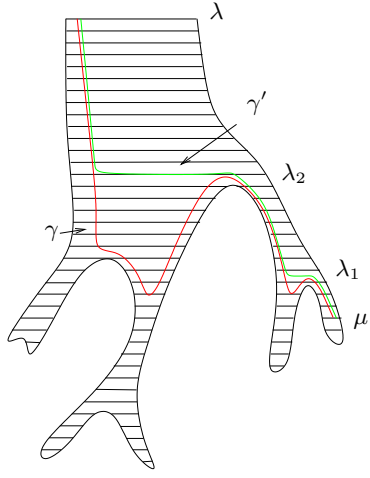}
\relabel {g}{$\gamma$}
\relabel {gp}{$\gamma'$}
\relabel {lambda}{$\lambda$}
\relabel {lambda1}{$\lambda_1$}
\relabel {lambda2}{$\lambda_2$}
\relabel {mu}{$\mu$}
\endrelabelbox}
\caption{An arc $\gamma$ in a foliation with one--sided branching
which intersects leaves below one endpoint can be straightened by
pushing in subarcs $\gamma_i$ which start and end on the same leaf $\lambda_i$
and only intersect leaves below $\lambda_i$.}
\end{figure}

To see this, first define a collection of subarcs
$\lbrace \gamma_i \rbrace$ to be the arcs whose initial points
are the local maxima of $\gamma$, and whose
end points are the first return of $\gamma$ to the leaf $\lambda_i$
containing the initial point of $\gamma_i$. By construction, two
subarcs $\gamma_i,\gamma_j$ intersect iff they are nested. Refine
the collection $\lbrace \gamma_i \rbrace$ by discarding innermost
subarcs until the $\gamma_i$ are all disjoint.

Since leaves of $\til{\F}$ are uniformly properly embedded below, and since
these subarcs all have length $\le t$, each $\gamma_i$ can be
straightened to the leafwise geodesic arc $\gamma_i'$ in $\lambda_i$ with
the same endpoints. This gives a {\em canonical} map from
$\gamma$ to a monotone arc $\gamma'$ with the same endpoints.
By assumption, $|\gamma|\le t$, so the arcs in $\gamma \cap \gamma'$
have total length $\le t$.
Moreover, the length of each subarc $\gamma_i'$ is
uniformly comparable to the length of $\gamma_i$. For, leaves of $\til{\F}$
are continuously varying hyperbolic planes, so leafwise geodesics between pairs of
points are unique and vary continuously; we will use this fact repeatedly.
It follows that the length of $\gamma'$ can be controlled in terms of $t$. Then $\gamma'$
can be perturbed (not canonically) an arbitrarily small amount to be transverse.
\end{pf}

\begin{lem}\label{quasiisometric_in_neighborhood}
For each $t$ there is a $k_t,\epsilon_t$ such that for each
leaf $\lambda$ of $\til{\F}$ the leaf $\lambda$ is 
$(k_t,\epsilon_t)$ coarsely quasi--isometrically embedded in
$N_t(\lambda) \cap \til{M}^\lambda$.
\end{lem}
\begin{pf}
Let $p,q \in \lambda$ for some $\lambda$, and suppose $\gamma$ is
a path in $N_t(\lambda) \cap \til{M}^\lambda$ joining $p$ to $q$.
We can decompose $\gamma$ into approximately $|\gamma|$
segments of length $\le 1$. Denote these segments $\gamma_i$. Denote
the two endpoints of $\gamma_i$ by $r_i,s_i$. Then $r_i$ and $s_i$
can be joined to points $p_i,q_i$ in $\lambda$ by paths of length
$\le t$. It follows that
$$d_\lambda(p_i,q_i) \le f(2t+1)$$
by lemma~\ref{uniformly_proper_below}
and therefore $p,q$ can be joined by a path in $\lambda$ of length
approximately $|\gamma|f(2t+1)$ with an error term of size $f(2t+1)$.
\end{pf}

\begin{defn}
A maximal closed saturated subset $P$ of 
$\til{M}$ whose interior is foliated as a product is called a {\em pocket}. If
$P$ is homeomorphic to $[0,1] \times \R^2$ we call it an {\em inessential
pocket}. Otherwise we call it an {\em essential pocket}.
\end{defn}

\begin{lem}\label{neighborhood_is_product}
Let $\F$ be a taut foliation of $M$ which has one--sided branching in the
negative direction. Suppose $P$ is foliated by leaves with a unique
uppermost leaf $\lambda$ and a collection of lowermost leaves
$\mu_i$, such that there is a $\delta$ with
$\mu_i \subset N_\delta(\lambda)$ for all $i$. Then $P$
is an inessential pocket and is quasi--isometric with its
inherited path--metric to either of its boundary components $\lambda$ or
$\mu$.
\end{lem}
\begin{pf}
By lemma~\ref{quasiisometric_in_neighborhood}, for any $\delta$, 
the inclusion of $\lambda$ into $\til{M}^\lambda \cap N_\delta(\lambda)$ 
is a quasi--isometry. We construct a homotopy 
$$R:\bigcup_i \mu_i \times I \to P$$
such that $R(\cdot,1):\bigcup_i \mu_i \to \lambda$, and the arcs $R(p,I)$ are 
transversals of length uniformly bounded by some $t$. For ease of notation,
we define $r(\cdot) = R(\cdot,1)$.

By lemma~\ref{straighten_transversals}, for every $p \in \mu_i$ there is
a transversal $\tau_p$ from $p$ to $\lambda$ of uniformly bounded length.
We let $T$ be a triangulation of $\bigcup_i \mu_i$ with small geodesic
simplices, and choose such a set of transversals for the vertices $T^0$ of
$T$. This defines $R|_{T^0}$. For an edge $e$ of $T^1$ between $p_1,p_2$,
the transversals $\tau_{p_1},\tau_{p_2}$ intersect the same leaves, and
therefore they can be joined by leafwise geodesic arcs. Parameterizing these
by arclength, this defines $R|_{T^1}$. Finally, $R$ can be extended over the
$2$--simplices by barycentric extension. Notice that for $\nu$ a leaf
intermediate between $\mu_i$ and $\lambda$, if $\tau_{p_1},\tau_{p_2}$
intersect $\nu$ at $r_1,r_2$ then the leafwise
distance between $r_1$ and $r_2$ in $\nu$ is uniformly bounded, since
$r_1$ and $r_2$ are joined by subarcs of $\tau_{p_1},\tau_{p_2}$ and the
edge in $\mu_i$ between $p_1,p_2$, and $\nu$ is uniformly properly embedded below
by lemma~\ref{uniformly_proper_below}. It follows that the fibers of $R$
have the desired properties.

Note that $r:\mu_i \to \lambda$ is {\em proper} for each $\mu_i$, since
both leaves are properly embedded in $\til{M}$, and $R$ only moves points
a bounded distance. So the degree of $r$ on each $\mu_i$ is well--defined.

Let $q \in \lambda$ and $B_s(q) \subset \lambda$ be
the disk in $\lambda$ of (leafwise) radius $s$ about $q$. Suppose there
is a disk $D \subset \mu_i$ such that $B_s(q) \subset r(D)$ and
$r(\partial D)$ winds around $B$ with nonzero degree. In this
case we say that $r(D)$ {\em covers} $B$ with degree equal to the winding
number of $r(\partial D)$ about $B$. Then $R(\partial D, I)$
separates $\mu_j$ from $q$. Moreover, since $\lambda$ is quasi--isometrically
embedded in $P$, $R(\partial D,I)$ is a definite distance away from $q$.
It follows that there is a {\em uniform} $s$ such that if $q,D$ exist with the
properties above, $q$ is not in the image $r(\mu_j)$ for any $j \ne i$.
But if $r|_{\mu_i}$ has nonzero degree for some $i$, then for {\em every} point
$q \in \lambda$ there is a $D \subset \mu_i$ which covers $B_s(q)$. It follows
that if $r|_{\mu_i}$ has nonzero degree for some $i$, then there are no $\mu_j$
with $j \ne i$. In this case, the region $P$ does not branch, and is foliated
as a product.

It suffices to show that there is some $\mu_i$ which is mapped by $r$ with
nonzero degree.

We make the observation that if $\nu$ is an intermediate leaf in $P$ bounded
below by some subset $\lbrace \mu_j \rbrace$ of the $\mu_i$, then the subset 
$P_\nu = \til{M}^\nu \cap P$ satisfies the hypothesis of the lemma.
Let $r_\nu:\mu_j \to \nu$ be the retraction
defined by the appropriate level $R(\cdot,t)$ of $t$. The degree of $r_\nu$
restricted to $\mu_j$ is well--defined, of course.
Let $\mu_j$ be a lowermost leaf, and let $\nu$ be a leaf which is within $\epsilon$
of $\mu_j$ at some point $q \in \nu$. If $\epsilon$ is sufficiently small, $\mu_j$
and $\nu$ are very close on big compact subsets, so there is a big disk $D \subset \mu_j$
which covers $B_s(q)$ with degree one. If $r_\nu:\mu_j \to \nu$ has
degree $\ne 1$, there is some point $p \in \mu_j \backslash D$
which maps to $q$. But $p$ is separated from $q$ by $R(\partial D,I)$, which is
a definite distance from $q$ if $\epsilon$ is small enough,
since $\nu$ is quasi--isometrically embedded in $P_\nu$. 
This is a contradiction. It follows that
$r_\nu$ has degree one. In particular, $r_\nu:\mu_j \to \nu$ is onto, so $\nu$
is contained within some uniform distance of $\lambda$ too.
If $\tau$ is a transversal from $\mu_j$ to $\lambda$,
we can break up $\tau$ into a finite number of intervals each of length less
than $\epsilon$. Let $\nu_i$ be the leaves of $\til{\F}$ at the endpoints of these
subintervals. It follows that in the sequence of maps
$$\mu_j \to \nu_1 \to \nu_2 \to \dots \to \nu_n \to \lambda$$
each map has degree one. Moreover, the composition $\mu_j \to \lambda$ is just $r$.
So $r$ has degree one, and $P$ is an inessential pocket. 
\end{pf}

\begin{exa}
If we make no restriction on the geometry of leaves of $\F$, it is possible to
produce a region $P$ foliated by $\F$ with one--sided branching which is
bounded above by a single leaf, which is bounded below by more than one leaf,
and which is quasi--isometric as a path metric space to the topmost leaf.

\begin{figure}[ht]
\centerline{\relabelbox\small \epsfxsize 3.0truein
\epsfbox{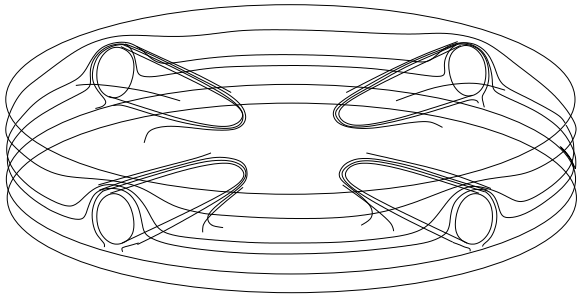}
\endrelabelbox}
\caption{This configuration cannot be realized as the topology of
an essential pocket in $\til{\F}$.}
\end{figure}

The pocket $P$ is obtained geometrically from $\H^2 \times I$ by drilling out neighborhoods
$N(r_i)$ of a finite number of geodesic rays $r_i$ in $\H^2 \times \frac 1 2$. The foliation
$\F$ looks like a product near the top leaf $\H^2 \times 1$.
As it moves in the negative direction,
the leaves bend around in a neighborhood of the disks $\partial N(r_i)$. The
interior of $P$ is foliated as a product, and the lowermost leaves consist of the
$\partial N(r_i)$ and the lowermost leaf $\H^2 \times 0$. See figure 2.
\end{exa}

\begin{thm}\label{minimal_set_inessential}
Let $\F$ be a taut foliation with one--sided branching. Let $\Lambda$ be
a minimal set. Then the path--metric completion of
each complementary region to $\Lambda$ is an
inessential pocket.
\end{thm}
\begin{pf}
Let $P$ be the closure of a complementary region to $\til{\Lambda}$. Then
$P$ has an uppermost leaf $\lambda$, and some lowermost leaves $\mu_i$.

Suppose $P$ is not an inessential pocket. Then by 
lemma~\ref{neighborhood_is_product} we can conclude that $\mu_i$ is
not contained in a bounded neighborhood of $\lambda$ for some $i$.
Let $p \in \mu_i$ be chosen with $d(p,\lambda) \ge s$. Join $p$ to $\lambda$
by a transversal $\tau$ with $|\tau|=t$. Let $\sigma$ be a small transversal
emanating in the negative direction from $p$ with length $a$. By
lemma~\ref{big_ball_dense} there is a uniform $r$ so that any point
in $\Lambda$ can be joined by a leafwise arc of length $\le r$ to some
point within $\epsilon$ of $\pi(p)$, the projection of $p$ to $M$.
Since $\mu_i$ is isolated in $\til{\Lambda}$
on the positive side, any point in $\Lambda$ can be joined by a leafwise
arc of length $\le r$ to some point within $\epsilon$ of $\pi(p)$ on the
{\em negative} side of $\mu_i$. It follows that we can arrange for any point in
$\Lambda$ to be joined by a leafwise arc of length $\le r$ to some point
in the projection of $\sigma$.

Now pick another point $q$ in $\mu_i$ which is not contained within
distance $u$ of $\lambda$, where $u \gg t,s$. Then we can join $q$ by an arc
in $\mu_i$ of length $\le r$ to some point $p' \in \alpha(\sigma)$ for
some $\alpha$. Now
$\alpha(p) \in \til{\Lambda}$ since $p \in \til{\Lambda}$.
On the other hand, $\alpha(p)$ is on the positive side of
$\mu_i$, by construction, and therefore $\alpha(p) \ge \lambda$.
Furthermore, $\alpha(p)$ is within distance $a$ of $p'$, and
therefore within distance $a+r$ of $q$. This violates the
choice of $q$, since $q$ is at least distance $u$ from $\lambda$.
This contradiction establishes the theorem.
\end{pf}

\begin{rmk}\label{short_special_case}
We do not actually need the full strength of lemma~\ref{neighborhood_is_product}
to prove this theorem. In the context of this theorem, a short proof of
this lemma is as follows:

The region $P$ is the universal cover of some complementary region
$N$ to $\Lambda$ in $M$. Since $\F$ is co--oriented, $\partial N$
consists of a single uppermost leaf $l$ covered by $\lambda$,
and some lowermost leaves $m_i$ covered by the $\mu_j$ (with possibly
different indices). We would like to show that $N$, and therefore
$P$, is a product region between $l$ and a single leaf $m$.
If not, then $N$ can be decomposed nontrivially into a {\em gut region}
$G$ and a collection of interstitial regions $R_i$ which are
$I$--bundles over noncompact surfaces $\Sigma_i \subset \cup_j m_j$.
Since $l$ is the unique uppermost leaf in the pocket, every
$R_i$ interpolates between the $\cup_j m_j$ and $l$. It follows
that if $N$ is not a product region, the gut $G$ is not homotopic
into $l$ by a homotopy which collapses $\partial G$ along the $I$
fibers. But in this case, there is some $\alpha \in \pi_1(G)$ which
is not conjugate into $\pi_1(l)$. However, $\lambda$ is the unique
uppermost leaf of $P$, and is therefore stabilized by $\pi_1(N)$.
It follows that $G$ is homotopic into $l$ as indicated. A sublamination of
a taut foliation is essential, so complementary regions are irreducible; therefore
$N$ is a product pocket, by standard $3$--manifold topology.
It follows that $P$ has a single lowermost
leaf $\mu$, which covers $m$, and $P$ is an inessential pocket.
\end{rmk}

\begin{rmk}
One--sidedness is important for this theorem. If $\Lambda$ is a co--orientable
genuine lamination of a manifold $M$ and $P$ is the cover of a complementary region,
every lowermost leaf of $\partial P$ is contained in a bounded neighborhood of
a {\em union} of uppermost leaves, and conversely. For, if $P$ covers a complementary
region $N$ downstairs, $N$ decomposes into a gut $G$ which is compact, and
interstitial regions $R_i$ of bounded thickness.
Either $N$ is a compact product, or every point
in a lowermost leaf (say) can be joined by a path of bounded length (less than
the diameter of $G$ to an interstitial region. The interstitial regions run from
lowermost leaves to uppermost leaves; they have bounded thickness, so this
path can be joined by a vertical arc of some $R_i$ of bounded length
to reach an uppermost leaf. On the other hand, if some
lowermost leaf is contained in a bounded neighborhood of an uppermost leaf or
vice versa, then the argument above would show that $P$ is an inessential
pocket.
\end{rmk}

\begin{cor}
If $\F$ is a taut foliation with one--sided branching, then $\til{\F}$ contains
no essential pockets.
\end{cor}
\begin{pf}
The maximality of a pocket implies that it is contained in the
complement of a minimal set. Therefore such a pocket is inessential.
\end{pf}

In the sequel we refer to inessential pockets merely as pockets, since there
are no other kinds.

Since $\F$ has branching, there is a uniform $t$ such that no pocket
contains an embedded ball of radius $t$.

\begin{thm}\label{blow_down_pockets}
A taut foliation $\F$ with one--sided branching is obtained from another
taut foliation $\G$ with one--sided branching and every
leaf dense (and therefore no pockets) by blowing up and 
perturbing some leaves.
\end{thm}
\begin{pf}
Let $P$ be a complementary region to $\til{\Lambda}$, where $\Lambda$ is
a minimal set for $\F$. By theorem~\ref{minimal_set_inessential}, $P$ is
an inessential pocket with highest leaf $\lambda$ and lowest leaf $\mu$,
and $\mu \subset N_\delta(\lambda)$ for some $\delta$. We claim $\lambda \subset
N_\delta(\mu)$ for some possibly different $\delta$. Repeating the argument in
lemma~\ref{neighborhood_is_product}, the map
$R(\cdot,1):\mu \to \lambda$ constructed in that lemma is degree
one, and therefore onto. On the other hand, each point in $\mu$ is moved by $R$ along
an arc in $P$ of uniformly bounded length. It follows that every point in
$\lambda$ is joined to a point in $\mu$ by an arc of uniformly bounded length, and the
claim is proved. By lemma~\ref{quasiisometric_in_neighborhood}, the pocket
$\til{M}^\lambda_\mu$ is quasiisometric to $\lambda$. On the other hand, for each
leaf $\nu \subset P$, the pocket $\til{M}^\nu_\mu$ satisfies $\mu \subset
N_\delta(\nu)$. For, $\nu$ separates $\mu$ from $\lambda$, and
$\mu \subset N_\delta(\lambda)$.
We can re--parameterize the $t$ variable in $R(\cdot,t)$ so that at time $1$,
$R(\cdot,1):\mu \to \nu$. By lemma~\ref{neighborhood_is_product}
$R(\cdot,1):\mu \to \nu$ is degree one and therefore onto. It
follows that the image of $\mu$ sweeps through all of $P$ under the homotopy $R$. 
So $P \subset N_\delta(\mu)$ for some $\delta$.

We want to show that $\mu$ is quasi--isometrically embedded in $P$. Since
$P \subset N_\delta(\mu)$, it follows that {\em either} $\mu$ is quasi--isometrically
embedded, or it is not uniformly properly embedded. In the latter case,
there are a sequence of pairs of points
$\lbrace p_i,q_i\rbrace$ in $\mu$ with $d_P(p_i,q_i) 
= \text{const.}$ but $d_\mu(p_i,q_i) \to \infty$. After choosing a subsequence and
translating by elements of $\pi_1(M)$, we may assume $\lbrace p_i,q_i\rbrace \to
\lbrace p,q\rbrace$ lying on distinct leaves $\mu_1,\mu_2$.
Since leaves are quasi--isometrically embedded in their $\epsilon$--neighborhoods
for {\em some} $\epsilon$, it follows that $\mu$ stays a definite distance away
from $\lambda$ near $p_i,q_i$ and therefore $p,q$ are contained in leaves of
$\til{\Lambda}$ which are isolated in $\til{\Lambda}$ on the positive side.

We can join $p_i,q_i$ by an arc $t_i \subset P$, and arrange for $t_i$ to be bounded
away from $\lambda$. The limit $t$ will be bounded away from $\til{\Lambda}$. It follows
that $\mu_1,\mu_2$ are lowermost leaves of a complementary region to $\til{\Lambda}$,
contrary to lemma~\ref{minimal_set_inessential}. This contradiction establishes that
$\mu$ is quasi--isometrically embedded in $P$, and therefore $R(\cdot,1):\mu \to \lambda$ is
a quasi--isometry. 

As in \cite{dC99b}, $P$ covers a complementary region to $\Lambda$ which is a
(possibly perturbed) blow--up of some leaf of a foliation $\G$, 
which can be recovered by simultaneously blowing down all such pockets.
In more detail: for each (inessential) pocket $P$ bounded by leaves
$\mu,\lambda$, we know that 
$$\text{stab}(\lambda) = \text{stab}(\mu) = \text{stab}(P)$$
and $P$ can be parameterized $\text{stab}(P)$--equivariantly as
$\mu \times I$ where the $I$ factors are of bounded length (and necessarily
must taper off towards each end of $\mu/\text{stab}(P)$). So we can
include $P/\text{stab}(P)$ homeomorphically to $M$ and take the quotient
$\mu \times I \to \mu$. This can be done simultaneously for all pockets, and
the quotient is homeomorphic to $M$ and the image of $\Lambda$ is a
minimal taut foliation $\G$. The foliation $\G$ has one--sided
branching if $\F$ has, since the actions on the leaf space of
the universal cover are monotone equivalent.
\end{pf}

\subsection{Action on the leaf space}

The main goal of this section is to analyze the compressibility of
the action of $\pi_1(M)$ on the leaf space $L$ of $\til{\F}$. 

\begin{defn}
An interval $I \subset L$ is {\em incompressible} if
there is no $\alpha$ which takes $I$ to a proper
subset of itself.
\end{defn}

If $\F'$ is obtained from $\F$ by blowing up a leaf, no interval of $L$ 
corresponding to a blown--up pocket can be taken to a proper subset of itself by
the action of $\pi_1(M)$. On the other hand, for $L$ the leaf space of a
{\em minimal} foliation with one--sided branching, we will see that incompressibility
of a subinterval puts very strong restrictions on the geometry and topology of $L$.

The incompressibility of some interval $I \subset L$ for $\F$ minimal
turns out to be a coarse version of an invariant measure.
One can make this precise: an invariant measure on
$L$ gives rise to a coclosed $1$--cocycle and therefore an element of 
$H^1(\pi_1(M),\R)$. By contrast, the existence of an incompressible interval
gives rise to an {\em approximately coclosed} $1$--cocycle whose coboundary is a
representative of an element of the {\em bounded cohomology group} $H^2_b(\pi_1(M),\R)$.

\begin{thm}\label{action_dichotomy}
Let $L$ be a simply--connected non--Hausdorff one--manifold
with one--sided branching in the negative direction; that is, with
respect to an orientation on $L$ and hence a partial order, any two
elements of $L$ have a common upper bound. Suppose
$G$ acts on $L$ in an orientation preserving way with dense
orbits. Then one of the following two possibilities occurs:
\begin{itemize}
\item{For any two embedded intervals $I,J \subset L$ there is an
$\alpha \in G$ such that $\alpha(I) \subset J$.}
\item{There is a $G$--equivariant map $i:L \to \R$ whose restriction
to any embedded interval $I \subset L$ is injective and order--preserving,
and the induced action of $G$ on $\R$ is conjugate to a subgroup of 
$\til{\hom(S^1)}$.}
\end{itemize}
\end{thm}
\begin{pf}
For each embedded interval $I \subset L$ we use the notation
$I^\pm$ for the highest and lowest points of $I$ respectively. Note that
the only embedded intervals in $L$ are totally ordered.
Suppose there exists an incompressible interval $I$ --- i.e. one for which no
$\alpha \in G$ takes $I$ properly inside itself. Then of course, no $\alpha$
takes $I$ properly outside itself either.
Every $p \in L$ is contained in some translate of $I$, by minimality.
For each $p \in L$ we define
$$Z(p) = \sup_{\alpha(I^-)\le p} \alpha(I^+)$$

Note that the supremum of a collection of points in $L$ is well--defined if
they have a common upper bound, since $L$ branches only in the negative 
direction. Moreover, all the $\alpha(I^+)$ are bounded by $\beta(I^+)$ for
$\beta$ any element with $\beta(I^-)>p$, so $Z(p)$ exists.

We enumerate some properties of $Z$:
\begin{enumerate}
\item{If there exists $\alpha$ so that $p = \alpha(I^-)$ then
$Z(p) = \alpha(I^+)$.}
\item{The map $Z$ is {\em order--refining}; i.e. if $x>y$ then $Z(x)>Z(y)$. For, by
minimality there is some $x>z_1>z_2>y$ with $z_1 = \alpha_1(I^-)$ and $z_2 = \alpha_2(I^-)$,
and by incompressibility of $I$, $Z(x)\ge \alpha_1(I^+)>\alpha_1(I^-)\ge Z(y)$.}
\item{Similarly, $Z(y)>y$ for any $y$, since by minimality there is some $\alpha$ with
$\alpha(y)$ contained in the interior of $I$.}
\item{For any $y$, the sequence $Z^n(y)$ escapes in the positive direction in $L$. For,
if $x>y$ is arbitrary, the interval $[x,y]$ can be covered by a finite number $m$ of translates
of $I$, by minimality. It follows that $Z^m(y)\ge x$ and the claim follows.}
\end{enumerate}

By naturality, $Z$ commutes with the
action of $G$ on $L$. We can define a total order on a quotient of
$L$ by declaring $x \le y$ iff $Z^n(x) \le Z^n(y)$ for some sufficiently large $n$.
This defines a projection to a totally ordered set $L'$, which is injective
on every totally ordered subset of $L$. We can 
complete $L'$ with respect to the order topology to get a closed subset $L''$
homeomorphic to $\R$, and therefore obtain a $G$--equivariant
map $i:L \to \R$. The projection of $Z$ by $i_*$ is an increasing
map $\R \to \R$ without fixed points. Since $Z$ is continuous, its projection by
$i_*$ is also continuous. Hence $Z$ projects to a
translation of $\R$ commuting with the action of every element of $G$ on
$\R$; in particular, this action is conjugate to a subgroup of
$\til{\hom(S^1)}$.

Notice that each ordered interval $I \subset L$ is mapped injectively
into $\R$ by $i$. For, if $r,s \in I$ are any two points, there are
uncountably many points between the image of $r$ and $s$ in $L'$ and
therefore also in $\R$.

If by contrast for {\em every} interval $I$ there is some $\alpha$ taking
$I$ properly inside itself, we can find a sequence of elements $\alpha_i$
so that for each $i$, $\alpha_{i+1}(I) \subset \alpha_i(I)$
and furthermore, $\cap_i \alpha_i(I) = p$ for some single point $p$.
If $J$ is any other interval,
minimality implies that there is a $\beta \in G$ with $\beta(p)$ in the
interior of $J$. Then $\beta\alpha_i(I) = J_i$ is an open interval containing
$\beta(p)$, such that $\cap_i J_i = \beta(p)$. Fix an $i$, and let $K = J \cap J_i$.
Since $L$ is a $1$--manifold, $K$ is an open interval containing $\beta(p)$ in its
interior. Then $J_j \cap J = J_j \cap K$ for $j>i$, and this is a nested sequence of
open intervals whose intersection is $\beta(p)$. In particular, for some sufficiently
large $j$, $\beta\alpha_j(I) = J_j \subset K \subset J$.
\end{pf}

\begin{rmk}
A useful way to think of the action of the monotone map $Z$ on $L$ is
like zipping up a zipper: two intervals $I,J$ with common highest point
but incomparable lowest point are ``zipped'' into a single interval
by comparing $Z(I)$ and $Z(J)$. The ``teeth'' of the zipper are the
points in $I$ and $J$ which are interspersed together in some unique
way determined by the topology of the action of $G$.
\end{rmk}

\begin{rmk}
Associated to a group action of $G$ on $\R$ which is conjugate
into $\til{\hom(S^1)}$,
there is a 2--cocycle $c$ defined as follows. Define a $1$--cocycle $s$
whose value on $(g,h)$ with $g,h \in G$
is the rotation number of $g^{-1}h$. Then $c = \delta s$. It is easy to
see that $\|c\|_\infty <\infty$; that is, $c$ is a {\em bounded} 2--cocycle
and determines, except in very elementary cases,
a nontrivial element of $H^2_b(G,\R)$, as suggested earlier. 
\end{rmk}

If a foliation $\F$ has associated to it an action of $\pi_1(M)$ of the 
kind described by the second case we say that $\F$ is obtained from a 
{\em one--sided branched slithering} of $M$ over $S^1$.

\begin{exa}\label{bs_example}
Let $C$ be a $2$--complex obtained from a cylinder by gluing one end to
the other by an $n$--fold cover. Then 
$$\pi_1(C) = \langle x,t \; | \; txt^{-1} = x^n \rangle$$ 
is a Baumslag--Solitar group. $C$ has a foliation by circles which
lifts to a foliation of $\til{C}$ by ``horizontal'' lines. The leaf space
of this foliation has one--sided branching; moreover there is a natural
map of $\til{C}$ to $\R$ (thought of as the universal cover of a circle 
representing the generator $t$) on which $\pi_1(C)$ acts by translations.
In this example, the action on the leaf space is not dense. We can fix this
in the following way: glue an annulus $A$ onto the $2$--complex $C$ in such
a way that both ends of the annulus are identified with a loop representing
the generator $t$. This gives a new $2$--complex $D$ with
$$\pi_1(D) = \langle x,y,t \; | \; txt^{-1} = x^n, [t,y] = \text{id} \rangle$$
The annulus $A$ can be foliated by intervals in such a
way that holonomy around the annulus induces an arbitrary homeomorphism
of the circle representing $t$. The leaf space of the universal cover of
this foliation still has one--sided branching, but now every leaf of the
foliation of $D$ is dense.

By appropriately thickening this two complex, we can realize it as a foliation
of a $4$--manifold with boundary by $3$--dimensional leaves which are all
infinite genus handlebodies. It is easy to see how to find a manifold
neighborhood of this $2$--complex in the cylinder piece, but slightly more
tricky where the two end circles are glued. Imagine a neighborhood of this
circle being obtained from two pieces: one which is just a product 
$S^1 \times I$ and one consisting of a singular fiber of a Seifert fibration
together with a spiral of nonsingular fibers converging to the singular
fiber. Once these are thickened, let $A_1,A_2$ be two annuli in the
boundary of this $4$--manifold transverse to the foliation; the annulus
$A$ can be added by gluing $A_1$ to $A_2$ with a ``twist'' in the circle
direction.
\end{exa}

This example cannot arise in the context of $3$--manifolds. Firstly, because
of the results which we have already proved in section 2.2, but
at a more basic level, Baumslag--Solitar groups cannot embed in the fundamental
groups of $3$--manifolds, by a theorem of Shalen (\cite{Shalen_BS}).

If $\F$ does not have every leaf dense, 
then we can identify a minimal set $L'$ for the leaf
space $L$, and quotient out the space $L$ by contracting all complementary
regions to points. The quotient leaf space is still a simply
connected non--Hausdorff $1$--manifold with one--sided branching and an
action of $\pi_1(M)$, except now every leaf is dense. 

We now show that taut foliations of $3$--manifolds cannot arise from a
one--sided branched slithering over $S^1$.

\begin{thm}\label{no_branched_slithering}
Let $\F$ be a taut foliation of a $3$--manifold 
with one--sided branching. Then $\F$ does not arise
from a branched slithering over $S^1$.
\end{thm}
\begin{pf}
Let $\phi:L \to \R$ be the equivariant order--preserving
map guaranteed by theorem~\ref{action_dichotomy}. Let 
$I \subset \R$ be one unit of the map $Z$, and let $J\subset L$ correspond to
a maximal connected pocket of leaves of $\til{\F}$ which projects onto $I$.
By abuse of notation, we also denote by $J$ the corresponding saturated
subset of $\til{M}$.

Then $J$ has an uppermost leaf $\lambda$ and is bounded below 
by some collection of lowermost leaves $\mu_i$. Every point 
$p \in \til{M}$ can be joined by some transversal to a point
one unit of the branched slithering below it. Since the definition
of this unit only depends on the projection of $p$ to $M$, there is
a finite supremum $\delta$ on the length of a shortest such transversal.
Hence in particular, 
$$\lambda \subset N_\delta\Bigl(\bigcup_i \mu_i\Bigr)$$ 
for some $\delta$. 

Now pick some $p \in \mu_i$ within distance $\delta$ from $\lambda$. 
We know by lemma~\ref{big_ball_dense} that for any $\epsilon$ there 
is a $t$ such that any ball of radius $t$ in a leaf of $\Lambda$ 
is an $\epsilon$--net in $\Lambda$. Suppose for some $\mu_i$ that 
$\mu_i$ is not contained in a bounded neighborhood of $\lambda$. Then
we can find $q \in \mu_i$ a very long distance from $\lambda$, and an $\alpha$
such that $\alpha(p)$ is within $\epsilon$ of $r$ on the positive side, where 
$r \in \mu_i$ is within distance $t$ of $q$ in $\mu_i$. We must have 
$\alpha(\lambda) > \lambda$, by incompressibility of $J$. But now $q$ 
is within distance $t + \delta + \epsilon$ of $\lambda$, contrary to assumption. 

It follows that $\mu_i$ is contained in a bounded neighborhood of 
$\lambda$, where the bound is independent of $i$. By 
lemma~\ref{neighborhood_is_product}, this 
implies that the space between $\lambda$ and $\bigcup_i \mu_i$ is 
foliated as a product; in particular, there is only one lowermost leaf
$\mu_i = \mu$. 
By the definition of $J$, the translates of $J$ cover
$\til{M}$. It follows that $\til{\F}$ is foliated as a product, and
$\F$ is $\R$--covered, contrary to assumption. 
\end{pf}

This leads to the following corollary:

\begin{cor}\label{compressible_action}
Let $\F$ be a taut minimal foliation of a $3$--manifold
with one--sided branching. Then for any two embedded intervals $I,J$ 
there is an $\alpha$ in $\pi_1(M)$ such that $J \subset \alpha(I)$.
\end{cor}
\begin{pf}
Apply theorem~\ref{action_dichotomy} and theorem~\ref{no_branched_slithering}.
\end{pf}

\begin{rmk}
The proof of theorem~\ref{no_branched_slithering} uses the full power of
lemma~\ref{neighborhood_is_product}, and not merely the special case whose
proof is sketched in remark~\ref{short_special_case}.
\end{rmk}

\section{The pinching lamination}

Throughout this section we will assume that all our foliations
$\F$ are taut minimal foliations with one--sided branching, by appealing
to theorem~\ref{blow_down_pockets}.

\subsection{Weakly confined directions}

\begin{defn}
Define $E_\infty$ to be the circle bundle over $L$
whose fiber over a leaf $\lambda$ is $S^1_\infty(\lambda)$.
For each $v \in UT_p\lambda$, there is a unique $e(v) \in E_\infty$ which
is the endpoint of the geodesic ray in $\lambda$ which starts out at $p$
tangent to $v$. Call this the {\em endpoint map}.
Topologize $E_\infty$ by declaring that the natural embedding
by the endpoint map of $UT\F|_\tau$ for a transversal
$\tau$ is a homeomorphism onto its image.
\end{defn}

\begin{lem}
The topology on $E_\infty$ is well--defined. That is, if $\tau_1,\tau_2$ are
two transversals to $\til{\F}$ which project to the same interval $I \subset L$,
the map $e^{-1} e: UT\F|_{\tau_1} \to UT\F|_{\tau_2}$ determined by the endpoint map
is a homeomorphism.
\end{lem}
\begin{pf}
If $p,q \in \H^2$ are any two points, the map between the visual circles of
$p$ and $q$ is determined by the {\em local} geometry of the configuration of
the pair of points. That is, it is given by parallel transport of the canonical
flat $PSL(2,\R)$ connection on $UT\H^2$ along the geodesic arc from $p$ to $q$.
If $\tau_1,\tau_2$ are a pair of transversals as in the hypotheses
of the lemma, then leafwise parallel transport of the leafwise flat $PSL(2,\R)$ connection
on $UT\F$ along leafwise geodesic arcs is continuous and defines a homeomorphism
between $UT\F|_{\tau_1}$ and $UT\F|_{\tau_2}$.
\end{pf}

The ``trick'' to this lemma is the observation that visual circles can be compared
by using {\em local} data; we do not need to compare the asymptotic geometry of rays
emanating from points on $\tau_1,\tau_2$, which is just as well since we have no
control over how these rays vary in the large as subsets of $\til{M}$.

The space $E_\infty$ carries a natural foliation by circles, the circles
at infinity of leaves of $L$. We will tame the topology of
$E_\infty$ by finding certain arcs transverse to this foliation, which detect
asymptotic features of the geometry of leaves of $\til{\F}$.

\begin{defn}
Let $\F$ be a taut foliation with one--sided branching in the negative
direction, and let $\lambda$ be a leaf of $\til{\F}$. A point
$p \in S^1_\infty(\lambda)$ is {\em weakly confined below} if there is
a map $H:I \times \R^+ \to \til{M}$ with the following properties:
\begin{itemize}
\item{Each $H(I \times \text{point})$ is a transversal to $\til{\F}$.}
\item{There is a $k\ge 1$ such that each 
$H(\text{point} \times \R^+)$ is a $k$--quasigeodesic in some
leaf of $\til{\F}$. Moreover, $H(1 \times \R^+) \subset \lambda$. We refer to these
as {\em horizontal rays} in the marker.}
\item{There is an $\epsilon$ such that each transversal $H(I \times \text{point})$
has diameter bounded above by $\epsilon$. We refer to these as {\em vertical intervals}
in the marker.}
\end{itemize}
The map $H$ is called a {\em marker}, or a {\em $(k,\epsilon)$--marker} if we wish to
specify $k$ and $\epsilon$. We say that the marker {\em certifies} the point $p \in
S^1_\infty(\lambda)$.
\end{defn}

We will typically be interested in $(k,\epsilon)$--markers for which $\epsilon$
is small enough so that leaves of $\til{\F}$ are quasi--geodesically embedded in
their $2\epsilon$--neighborhoods.

\begin{lem}\label{straighten_markers}
A marker can be straightened so that every $H(\text{point} \times \R^+)$ is
geodesic in its leaf.
\end{lem}
\begin{pf}
For a quasigeodesic line, ray or segment $\alpha$ in $\H^2$, let $\alpha^s$ denote
the geodesic segment with the same (possibly ideal) endpoints as $\alpha$.
If $\gamma$ is a $k$--quasigeodesic ray in $\H^2$, then for any $\epsilon$ there is
a $t$ such that if $\sigma$ is a segment of $\gamma$ of length $\ge t$, and $\sigma'$ is
the middle third segment of $\sigma$, then $\sigma^s$ and $\gamma^s$ are $\epsilon$--close
along the segments corresponding to $\sigma'$. More generally, if an interval
$\sigma$ has an endpoint in common with $\gamma$, the geodesic interval $\sigma^s$ is
$\epsilon$--close to $\gamma^s$ away from the last third segment of $\sigma^s$.
Notice that we are always comparing the geodesic straightenings of a quasigeodesic
ray and a {\em finite} segment along the middle third, or initial two thirds of the
finite segment, so this always makes sense.

If $H$ is a marker, let $H^s$ denote the set of leafwise rays which are obtained
by leafwise straightening the rays $H(\text{point} \times \R^+)$.
For each $p \in I$, let $\gamma_p$ be the quasigeodesic ray $H(p \times \R^+)$.
If $\sigma_p,\sigma_q$ are segments of $\gamma_p,\gamma_q$ which are sufficiently
close in the Hausdorff topology, the leafwise straightenings $\sigma_p^s,\sigma_q^s$ can
also be made arbitrarily close. If $p_i \to p$, we can find initial segments of
$\gamma_{p_i}$ with length $\to \infty$ which converge to $\gamma_p$ and whose
leafwise straightenings $\gamma_{p_i}^s$ are eventually
$\epsilon$--close to $\gamma_p^s$ and to $\gamma_{p_i}^s$ away from their last
thirds. We reiterate that the $\gamma_{p_i}$ are always finite segments, so it
makes sense to say that the straightened initial two thirds of this segment stay close
to the straightening of $\gamma_p$.
It follows that the family $\gamma_p^s$ varies continuously as a function of
$s$ in the compact--open topology. Points in the markers are moved a uniformly bounded
distance by this leafwise straightening, so the diameters of the transversals are uniformly
bounded.
\end{pf}

The slight subtlety with this lemma is that we cannot assume that the leafwise rays
$\gamma_p$ vary continuously in the Hausdorff topology. It is certainly possible for
there to exist $\delta>0$ and $p_i \to p$ such that $\gamma_{p_i}$ is not contained in
the $\delta$--neighborhood of $\gamma_p$ for any $i$.

From now on we will only consider markers which are leafwise geodesic.
It follows that if $H:I \times \R^+ \to \til{M}$ is a marker, there is an associated
map $e(H):I \to E_\infty$ which is given by the leafwise endpoint of the geodesic
ray in the image of the marker. By abuse of notation, we will also refer to these
transversals at infinity as markers.

\begin{thm}\label{markers_exist}
For $\F$ minimal, for every leaf $\lambda \in \til{\F}$ there is some 
point which is weakly--confined below in $S^1_\infty(\lambda)$. Moreover, the 
certifying marker can be taken to be a $(1,\epsilon)$--marker for any $\epsilon>0$.
\end{thm}
\begin{pf}
Let $\gamma$ be a nontrivial curve contained in a leaf $\lambda$ of
$\F$. Such a $\gamma$ exists, since if every leaf of $\F$
is a plane, the entire foliation is by planes and $M$ is $T^3$ foliated
by geodesic planes of maximally irrational slope, by a theorem of
Rosenberg (\cite{Rosenberg}). Now for either side of
$\gamma$, there is a choice of orientation on $\gamma$ such that
holonomy through $\gamma$ is weakly contracting on that side. Here we
say that holonomy through $\gamma$ is 
{\em weakly contracting} (resp. {\em weakly expanding}) on the negative side of
$\lambda$, say, if for some transversal $\tau$ to $\F$ intersecting $\lambda$
at $p \in \gamma$, there are a sequence of points $p_i \to p$ with
$p_{i+1} > p_i$ such that if $H:\tau \to \tau$ denotes the holonomy around
$\gamma$, where it is defined, $H(p_i) \ge p_i$ (resp. $H(p_i) \le p_i$)
for all sufficiently large $i$.
Clearly, holonomy around any loop is either weakly contracting or
weakly expanding on either side, and by possibly reversing the orientation
of $\gamma$, we can take it to be weakly contracting. Let $\tau_i$ denote
the segment between $p$ and $p_i$ where $i$ is sufficiently large as above.
Lift to $\til{M}$, and let $R:I \times \R^+ \to \til{M}$ be defined by
holonomy transport of $\tau_i$ in the positive direction along the lift of
$\gamma$. Then $R$ is a marker. By restricting to a sufficiently small
segment of $\tau$, we can arrange for $R$ to be a $(k,\epsilon)$--marker
for any $\epsilon$. By lemma~\ref{straighten_markers}, $R$ can be taken
to be a $(1,\epsilon)$--marker for some better initial choice of $\epsilon$.

Since $\F$ is minimal, it follows that there is a weakly--confined direction
below in the circle at infinity of an arbitrary leaf.
\end{pf}

\begin{rmk}
An alternative approach to this theorem uses harmonic transverse measures
for foliations, as constructed in \cite{lG83}. This is Thurston's approach
in \cite{wT98}
\end{rmk}

The most important feature of markers is that they {\em do not coalesce}.
That is, if $H_1,H_2$ are markers which intersect the same interval
of leaves in $L$, if their endpoints are disjoint in some leaf, they are
disjoint in all leaves. We state this as a lemma:

\begin{lem}\label{markers_disjoint}
Let $H_1,H_2:I \times \R^+ \to \til{M}$ be $(1,\epsilon)$--markers. 
Then either the endpoint markers $e(H_1), e(H_2)$ agree up to the
range of their image, or they are disjoint. Furthermore, if the
$e(H_1),e(H_2)$ intersect incomparable leaves, they are disjoint.
Here $\epsilon$ is a sufficiently small positive number, which also has
the property that leaves of $\til{\F}$ are quasigeodesically
embedded in their $2\epsilon$--neighborhoods.
\end{lem}
\begin{pf}
Suppose $e(H_1)=p_1$ and $e(H_2)=p_2$ are distinct points in
$S^1_\infty(\lambda)$. Since every leaf of $\til{\F}$ is uniformly properly
embedded in its $\epsilon$ neighborhood for {\em some} $\epsilon$, there is a compact subset
$K \subset \til{M}^\lambda$ such that the geodesic rays
$H_1(1,\R^+)$ and $H_2(1,\R^+)$ do not
come within distance $t$ of each other in $\til{M}$ outside $K$, for {\em some} $t$.
Notice that $t$ here cannot be arbitrary.

On the other hand, if $e(H_1) \cap e(H_2)$ is non empty, there is some 
$q \in S^1_\infty(\nu)$ for some leaf $\nu$
in the image of both. It follows that the leafwise
geodesic rays in $\nu$ in the image of $H_1$ and $H_2$ are both asymptotic to $q$,
and therefore they come arbitrarily close, contrary to the definition of $K$.
it follows that if $e(H_1),e(H_2)$ are disjoint in some leaf
they are disjoint in all other leaves.

Finally, if $e(H_1)$ and $e(H_2)$ intersect incomparable leaves 
$\mu_1,\mu_2 < \lambda$, they are disjoint in all comparable leaves.
For, if they intersect in $S^1_\infty(\lambda)$, the ends of $H_1,H_2$ are
contained within an arbitrarily small neighborhood of each other. i.e. there
exist sequences of pairs of points $p_1^i,p_2^i$ contained in $H_1,H_2$ respectively
which exit every compact set, and for which the distance between $p_1^i,p_2^i$ in
$\til{M}$ converges to $0$.
On the other hand, for sufficiently large $i$, the vertical intervals in $H_1,H_2$ have
length $\le \epsilon$, by the definition of markers. In particular, the points
$p_1^i,p_2^i$ are within $\epsilon$ of points $q_1^i,q_2^i$ which are contained in
incomparable leaves. But our choice of $\epsilon$ was such that leaves are quasi--geodesically
embedded in their $2\epsilon$ neighborhoods. In particular, this
implies that incomparable leaves never come within $2\epsilon$ of each other, which
gives us a contradiction. This contradiction finishes the proof of the lemma.
\end{pf}

\subsection{Ergodic theory at infinity}

The following material is taken more or less directly from \cite{dC99b}.

\begin{defn}
Let $\F$ be a foliation of the compact manifold $M$. A probability measure $\mu$
on $M$ is {\em harmonic} if
$$\int_{M} \Delta_\F f d\mu = 0$$
for every function $f$ which is measurable on $M$ and smooth in the leafwise direction.
here $\Delta_\F$ denotes the {\em leafwise Laplacian}, which governs the heat flow in
the manifold where heat is restricted to flow along leaves of the foliation.
\end{defn}

Garnett's theorem from \cite{lG83} says that nontrivial harmonic measures exist
for foliations:

\begin{thm}[Garnett]
If $\F$ is an arbitrary foliation of a compact space,
then it admits a non--zero harmonic measure. If $\F$ is
minimal, this measure decomposes locally into the
product of leafwise Riemannian measure with a transverse
measure of full support which is infinitesimally harmonic, as a function
of the leaf metric.
\end{thm}  

The following theorem is proved in \cite{dC99b}:

\begin{thm}\label{circle_alternative}
Let $\F$ be a minimal foliation with hyperbolic leaves of a
compact $3$--manifold $M$. Let $\Sigma \subset E_\infty$ be a closed
$\pi_1(M)$--invariant subset. Then either the subset $\Sigma \cap S^1_\infty(\lambda)$
is dense for almost every fiber $S^1_\infty(\lambda)$, or the subset
$\Sigma \cap S^1_\infty(\lambda)$ is empty or consists of a single point for almost every
fiber $S^1_\infty(\lambda) \subset E_\infty$.
\end{thm}

The proof is a straightforward application of Garnett's theorem, and
the only hypotheses are that $\F$ is minimal with hyperbolic leaves, and
$M$ is compact. In particular, the kind of branching of $\til{\F}$ is
irrelevant (though of course it might be relevant to prove minimality).

\subsection{The pinching lamination}

We will show that markers between incomparable leaves cannot ``link at infinity''.

\begin{lem}\label{markers_unlinked}
Suppose $\lambda > \mu,\nu$ with $\mu,\nu$ incomparable. If $p,q \in S^1_\infty(\lambda)$
are endpoints of markers from $\mu$ to $\lambda$ and $r,s \in S^1_\infty(\lambda)$ are
endpoints of markers from $\nu$ to $\lambda$ then the pair $p,q$ does not
link the pair $r,s$.
\end{lem}
\begin{pf}
Suppose $\lambda > \mu,\nu$ with $\mu,\nu$ incomparable.
Say $R_1,R_2$ are a pair of markers between $\lambda$ and $\mu$ and $S_1,S_2$
are a pair of markers between $\lambda$ and $\nu$. Then these markers can be
extended to bi--infinite rectangles $R:I \times \R \to \til{M}$ and
$S:I \times \R \to \til{M}$ which have ends $R_1,R_2$ and $S_1,S_2$
respectively. If the endpoints of $R$ in $S^1_\infty(\lambda)$
link the endpoints of $S$ in $S^1_\infty(\lambda)$, the geodesics $R_\lambda$
and $S_\lambda$ in $\lambda$ must cross somewhere.
In particular, $S \cap R$ is nonempty.
These subsets of $\til{M}$ are both closed, so their intersection is
closed. Either the intersection stays in a compact subset of 
$R \cup S$ or else the intersection goes off to infinity.

Since $R_\lambda$ and $S_\lambda$ intersect transversely, for every $t$
there is a compact subspace $K_t \subset \lambda$ so that points in
$R_\lambda \cap (\lambda \backslash K_t)$ do not come within distance $t$
in $\lambda$ of points in $S_\lambda \cap (\lambda \backslash K_t)$.
Since $\lambda$ is uniformly properly embedded in $\til{M}^\lambda$,
the geodesics $R_\lambda$ and $S_\lambda$ come close to each other in
$\til{M}^\lambda$ only near a compact subset.
On the other hand, both $R$ and $S$ are contained in a bounded neighborhood
in $\til{M}^\lambda$ of $R_\lambda$ and $S_\lambda$ respectively. It
follows that their intersection is compact. The intersection is just
an arc transverse to $\til{\F}$, and therefore its endpoint 
is on the bottom of either $R$ or $S$ or both.
In any case, this says that $\lambda$ and $\nu$ are comparable, contrary to
assumption.
\end{pf}

For $\lambda > \mu$, let $M(\lambda,\mu) \subset S^1_\infty(\lambda)$ denote
the set of endpoints of markers from $\mu$ to $\lambda$. We have shown that
for $\lambda > \mu,\nu$ incomparable, the sets $M(\lambda,\mu)$ and
$M(\lambda,\nu)$ are unlinked.

\begin{defn}
For $r \subset L$ a properly embedded copy of $\R$, let $C_\infty(r) = E_\infty|_r$
denote the restriction of the circle bundle $E_\infty$. $C_\infty(r)$ is topologically
a cylinder. A {\em long marker} is a map $\tau:\R \to E_\infty$ with the
following properties:
\begin{enumerate}
\item{The image of $\tau$ is transverse to the circles of $E_\infty$, and the
projection of $\tau$ to $L$ is an embedding onto a properly embedded copy $r \subset L$
of $\R$.}
\item{The transversal $\tau(\R) \subset C_\infty(r)$ is the limit in the
Hausdorff topology of a sequence of $(1,\epsilon)$--markers $\tau_i \subset C_\infty(r)$.}
\end{enumerate}
\end{defn}

If $L' \subset L$ is a submanifold (possibly with boundary),
we can also define a {\em long marker relative to $L'$}. For $r \subset L'$
which is a properly embedded copy of $\R$ (or, if $L'$ has boundary, $r$ could be
a ray or interval with boundary contained in $\partial L'$). A long marker relative
to $L'$ is a map $\tau:r \to C_\infty(r)$, transverse to the foliation by circles,
such that the composition with the projection $C_\infty(r) \to r$ is the identity,
and which is a limit in the Hausdorff topology of a sequence of markers, as above.

Let $\lambda$ be a leaf of $L$, and $L_\lambda$ the subset of $L$ corresponding
to the leaves in $\til{M}_\lambda$. Then the restriction of a long marker to
$\til{M}_\lambda$ is a long marker relative to $L_\lambda$.

To say that $\tau(\R)$ is a limit of the $\tau_i$ is the same as saying that
for any segment $\sigma \subset \tau(\R)$ and any neighborhood $N(\sigma)$ of
$\sigma$ in $C_\infty(r)$, there is an
$i$ such that for all $j>i$, the marker $\tau_j$ intersects the same leaves of $L$
that $\sigma$ intersects, and furthermore its intersection with the subset of
$C_\infty(r)$ lying over the projection of $\sigma$ to $L$ is contained in $N(\sigma)$.

\begin{defn}
Let $\lambda \in L$ be some leaf. A {\em filter below $\lambda$} is a {\em maximal} union
of pairwise incomparable leaves $\mu_i$ with each $\mu_i<\lambda$.
\end{defn}

Imagine $\lambda$ as the source of a forking river, which flows in the negative
direction in $L$ in all possible ways. Then a filter is a collection of places one
would need to dam the river to completely stop it, in such a way that every dam
blocks some water. For instance, if $L$ is homeomorphic to $\R$, any $\mu < \lambda$
by itself is a filter. 

\begin{thm}\label{dense_markers}
Suppose that $\F$ has one--sided branching in the negative direction. 
Then for every $\lambda > \mu$, the set of long markers between
$\lambda$ and $\mu$ is dense in $S^1_\infty(\mu)$.
Moreover, if $\mu_i$ is a filter below $\lambda$, the set of long markers
between $\lambda$ and some $\mu_i$ is dense in $S^1_\infty(\lambda)$.
\end{thm}
\begin{pf}
For the sake of brevity, for the duration of this proof we will refer
to both the transversals $e(H)$ to $E_\infty$ and the maps $H:I \times \R^+ \to \til{M}$
as markers. It should be clear from context which we mean in each case.

By theorem~\ref{markers_exist} we know that the set of 
markers intersects every leaf of $\til{\F}$. By
theorem~\ref{circle_alternative}, this set is either dense or intersects
almost every circle $S^1_\infty(\lambda)$ in a single point.
Let $\tau \subset E_\infty$ be a marker, and let $I \subset L$ be the
interval of leaves in $L$ which it intersects.
We know by corollary~\ref{compressible_action} that for any $J \subset L$
there is some $\alpha \in \pi_1(M)$ with $\alpha(J) \subset I$.
It follows that $\alpha^{-1}(\tau)$ is a marker which projects to
an interval in $L$ containing $J$. On the other hand, $J$ was arbitrary.
Let $J_1,J_2$ be two intervals in $L$ with a common uppermost leaf
$\lambda$ and incomparable lowermost leaves $\mu_1,\mu_2$. By construction,
there are markers $\tau_1,\tau_2$ in $E_\infty$ from $S^1_\infty(\lambda)$
to $S^1_\infty(\mu_1),S^1_\infty(\mu_2)$ respectively. By
lemma~\ref{markers_disjoint}, these markers are disjoint in $E_\infty$.
It follows that there is an interval of leaves of $L$ which intersect
at least two $(1,\epsilon)$--markers, and therefore {\em every} leaf
intersects at least two $(1,\epsilon)$--markers, and the set of
$(1,\epsilon)$--markers is dense in $E_\infty$.

Even stronger, this set of markers is dense in $S^1_\infty(\lambda)$ for
every $\lambda$. For, let $K$ be a compact fundamental domain for $M$.
For every point $p \in K$ in a leaf $\lambda$ of $\til{\F}$, there are at 
least two markers in $S^1_\infty(\lambda)$. So there is an upper bound
$\theta < 2\pi$ on the visual angle between any two markers, as seen from any
point $p \in M$. If there were some gap $I \subset S^1_\infty(\lambda)$
for some leaf $\lambda$ which did not intersect any marker, then there
would be a sequence of points $p_i \to p$ in $\lambda$ for which the
apparent visual size of this gap would converge to $2\pi$, contrary to
the existence of $\theta$. It follows that the set of markers is dense in
each $S^1_\infty(\lambda)$.

In particular, for any $n$
there are a pair of leaves $\lambda > \mu$ for which there are
at least $n$ markers from $\lambda$ to $\mu$. By corollary~\ref{compressible_action},
for {\em every} pair of leaves $\lambda > \mu$ there are {\em infinitely} many markers
from $\lambda$ to $\mu$.

Now, if $\tau \subset E_\infty$ is a marker projecting to $I \subset L$,
we can find a sequence of elements $\alpha_i \in \pi_1(M)$ such that
$\alpha_i(I) \subset \alpha_{i+1}(I)$ and $\bigcup_i \alpha_i(I) = r$, 
a properly embedded line in $L$ (i.e. $r$ is unbounded both above and
below).

Since the $\alpha_i(\tau)$ are all contained
in $E_\infty|_r$ which is locally compact, and since they are all constrained to
be disjoint from or contained in 
a dense set of transversals which are dense in every circle,
we can find a subsequence which converges in $E_\infty|_r$ to a {\em long marker}
$\tau_\infty$. We claim the translates of $\tau_\infty$ are dense in $E_\infty$.

For, if not, by theorem~\ref{circle_alternative} the set of long markers
intersects almost every leaf of $E_\infty$ in at most one point, and therefore,
since each long marker intersects every leaf above its initial point,
the set of long markers must intersect every $S^1_\infty(\lambda)$ in {\em exactly}
one point. Now, suppose $\mu_i$ are a collection of incomparable leaves below
$\lambda$. Then the long markers between $\mu_i$ and $\lambda$ must intersect
$\lambda$ in a unique point $q_i$ in $S^1_\infty(\lambda)$. Since the restriction
of these long markers to $\til{M}_\lambda$ are still long markers relative to
$L_\lambda$, they must agree, and therefore the $q_i$ are equal to some $q$, for all $i$.

But how can this be possible? Remember $\tau_\infty$, and therefore every long marker,
is a limit of a sequence of genuine $(1,\epsilon)$--markers, and $(1,\epsilon)$--markers
are {\em disjoint} (or agree on the circles they both intersect)
in $E_\infty$. So for each $j$,
$q$ is a limit of a sequence of points $q_i^j \in S^1_\infty(\lambda)$ 
with $j \to \infty$ which are the
endpoints of actual markers from $\mu_i$ to $\lambda$. If $q$ is the endpoint of an actual
marker, $q$ could be in $M(\lambda,\mu_i)$ for some $i$. But there is at most one $i$
for which this could happen. For every other $i$, infinitely many of the $q_i^j$ are 
{\em distinct}. By lemma~\ref{markers_unlinked}, we know that $M(\lambda,\mu_k)$ and
$M(\lambda,\mu_l)$ are {\em unlinked} for $l \ne k$. It follows that if $q_i^j$ contain
infinitely many points approaching $q$ from the right, say, then for $k \ne i$,
the sequence $q_k^j$ must approach $q$ from the {\em left}. But we can certainly find
at least $4$ incomparable $\mu_i$ below $\lambda$. It follows that either the $q_i$
are not all the same, or else some pair of points in $M(\lambda,\mu_i)$ links some pair of
points in $M(\lambda,\mu_j)$ for some $i \ne j$. The second contradicts
lemma~\ref{markers_unlinked}, so not all $q_i$ are the same. But this implies that the
$\tau_\infty$ are dense in $E_\infty$, and therefore dense in $S^1_\infty(\mu)$ for every
leaf $\mu$. Since every long marker continues indefinitely in the positive direction for all
time, there are a dense set of long markers in $S^1_\infty(\mu)$ from $\mu$ to $\lambda$ for
any pair $\lambda > \mu$. Since long markers can be approximated by genuine markers, there
are a dense set of genuine markers in $S^1_\infty(\mu)$ from $\mu$ to $\lambda$.

Conversely, since every long marker continues indefinitely in some negative direction for
all time, there are a dense set of long markers in $S^1_\infty(\lambda)$ from $\lambda$
to some $\mu_i$ for any filter $\mu_i < \lambda$. Again, since long markers can be
approximated by genuine markers, there are a dense set of genuine markers in
$S^1_\infty(\lambda)$ from $\lambda$ to some $\mu_i$ for any filter $\mu_i < \lambda$,
as required.
\end{pf}

Notice the asymmetry of the conclusion of this theorem. It is definitely
{\em not} true that the set of markers from $\mu$ to $\lambda$ is
dense in $S^1_\infty(\lambda)$ for $\lambda > \mu$.

In the sequel, we will assume all our genuine markers are $(1,\epsilon)$--markers without
further comment.

\begin{defn}
Suppose $\F$ has one--sided branching in the negative direction.
For $\lambda > \mu$ leaves of $\til{\F}$, define
$\Lambda^+(\lambda,\mu) \subset \lambda$ to be the
boundary of the convex hull of limit points $p$ in $S^1_\infty(\lambda)$
of {\em long markers} between $\lambda$ and $\mu$. Now set
$$\Lambda^+(\lambda) = \overline{\bigcup_{\lambda > \mu} \Lambda^+(\lambda,\mu)}$$
$$\til{\Lambda}^+ = \bigcup_\lambda \Lambda^+(\lambda)$$
We call $\Lambda^+(\lambda)$ the {\em pinching lamination} of $\lambda$.
\end{defn}

In order to justify this terminology, we will show in the sequel the following facts:
\begin{enumerate}
\item{For each leaf $\lambda$ of $\til{\F}$, the set $\Lambda^+(\lambda)$ is a
geodesic lamination. (theorem~\ref{pinching_defined})}
\item{The sets $\Lambda^+(\lambda)$ vary continuously with $\lambda$ in the
Hausdorff topology as a function of $\lambda$. They sweep out a branched
leafwise--geodesic lamination $\til{\Lambda}^+$ of $\til{M}$ which
can be equivariantly split open leafwise to an actual lamination 
which intersects leaves of $\til{\F}$ in $k,\epsilon$
quasi--geodesics, for any $k>1,\epsilon>0$, and covers a genuine lamination
$\Lambda^+$ in $M$. (lemma~\ref{produces_genuine_lamination})}
\end{enumerate}

The idea is that the lamination $\til{\Lambda}^+$ describes how the leaves of
$\til{\F}$ ``pinch'' off as we move in the negative direction. A higher
leaf $\lambda$ cobounds a region in $\til{M}$ together with
a disjoint collection of lower leaves $\coprod_j \mu_j$
which are obtained topologically from $\lambda$ by ``pinching'' $\lambda$ along the
leaves of $\Lambda^+(\lambda,\mu_j)$ for each $j$.

For $\lambda > \mu$, let $M_l(\lambda,\mu) \subset S^1_\infty(\lambda)$ denote
the set of endpoints of long markers between $\mu$ and $\lambda$. So
$\Lambda^+(\lambda,\mu)$ is the boundary of the convex hull of $M_l(\lambda,\mu)$,
for each pair $\lambda > \mu$. Notice that for each $\lambda,\mu$, the
set $M_l(\lambda,\mu)$ is closed, and is contained in the closure of $M(\lambda,\mu)$.

\begin{thm}\label{pinching_defined}
Suppose $\F$ has one--sided branching in the negative direction. Then
for any $\lambda$, $\Lambda^+(\lambda)$ is a geodesic lamination in $\lambda$.
\end{thm}
\begin{pf}
Suppose $\lambda > \mu > \nu$.
Any marker, and hence any long marker between $\lambda$ and $\nu$ actually contains
a submarker between $\lambda$ and $\mu$, so we have containment
$$M_l(\lambda,\nu) \subset M_l(\lambda,\mu)$$
and therefore $\Lambda^+(\lambda,\mu)$ and $\Lambda^+(\lambda,\nu)$ do not intersect
transversely.

On the other hand, if $\lambda> \mu,\nu$ incomparable, lemma~\ref{markers_unlinked}
implies that $M(\lambda,\mu)$ and $M(\lambda,\nu)$ are unlinked, and therefore
$M_l(\lambda,\mu)$ and $M_l(\lambda,\nu)$ are unlinked since they are contained
in the closures of the $M(\lambda,\cdot)$. It follows that
$\Lambda^+(\lambda,\mu)$ and $\Lambda^+(\lambda,\nu)$ do not intersect transversely.
\end{pf}

\begin{rmk}
It is possible to define the $\Lambda^+(\lambda,\mu)$ in terms of the endpoints of
{\em all} markers between $\lambda$ and $\mu$,
but it will be more technically convenient in the sequel to work with long markers.
\end{rmk}

\subsection{Monotone maps and laminar relations}

\begin{defn}
A geodesic lamination $\Lambda$ of $\H^2$ determines and is determined
by a relation $R_\Lambda$ on $S^1_\infty$, i.e. a subset of
the space of pairs of {\em distinct} points $S^1_\infty \times S^1_\infty - \text{diagonal}$
with the following properties:
\begin{enumerate}
\item{The relation is {\em symmetric}. i.e. $x R_\Lambda y$ iff $y R_\Lambda x$ for each
pair $x,y$.}
\item{The subspace of $S^1_\infty \times S^1_\infty - \text{diagonal}$
consisting of pairs $x,y$
with $x R_\Lambda y$ is {\em closed}}
\item{If $x R_\Lambda y$ and $z R_\Lambda w$ then the pair $x,y$ does not
link the pair $z,w$ in $S^1_\infty$.}
\end{enumerate}
We refer to such a relation as a {\em laminar relation of $S^1_\infty$}, or by
abuse of notation, as a {\em lamination of $S^1$}.
The leaves of $\Lambda$ are in $1$--$1$ correspondence with unordered pairs
$x,y$ with $x R_\Lambda y$.
\end{defn}

Notice that if $R_\Lambda$ is a relation on $S^1_\infty$ satisfying conditions (1) and
(3) above but not necessarily (2), its closure in $S^1_\infty \times S^1_\infty$ is a
laminar relation. 

Notice too that a laminar relation is {\em not} an equivalence relation. That is,
we do {\em not} require that this relation should be {\em transitive}. Firstly,
this is forced on us by the fact that we consider the relation as a subset of
the product {\em minus the diagonal}. But secondly, and more importantly,
any four distinct points in $S^1$ contains linking pairs, so transitivity would
be to a large extent incompatible with the ``no--linking'' condition.

\begin{defn}
A {\em monotone map} $\phi:S^1 \to S^1$ is a degree one map which {\em does not reverse}
the circular order on triples of points. That is, if $(x,y,z) \in (S^1)^3$
is a positively ordered triple of points, then $(\phi(x),\phi(y),\phi(z))$ is
not a negatively ordered triple of points, and vice versa. Note here that the
source and the target circles are distinct, and should not be thought of as the
same circle.
\end{defn}

A monotone map need not be injective. That is, if $(x,y,z)$ are
a positively ordered (necessarily distinct) triple, the triple
$(\phi(x),\phi(y),\phi(z))$ need not be distinct, and therefore need not be
positively ordered. A monotone map can be uniquely described as a quotient
map, where $\bigcup_i K_i \subset S^1$ is a countable union of disjoint closed intervals,
and the map $\phi:S^1 \to S^1$ is the quotient map $S^1 \to (S^1/\sim) \cong S^1$
given by quotienting each $K_i$ to a point.

If $\Lambda$ is a geodesic lamination of $\H^2$, determined by a laminar relation
$R_\Lambda$ of $S^1_\infty$, and $\phi: S^1_\infty \to S^1_\infty$ is a monotone
map, then the image of $\phi(R_\Lambda)$ in $S^1_\infty \times S^1_\infty - \text{diagonal}$
is a laminar relation of $S^1_\infty$, where
$$\phi:S^1_\infty \times S^1_\infty - \text{diagonal} \to S^1_\infty \times S^1_\infty$$
is defined by
$$\phi(x,y) = (\phi(x),\phi(y))$$
Conversely, if
$R_\Lambda$ is a laminar relation of $S^1_\infty$, we can {\em pull back}
$R_\Lambda$ by $\phi$ to produce a new laminar relation $\phi^{-1}(R_\Lambda)$ as
follows:
\begin{enumerate}
\item{If $xR_\Lambda y$, and both $\phi^{-1}(x),\phi^{-1}(y)$ consist of single
points $p,q$ then $p\phi^{-1}(R_\Lambda)q$.}
\item{If $xR_\Lambda y$ where $\phi^{-1}(x) = p$ and $\phi^{-1}(y) = I$ for some point $p$
and interval $I$, then
$p\phi^{-1}(R_\Lambda)I^+$ and $p\phi^{-1}(R_\Lambda)I^-$.}
\item{If $xR_\Lambda y$ where $\phi^{-1}(x) = I$ and $\phi^{-1}(y) = J$ for some intervals
$I,J$ then $I^+\phi^{-1}(R_\Lambda)J^-$ and $J^+\phi^{-1}(R_\Lambda)I^-$.}
\end{enumerate}

Notice here that by our convention, $xR_\Lambda y$ implies that $x$ and $y$ are
{\em distinct}.

Here $I^+$  and $I^-$ 
denote the anticlockwisemost and clockwisemost endpoints of an interval
$I \subset S^1$ respectively. Said more informally, a circle is obtained from
another circle under the inverse of a monotone map by {\em blowing up} a countable
collection of points to intervals.

\begin{figure}[ht]
\centerline{\relabelbox\small \epsfxsize 4.0truein
\epsfbox{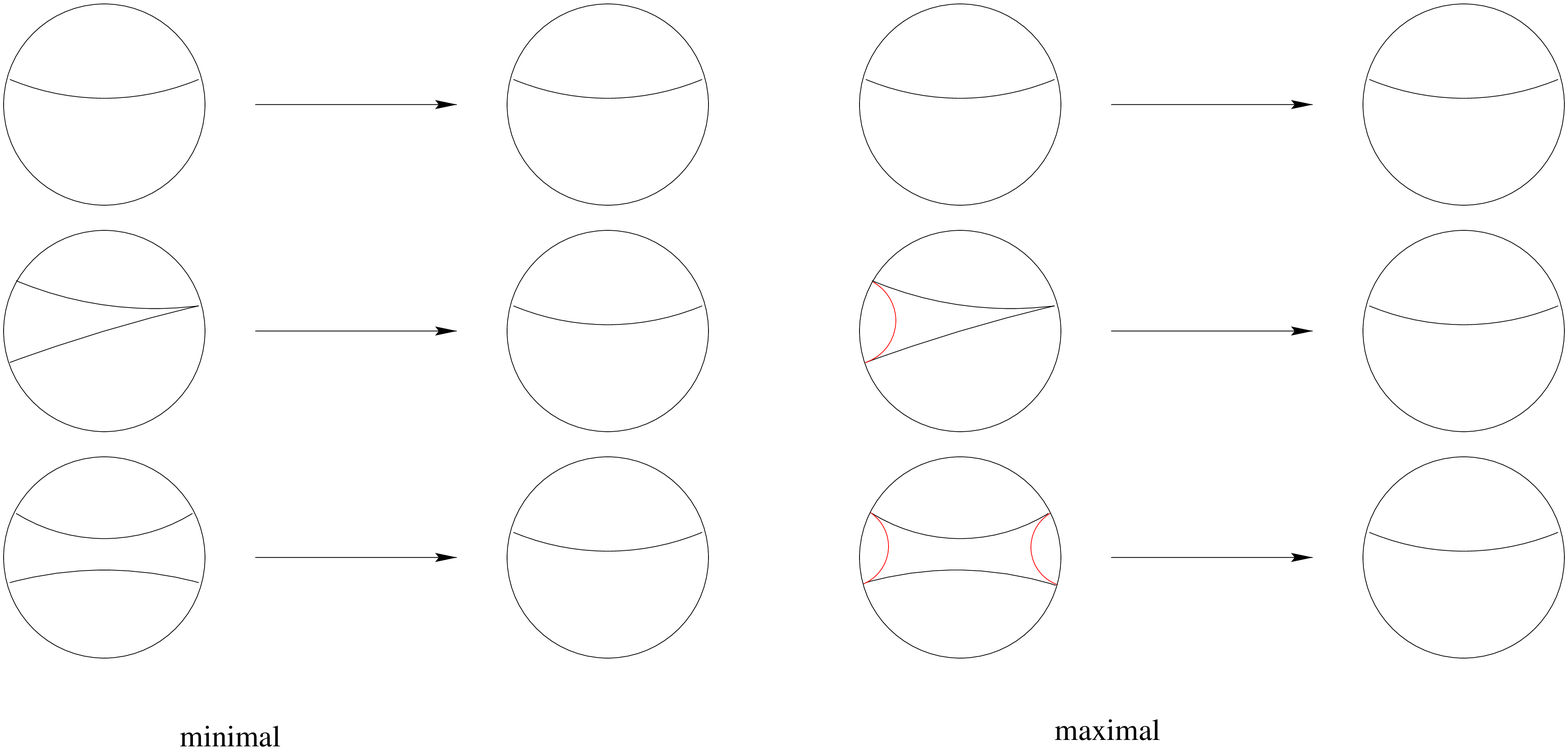}
\relabel {minimal}{minimal pullback}
\relabel {maximal}{maximal pullback}
\endrelabelbox}
\caption{A laminar relation on $S^1$ can be pulled back under
a monotone map to another laminar relation. There are at least two obvious
canonical pullbacks, which we refer to as the {\em minimal} pullback,
and the {\em maximal} pullback.
By convention, we will assume our pullbacks are {\em minimal}.}
\end{figure}

Notice that if $R_\Lambda$ is the trivial laminar relation on $S^1$ (corresponding to the
empty subset of $S^1 \times S^1 - \text{diagonal}$),
the preimage $\phi^{-1}(R_\Lambda)$ under a monotone map will also be trivial.
Another possible definition of the pullback of a lamination would add
$x \phi^{-1}(R_\Lambda) y$ for $x,y$ endpoints of a nontrivial interval $\phi^{-1}(p)$
for any point $p$.
If we want to distinguish these
pullbacks, we will refer to them as the {\em minimal pullback} and the
{\em maximal pullback} respectively. These are just names --- they do not
denote any special property of these canonical pullback laminations
with respect to the lattice of inclusions of laminar relations with image $R_\Lambda$.
In the sequel, by convention, $\phi^{-1}(R_\Lambda)$
will always refer to the minimal pullback. The difference between these two
pullbacks is illustrated in figure 3.

As we have just remarked, it is clear that the composition 
$$R_\Lambda = \phi \circ \phi^{-1}(R_\Lambda)$$
On the other hand, under the composition $\phi^{-1} \circ \phi$ some leaves of
the lamination might be lost, and one obtains in general a sublamination.
If $\phi_i:S^1 \to S^1$ are a sequence of monotone maps converging to the identity and
if $R_\Lambda$ is any laminar relation on $S^1$, then the diameters of
the intervals in $S^1$ which are collapsed by the $\phi_i$
converge to $0$, and the laminar relations converge:
$$\lim_{i \to \infty} \phi_i^{-1} \circ \phi_i(R_\Lambda) = R_\Lambda$$

It is also clear that if $\phi_1:S^1_0 \to S^1_1$ is a monotone map,
and $\phi_2:S^1_1 \to S^1_2$ is another monotone map, for any laminar
relation $R_\Lambda$ of $S^1_0$, the composition is associative:
$$(\phi_2 \phi_1)(R_\Lambda) = \phi_2(\phi_1(R_\Lambda))$$

The long markers between leaves of $\til{\F}$ give rise to a coherent family of
monotone maps between the circles $S^1_\infty(\lambda)$ which relate the laminations
$\Lambda^+(\lambda)$. For the sake of brevity, in the sequel we will use the
same notation for a geodesic lamination of a hyperbolic plane and for the associated
laminar relation on its ideal boundary. 

\begin{thm}\label{universal_circle}
Suppose $\F$ has one--sided branching in the negative direction. Then there
is a universal circle $S^1_\u$ with a natural action of $\pi_1(M)$ on it.
The action of $\pi_1(M)$ leaves invariant a laminar relation 
$\Lambda^+_\u$ of $S^1_\u$. Moreover, 
for each leaf $\lambda$ of $\til{\F}$ there is a monotone map 
$$\phi_\lambda:S^1_\u \to S^1_\infty(\lambda)$$
such that $\phi_\lambda(\Lambda^+_\u) = \Lambda^+(\lambda)$, and such that
for $\lambda > \mu$, there is a monotone map
$\pi_\mu^\lambda:S^1_\infty(\lambda) \to S^1_\infty(\mu)$ with the property
$$\phi_\mu = \pi_\mu^\lambda \circ \phi_\lambda$$
\end{thm}
\begin{pf}
For $\lambda > \mu$, we know the set of long markers between $\lambda$ and
$\mu$ are dense in $S^1_\infty(\mu)$. These markers define a map from a subset of
$S^1_\infty(\lambda)$ to a dense subset of $S^1_\infty(\mu)$, by sending one end of
the marker to the other. Notice that long markers do not cross. For, genuine markers
do not cross, and long markers are Hausdorff limits of genuine markers. It follows
that if two long markers cross, they could be approximated near the crossing
by genuine markers which would be forced to cross, contradicting the (stronger)
lemma~\ref{markers_disjoint}. It follows that the map defined by long markers
preserves the circular ordering on triples of
points in $S^1_\infty(\lambda)$ and $S^1_\infty(\mu)$. For, if $H_1,H_2,H_3$ are
three long markers, we can find a $1$--parameter family $\Delta_I$ of ideal triangles
in the interval $I \subset L$ from $\lambda$ to $\mu$ which are asymptotic in each
leaf to the endpoints of the $H_i$. Since markers do not cross at
infinity, the triangles might degenerate, but their orientation cannot
reverse, relative to the orientations on the leaves of the $\til{\F}$.
Consequently the circular orderings of
the endpoints of $\Delta_\lambda$ and $\Delta_\mu$ agree if they are both nondegenerate.

It follows that there is a unique continuous extension of this map to a
monotone map $\pi_\mu^\lambda$ from $S^1_\infty(\lambda)$ to 
$S^1_\infty(\mu)$ obtained by quotienting out complementary intervals
to $M_l(\lambda,\mu)$.
Said another way, the intervals $I \subset S^1_\infty(\lambda)$ quotiented out
by $\pi_\mu^\lambda$ do not intersect any long marker from $\mu$ to $\lambda$ in their
interior, but their endpoints can be approximated by long markers from $\mu$ to $\lambda$.
The laminar relation $\Lambda^+(\lambda,\mu)$ can be seen to be the pullback of
the trivial laminar relation of $S^1_\mu$ under $(\pi_\mu^\lambda)^{-1}$.

By theorem~\ref{dense_markers}, the set of long markers is dense in $E_\infty$.
If $\lambda > \mu > \nu$, the set of long markers from $\nu$ to $\lambda$ restricts to
the set of long markers from $\nu$ to $\mu$, and therefore these monotone
maps satisfy a cycle condition
$$\pi_\nu^\lambda = \pi_\nu^\mu \circ \pi_\mu^\lambda$$
on a set of points whose image is dense in $S^1_\infty(\nu)$, and therefore
by the uniqueness of monotone extensions, the cycle condition is satisfied everywhere.

It follows that we can obtain $S^1_\u$ as an inverse limit of the $\pi_\mu^\lambda$.
By the definition of inverse limit, for each $\lambda$ there is a map
$\phi_\lambda:S^1_\u \to S^1_\infty(\lambda)$ satisfying
$\phi_\mu = \pi^\lambda_\mu \circ \phi_\lambda$ for each pair $\lambda > \mu$.

Now, for $\lambda > \mu$, we claim $\pi_\mu^\lambda(\Lambda^+(\lambda)) = \Lambda^+(\mu)$
where these should be thought of as laminar relations of the corresponding circle
at infinity. For, if $\lambda > \nu$ is arbitrary, one of the following possibilities
occurs:
\begin{enumerate}
\item{The leaves $\mu,\nu$ are incomparable. In this case, the long markers from $\nu$ to
$\lambda$ are in a complementary interval in $S^1_\infty(\lambda)$ from
$\mu$ to $\lambda$, and therefore $\pi_\mu^\lambda(\Lambda^+(\lambda,\nu))$ is the
trivial laminar relation of $S^1_\infty(\mu)$.}
\item{We have $\lambda>\nu>\mu$. Then every interval quotiented out under
$\pi_\nu^\lambda:S^1_\infty(\lambda) \to S^1_\infty(\nu)$ is a subinterval of an
interval quotiented out under $\pi_\mu^\lambda:S^1_\infty(\lambda) \to S^1_\infty(\mu)$.
It follows that $\pi_\mu^\lambda(\Lambda^+(\lambda,\nu))$ is the trivial laminar
relation of $S^1_\infty(\mu)$.}
\item{We have $\lambda > \mu > \nu$. Then the set of long markers from $\nu$ to $\lambda$
restricts to the set of long markers from $\nu$ to $\mu$. It follows that if $I$ is
a complementary interval to $M_l(\lambda,\nu)$, the image of $I$ under 
$\pi_\mu^\lambda$ is a complementary interval to $M_l(\mu,\nu)$. It
follows that the endpoints of $\pi_\mu^\lambda(I)$ bound a leaf of
$\Lambda^+(\mu,\nu)$. Conversely, the endpoints of every such leaf bound an
interval in the image of some $I$ as above. Thus 
$$\pi_\mu^\lambda(\Lambda^+(\lambda,\nu))= \Lambda^+(\mu,\nu)$$}
\end{enumerate}
It follows that a dense subset of $\pi_\mu^\lambda(\Lambda^+(\lambda))$ is contained
in $\Lambda^+(\mu)$, and moreover the image contains a dense subset of $\Lambda^+(\mu)$.
Since both laminar relations are closed, they are equal, and 
$$\pi_\mu^\lambda(\Lambda^+(\lambda)) = \Lambda^+(\mu)$$
for all $\lambda > \mu$.

Now, for $\lambda_i$ an increasing unbounded sequence of leaves, the monotone maps
$\phi_{\lambda_i}:S^1_\u \to S^1_\infty(\lambda_i)$ quotient out smaller and smaller
subintervals, by the definition of an inverse limit. It follows that the laminar relations
$\phi_{\lambda_i}^{-1}(\Lambda^+(\lambda_i))$ are an increasing union of compact
subsets of $S^1_\u \times S^1_\u$,
and converge in the Hausdorff topology to a laminar relation $\Lambda^+_\u$ of
$S^1_\u$. By construction, $\phi_\lambda(\Lambda^+_\u) = \Lambda^+(\lambda)$ for
each $\lambda$. Moreover we have
$$\lim_{i \to \infty} \phi_{\lambda_i}^{-1} \circ
\phi_{\lambda_i}(\Lambda^+_\u) = \Lambda^+_\u$$
converging as compact subsets of $S^1_\u \times S^1_\u$.

Now, any two increasing sequences of leaves are eventually comparable, and therefore
the natural action of $\pi_1(M)$ on the set of increasing unbounded sequences of
leaves induces an action on $S^1_\u$ which preserves the laminar
relation $\Lambda^+_\u$.
\end{pf}

\begin{lem}\label{produces_genuine_lamination}
Let $\Lambda_\u$ be an invariant laminar relation in $S^1_\u$. Then
for each leaf $\lambda$ of $\til{\F}$ there is a geodesic
lamination $\Lambda(\lambda)$ determined
by the laminar relation $\phi_\lambda(\Lambda_\u) \subset S^1_\infty(\lambda)$.
Then 
$$\til{\Lambda} = \bigcup_\lambda \Lambda(\lambda)$$
is a branched lamination, which can be split open to a
lamination of $\til{M}$ which intersects leaves
of $\til{\F}$ in $k,\epsilon$--quasigeodesics for any $k>1,\epsilon>0$, and
covers an essential lamination in $M$.
\end{lem}
\begin{pf}
To see that $\til{\Lambda}$ is a branched lamination, we must first show that
$\Lambda(\lambda)$ varies continuously with $\lambda$, in the geometric
topology. Let $\tau(t)$ for $t \in I$ be a transversal to $\til{\F}$,
intersecting leaves $\lambda_t$. Parameterize $\tau$ as the interval $[-1,1]$, so that
$\lambda = \lambda_0$ is an interior leaf. Let $C_\infty$ be the cylinder
$E_\infty|_\tau$. 

We claim that as $t_i \to 0$ for $t_i$ an increasing sequence, the set of long markers
from $\lambda_{t_i}$ to $\lambda$ converges
to all of $S^1_\infty(\lambda)$. For, by theorem~\ref{dense_markers}, for any filter
$\mu_i < \lambda$, the set of long markers from $\lambda$ to some $\mu_i$ is dense
in $S^1_\infty(\lambda)$. But for any $i$, there is a $t_j$ such that for $t_k>t_j$,
$\lambda_{t_k}>\mu_i$. This proves the claim. In particular, for any $\epsilon$,
there is an $i$ such that for $j>i$, the set of endpoints of markers from $\lambda_{t_j}$ to
$\lambda$ is an $\epsilon$--net in $S^1_\infty(\lambda)$, in the visual metric as seen
from $\tau(0)$.

Conversely, if $t_i \to 0$ for $t_i$ a decreasing sequence, the set of long markers
from $\lambda$ to $\lambda_{t_i}$ is eventually arbitrarily dense in 
$S^1_\infty(\lambda_{t_i})$. That is, for any $\epsilon$, there is an $i$ such
that for $j>i$, the set of endpoints of markers from $\lambda$ to $\lambda_{t_j}$ is
an $\epsilon$--net in $S^1_\infty(\lambda_{t_j})$, in the visual metric as seen
from $\tau(t_j)$. For, the set of long markers is dense in $S^1_\infty(\lambda)$,
so pick a finite set $M_k$ which are an $\epsilon/2$--net in $S^1_\infty(\lambda)$.
The visual angle between the endpoints of a pair of markers as seen
from $\tau(t)$ varies continuously as a function of $t$, so for some sufficiently
small $\delta$, when $t<\delta$, the visual angle between endpoints of
pairs of consecutive markers $M_k,M_{k+1}$ is at most $\epsilon$. Choosing
$\tau_i < \delta$ proves the claim.

In particular, the map 
$$\Pi_\lambda:S^1_\infty(\lambda_1) \times I \to C_\infty \text{ defined by }
\Pi_\lambda(\cdot,t) = \pi^{\lambda_1}_{\lambda_t}(\cdot)$$
is {\em continuous} at $0$. Since $0$ was an arbitrary choice in our parameterization
of $I$, this map is continuous everywhere.

It follows that the laminar relations 
$$\Lambda(\lambda_t) = \phi_{\lambda_t}(\Lambda_\u) =
\Pi_\lambda(\cdot,t)(\Lambda(\lambda_1))$$
vary continuously as a function of $t$. Since the corresponding geodesic laminations
of leaves of $\lambda_t$ are determined on any compact subset in a continuous
way by the laminar relations at infinity, and since the leaves of $\til{\F}$ vary
continuously on compact subsets, the geodesic laminations $\Lambda(\lambda_t)$ vary
continuously in the Hausdorff topology on compact subsets.

We conclude that the $\Lambda(\lambda)$ sweep out a lamination $\til{\Lambda}$ of
$\til{M}$, with the following proviso: it is possible for distinct leaves
of $\Lambda(\lambda)$ to be identified to a single leaf of
$\Lambda(\mu)$ for $\lambda > \mu$, so that $\til{\Lambda}$ might merely
be a branched lamination, whose branch locus is a $1$--manifold (i.e. the branch
locus does not cross itself) and leaves branch open only in the positive
direction. By the equivariance of the construction, $\til{\Lambda}$ covers
a branched lamination $\Lambda$ of $M$, which branches in the positive direction along
a $1$--manifold.

The fact that the
branch locus is a manifold implies that $\Lambda$ can be split open to a lamination;
but the fact that the branching is all in one direction means that this splitting
can be done canonically, arbitrarily close to the identity. We call the split
open lamination $\Lambda$, and its universal cover $\til{\Lambda}$, by abuse of
notation. Notice that the complementary regions to $\til{\Lambda}$ correspond to the
complementary regions to the geodesic lamination $\Lambda_\u$ of $\H^2$ bounded
by $S^1_\u$. Since the splitting is done in the negative direction, distinct complementary
regions to $\til{\Lambda}$ are not joined together by the splitting and the complementary
regions still correspond to complementary regions to $\Lambda_\u$.

Leaves carried by $\til{\Lambda}$ before the splitting are planes transverse to
$\til{\F}$ and intersecting leaves of $\til{\F}$ in geodesics. The leaves in
the split open lamination $\til{\Lambda}$ are also planes transverse to
$\til{\F}$, and may be taken to intersect leaves of $\til{\F}$ in $k,\epsilon$
quasigeodesics, for any uniform choice of $k>1$ and $\epsilon >0$. It follows
that leaves of $\Lambda$ are incompressible. Moreover, since $\til{M} = \R^3$,
complementary domains are irreducible. It is easy to see that end incompressibility
is satisfied, so $\Lambda$ is essential, as required.
\end{pf}

The object $\til{\Lambda}^+$ defined as the union of the $\Lambda^+(\mu)$ is
$\pi_1(M)$--equivariant, since its definition involves no choices.
Therefore it covers an object in $M$, which we refer to as
$\Lambda^+$. That is, $$\Lambda^+ = \pi(\til{\Lambda}^+)$$

\begin{thm}\label{pinching_lamination_genuine}
$\Lambda^+$ can be split open to a genuine lamination of $M$ transverse to $\F$, and
intersecting leaves of $\F$ in $k,\epsilon$ quasigeodesics for any $k>1,\epsilon> 0$.
\end{thm}
\begin{pf}
The intersection of $\til{\Lambda}^+$ with each leaf $\lambda$ is a
geodesic lamination whose corresponding laminar relation is the image
of $\Lambda^+_\u$ by the monotone map $\phi_\lambda$.

It follows from lemma~\ref{produces_genuine_lamination}
that the laminations $\Lambda^+(\lambda)$ vary continuously from
leaf to leaf, and sweep out a
(possibly branched) essential lamination of $\til{M}$ which is 
natural and therefore covers a (possibly branched)
essential lamination of $M$ which can be split open to give an essential
lamination, which by abuse of notation we call $\Lambda^+$.

$\Lambda^+$ could only fail to be genuine
if $\Lambda^+(\lambda)$ was actually a foliation for each $\lambda$.
But in this case for each $\lambda$ there is an interval $I$ in
$S^1_\infty(\lambda)$ where the endpoints of the leaves of the
foliation are all distinct. Let $p_i \in \lambda$ be a sequence of
points converging to a point in the interior of $I$. Then we can translate
the $p_i$ to a fixed fundamental domain of $M$ and extract a convergent
subsequence with the property that the foliation of the 
limiting leaf $\lambda'$ is by all the geodesics asymptotic in one
direction to a single point $q \in S^1_\infty(\lambda')$. We call this the
{\em asymptotic foliation} of $\H^2$.

Notice that this point is an endpoint of $\Lambda^+(\lambda',\mu)$
for all $\mu < \lambda'$; in particular, it is a limit of markers, and it is
clear therefore that the ``asymptotic point'' sweeps out a long marker.
As in the proof of theorem~\ref{dense_markers} this is
absurd for a foliation with branching. Hence the geodesic lamination
$\Lambda^+(\lambda)$ has complementary regions with $>2$ boundary components for
some (and therefore all) $\lambda$, and the same is true of $\til{\Lambda}^+$. In
particular, $\Lambda^+$ is genuine, as required.
\end{pf}

In the sequel we will distinguish $\Lambda^+$ before and after splitting by letting
$\Lambda^+_b$ denote the (possibly) branched lamination, and $\Lambda^+$ the result of
splitting it.

\begin{rmk}\label{actually_no_split}
In fact, we will see in corollary~\ref{no_tangential_annulus} that no splitting
is necessary for $\til{\Lambda}^+$ --- that is, the lamination we construct
has no branching, and does not need to be split open.
\end{rmk}

\begin{rmk}
It is actually not hard to show that the lamination produced by {\em any} invariant
laminar relation $\Lambda_\u$ of $S^1_\u$ is genuine, by following the strategy
of theorem 4.6.4 in \cite{dC99b}. Roughly, one shows that if $\Lambda_\u$ gives
rise to a lamination which is essential but not genuine, $\Lambda_\u$ is
the asymptotic foliation. It is not too hard to show in this case that
the corresponding asymptotic points $p(\lambda) \in S^1_\infty(\lambda)$ are
directions in which nearby leaves of $\til{\F}$ diverge at infinity; but this quickly implies
that leaves of $\til{\F}$ converge in {\em every other} direction at infinity,
and one obtains a contradiction by the methods of section 2. The proof is somewhat
streamlined for $\Lambda^+_\u$, and therefore we give a complete proof of this case in 
theorem~\ref{pinching_lamination_genuine}.
\end{rmk}

\begin{cor}
Suppose $\F$ is a taut foliation with one--sided branching of an atoroidal
$3$--manifold $M$. Then $\pi_1(M)$ is $\delta$--hyperbolic in the sense of
Gromov.
\end{cor}
\begin{pf}
This is an immediate corollary of theorem~\ref{pinching_lamination_genuine}
and the main result of \cite{dGwK98}.
\end{pf}

\begin{rmk}
The main result advertised in \cite{wT98} is the existence (for an arbitrary
oriented and co--oriented taut foliation) of a universal
circle $S^1_\u$ with an action of $\pi_1(M)$ on it, which maps monotonely
to each $S^1_\infty(\lambda)$. No proof of this fact is actually written
in \cite{wT98}, but a program to produce such a proof has been outlined 
by Thurston in several talks and personal communications. A sketch of an
alternative construction is also given in \cite{CalPhD}.

Theorem~\ref{universal_circle} substantially strengthens this 
claim of Thurston for foliations with one--sided 
branching; in particular, the existence of monotone maps 
$\pi^\lambda_\mu:S^1_\infty(\lambda) \to S^1_\infty(\mu)$ compatible with
the $\phi_\lambda$ demonstrates that our universal circle is the
{\em smallest possible} (i.e. there is a monotone $\pi_1(M)$--equivariant
map from any other universal circle to ours) and our construction makes
the relationship of this circle to the extrinsic geometry of $\til{\F}$
much more explicit.
\end{rmk}

\section{The bending lamination}

\subsection{Geometry of the pinching lamination}

We analyze in more detail the geometry of $\Lambda^+$. From this we
can get more detailed information about the representation of $\pi_1(M)$
in $\hom(S^1_\u)$, and thereby produce another lamination $\Lambda^-_\u$
of $S^1_\u$ which will complement $\Lambda^+_\u$. It is useful to think of this
complementary lamination as arising from the tension between the intrinsic and
extrinsic geometry of leaves $\lambda$ of $\til{\F}$ in the subspaces $\til{M}^\lambda$.
Leaves of $\til{\F}$ are uniformly properly embedded below (though not in
general {\em quasigeodesically} embedded below). Consequently, the leaves {\em bend}
in the negative direction, but they do not {\em pinch} and give rise to branching.

We use some technology here from the theory of genuine laminations. For details
of this theory, see \cite{dGwK98}. Basically, we use the fact that the complementary
regions to a genuine lamination can be partitioned into {\em finitely
many} compact {\em gut regions}, and
$I$--bundle {\em interstitial regions} which meet each other along compact
{\em interstitial annuli}. This partition is not unique, but there are several
canonical choices, including a {\em maximal gut} partition, where the $I$--bundles
are chosen to be as small as possible, and are all noncompact, and a 
{\em minimal gut} partition, where the union of the $I$--bundles is the
{\em characteristic} $I$--bundle of the complementary region, in the sense of
Jaco--Shalen and Johannson.

\begin{thm}\label{solid_torus_guts}
Let $\Lambda$ be a genuine lamination transverse to $\F$ intersecting
$\F$ in $(k,\epsilon)$--quasigeodesics. Then
the interstitial annuli of $\Lambda$ can be straightened to be
transverse or tangential to $\F$.
Moreover, either $M$ is toroidal and one can find a reducing
torus, or $\Lambda$ has solid torus guts.
\end{thm}
\begin{pf}
Let $\fG$ be a gut region complementary to $\Lambda$ bounded by some finite
collection $A_i$ of interstitial annuli. Let $\til{A}_i$ be a lift of
$A_i$ to $\til{M}$, and let $\gamma_i$ be the element of $\pi_1(M)$, acting
by translation on $\til{A}_i$ so that $\til{A}_i/\gamma_i = A_i$. Since
$\til{\F}$ does not branch in the positive direction, a fundamental
domain for the intersection
$\partial \til{A}_i \cap \lambda$ for $\lambda$ a boundary leaf of
$\til{\Lambda}$ can be isotoped in the positive direction until the
endpoints are on comparable leaves, and then straightened to be
transverse or tangential. 
This can be done with both sides of $\partial A_i$, and then
the annulus straightened leafwise on its interior to be transverse or
tangential. For, $\lambda$ is transverse to $\til{\F}$, and
therefore in $\lambda$ one can isotope in the positive direction.
This straightening is more or less the same as the procedure carried out in
lemma~\ref{straighten_transversals}.

Now the boundary $\partial \fG$ of the gut region is a compact
surface transverse to $\F$ or possibly tangent along a collection of
annuli, and therefore can be given a nonsingular one--dimensional
foliation. It follows that $\partial \fG$ is a torus or Klein bottle. If
this torus is essential, $M$ is toroidal. Otherwise $\partial \fG$ is 
compressible and bounds a solid torus, since $\pi_1(A_i) \to \pi_1(\partial \fG)$
is injective. One sees immediately that $\fG$ is a solid torus, with the
structure of a finite sided polygon bundle over $S^1$.
\end{pf}

\begin{rmk}
It is a special property of taut foliations with one--sided or no branching
that loops can be homotoped to be transverse or tangential. For foliations
with two--sided branching, one can always find loops which are not homotopic
to transverse loops. This fact is exploited in \cite{dC99c} to explore the
space of isotopy classes of taut foliations on a $3$--manifold with the
geometric topology.
\end{rmk}

\begin{rmk}
If $M$ is atoroidal, since the guts are solid tori for {\em any} choice of partition of
complementary regions into guts and interstitial regions,
it is well--known (see e.g. \cite{dGuO87}) that in this case the complementary
regions $G$ to $\Lambda$ are ideal polygon bundles over $S^1$, and the interstitial
regions are products $S^1 \times \R^+ \times I$ corresponding to the cusps of
the ideal polygons $\times S^1$. One typically says in this case that $\Lambda^+$ is
{\em very full}.
\end{rmk}

It follows that the complementary regions to $\Lambda^+$ are finite
sided ideal polygon bundles over $S^1$. We fix the following notation in the
sequel: $G$ will denote a complementary region to $\Lambda^+$, and $\til{G}$
will denote a complementary region to $\til{\Lambda}^+$ covering $G$ with gut $\fG$.
Since $G$ is a bundle, some finite cover of $G$ is topologically a product
$P \times S^1$ for some ideal polygon $P$, and $\til{G}$ is topologically
a product $P \times \R$. The boundary of the gut of $G$ is a solid torus
foliated by the intersection with $\F$ away from the annuli $A_i$ which
are tangent to $\F$. These annuli can be foliated as products, to give a nonsingular
foliation of the boundary of a gut $\fG$ of $G$. On each annulus of
$\Lambda^+ \cap \partial \fG$ the foliation is a product, which is transverse along
$\partial A_i$ where $A_i$ is transverse to $\F$, and is asymptotic to $\partial A_i$
where $A_i$ is tangent to $\F$. The foliations of the $A_i$ are products --- either
meridional in case $A_i$ is transverse to $\F$, or longitudinal if $A_i$ is
tangent to $\F$.

If $\lambda$ is a leaf of $\til{\F}$ intersecting the
gut part of $\til{G}$, the intersection
$\lambda \cap \til{G}$ is some (possibly infinite)
polygon $P'$. Its projection to $M$ is transverse to $\partial \fG$ away
from the $A_i$, but might spiral around the $A_i$ which are tangent to $\F$.
If an end of $P'$ stays in the gut $\fG$ of 
$G$, it must spiral towards some $A_i$. An end of $P'$ is topologically
a product; if it stays in $\fG$, it can't get too thin. A pair of geodesics in $\H^2$
which stay close together on some noncompact set must actually be asymptotic; it follows
that the end of such a $P'$ spirals around some blown--up branch locus of $\Lambda^+_b$. 
That is, the annulus $A_i$ is obtained by splitting open $\Lambda^+_b$ along a branch circle.
Conversely, if $A_i$ is a tangential annulus, the two boundary curves are isotopic, and
collapse to a branch circle when straightened to leafwise geodesics.
The end of $P'$ has infinite area in the
split open lamination, but had finite area in the branched lamination. It follows that
$P'$ is combinatorially equivalent to $P$. On the other hand, if $P'$ is contained in
the thin part of $\til{G}$, it is an infinite bigon, created when $\til{\Lambda}^+$ was
split open. 

\begin{figure}[ht]
\scalebox{.6}{\includegraphics{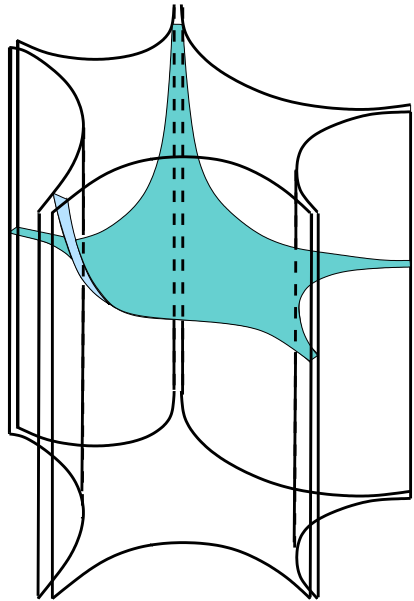}}
\caption{A complementary polygon $P' \subset \lambda$ intersects $G$
in a polygon which might be transverse to some interstitial annuli $A_i$
and spiral around some other $A_j$ corresponding to branch circles of $\Lambda^+_b$.}
\end{figure}

Notice, since $\F$ only branches in the negative direction, polygons
$P'$ of $\lambda \cap \fG$ spiral in the {\em positive direction} around the
tangential $A_i$, as in figure 4. Moreover, the foliation of 
$(\partial \fG \backslash \cup \text{tangential } A_i)$
is filled in with a foliation of $\fG$ by
disks of the form $P'$ as above. In particular, this is a {\em product} foliation away
from the tangential $A_i$. That is, the foliation of $\partial \fG$ is a product foliation
by circles turbularized in the positive direction along the tangential annuli, and
the foliations induced by $\F$ on each annulus of $\Lambda^+ \cap \partial \fG$ is
a product foliation by closed intervals.

Of course, in corollary~\ref{no_tangential_annulus}, we will show that $P'$
will actually be {\em transverse} to all the $A_i$, so that no turbularization is
necessary.

We define a {\em complementary region $P$ to $\Lambda^+_\u$} to be a complementary
region to a geodesic lamination of $\H^2$ determined by the laminar relation $\Lambda^+_\u$.
If $\phi:S^1_\u \to S^1$ is a monotone map, the image of the vertices of $P$
are the vertices of a complementary region to $\phi(\Lambda^+_\u)$, or else they
consist of either two or one ideal vertices. In any case, it makes sense to talk about
the image of a complementary region under a monotone map.

Our characterization of the intersections $\til{\fG} \cap \lambda$ for $\fG$ a gut of
$\Lambda^+$ and $\lambda$ a leaf of $\til{\F}$ has
implications for the monotone maps $\phi_\lambda:S^1_\u \to S^1_\infty(\lambda)$.
In particular, for each complementary polygon $P$ to $\Lambda_\u$, the map
$\phi_\lambda$ does one of the following three things:
\begin{enumerate}
\item{The map $\phi_\lambda$ does not collapse any boundary leaves, and there is a
complementary region to $\Lambda^+(\lambda)$ corresponding to $P$ which is an ideal
polygon with the same combinatorics as $P$. This happens for all sufficiently high leaves
$\lambda$.}
\item{The map $\phi_\lambda$ collapses $P$ to a single geodesic in $\Lambda^+(\lambda)$.
In the split open lamination, these bound an infinite rectangle. If $\lambda$ is
the uppermost leaf with this property, the corresponding geodesic in $\Lambda^+(\lambda)$ is a
branch line of the branch locus of $\til{\Lambda}^+$ before it is split open. 

In theorem~\ref{circle_branch_locus} we will show that 
this geodesic is stabilized by some non--trivial element of $\pi_1(M)$
and covers a branch circle in $\Lambda^+$ before splitting open. After splitting open,
this becomes a tangential interstitial annulus $A_i$. Together with an analysis of the
combinatorics of the action of $\pi_1(M)$ on $S^1_\u$, we will actually see that this
case cannot occur.}
\item{The map $\phi_\lambda$ collapses $P$ to a single ideal point, and there is no
corresponding complementary region or geodesic in $\Lambda^+(\lambda)$.}
\end{enumerate}

It follows that for each complementary polygon $P$ to $\Lambda^+(\lambda)$, and each
$\lambda > \mu$, the monotone map
$$\pi^\lambda_\mu:S^1_\infty(\lambda) \to S^1_\infty(\mu)$$
either preserves the combinatorics of $P$, collapses $P$ to a geodesic (or infinite
rectangle after splitting), or collapses $P$ to a single ideal point.

\begin{cor}\label{finite_sided}
If $M$ is atoroidal,
the complementary regions to $\Lambda^+_\u$ are finite sided ideal polygons.
\end{cor}
\begin{pf}
The complementary regions to $\Lambda^+_\u$ correspond to the complementary
regions to $\til{\Lambda}^+$, which are lifts $\til{G}$ of complementary regions $G$
to $\Lambda^+$. For $\lambda$ sufficiently high up, the
intersections $\lambda \cap \til{G}$ are finite sided polygons
$P'$ with the same combinatorics as $P$. It follows that for $\lambda_i$
an increasing unbounded sequence intersecting $\til{G}$ higher and higher up,
the complementary regions to $\Lambda^+(\lambda_i)$ which intersect $\til{G}$
are ideal polygons isomorphic to $P$. It follows that the complementary region
to $\Lambda^+_\u$ corresponding to $\til{G}$ is an ideal polygon isomorphic to $P$.
\end{pf}

\begin{cor}
If $M$ admits a taut foliation with one--sided branching, then
any homeomorphism $f:M \to M$ homotopic to the identity is isotopic to the
identity.
\end{cor}
\begin{pf}
If $M$ is toroidal it is Haken, and Waldhausen's theorem applies. Otherwise
one can construct $\Lambda^+$ which has solid torus gut regions. Now apply
the main theorem of \cite{dGwK93}.
\end{pf}

In the sequel we will assume without comment that $M$ is atoroidal, and therefore 
corollary~\ref{finite_sided} applies.

\begin{defn}
Let $X$ be a line transverse to $\til{\F}$. We say $X$ is {\em regulating
in the positive direction} if some (and therefore every) positive ray contained
in $X$ projects {\em properly} to $L$. Equivalently, if $X$ intersects a leaf
$\mu$ of $\til{\F}$, it intersects all $\lambda$ with $\lambda > \mu$. Similarly, $X$ is
{\em regulating in the negative direction} if some {\em negative} ray contained in
$X$ projects properly to $L$.

If $\alpha$ is a loop in $M$ transverse to $\F$, we say $\alpha$ is {\em regulating
in the positive direction} if every component of the preimage $\til{\alpha} \subset
\til{M}$ is regulating in the positive direction, in the sense above.
\end{defn}

Let $G$ be a complementary region to $\Lambda^+$, and let $\alpha$ be the
core circle of $G$. Some finite cover of $G$ is a product $P \times S^1$,
so some power of $\alpha$ fixes the vertices of the corresponding polygon $P$
in $S^1_\u$.

If $\til{G}$ is a lift of $G$, if a leaf $\mu$ of $\til{\F}$ intersects
the gut of $\til{G}$, the same is true for all higher leaves $\lambda>\mu$. In
particular, the lift $\til{\alpha}$ of $\alpha$ is {\em regulating} in
the positive direction. 

Suppose $\Lambda^+$ is not minimal. Then a minimal sublamination $\Lambda^+_m$ has
solid torus guts, by theorem~\ref{solid_torus_guts}, and these must be finite unions
of solid torus guts of $\Lambda^+$ glued along finitely many isolated leaves.
It is clear that a {\em unique} such minimal sublamination exists.
It follows that $\Lambda^+_m$ is obtained from $\Lambda^+$ by removing these
finitely many isolated leaves. In terms of the universal circle, $\Lambda^+_\u$
is obtained from $(\Lambda^+_m)_\u$ by adding some diagonals to the complementary
ideal polygons.

If $\alpha$ is a core circle of a gut region of $\Lambda^+_m$, $\alpha$ is
also regulating in the positive direction. Moreover, a gut region $\fG$ of
$\Lambda^+_m$ is foliated as a bundle by polygons contained in leaves of $\F$ as
illustrated in figure 4. If the projection of the negative end of some lift
$\til{\alpha}$ to $L$ is proper, the interstitial annuli $A_i$ are all transverse
to $\F$, and isotopic to some finite cover of $\alpha$. If there were a diagonal
leaf $l$ of $\til{\Lambda}^+\backslash \til{\Lambda}^+_m$ in such a gut, the intersections
$l \cap \lambda$ for $\lambda \in \til{\F}$ would be leaves of $\Lambda^+(\lambda,\mu)$
for some $\lambda > \mu$, since $l$ is {\em isolated}, and therefore its intersections
cannot be a nontrivial limit of leaves of $\Lambda^+(\lambda,\cdot)$. Since $l$
is isolated, the uppermost leaf $\mu$ for
which $l \cap \lambda \in \Lambda^+(\lambda,\mu)$ is
constant as $\lambda$ varies. But $\til{\fG}$
is regulating in the negative direction, and therefore intersects leaves $\nu$ which
are either incomparable with $\mu$ or are below $\mu$; in either case,
$\Lambda^+(\nu,\mu)$ would be empty, contrary to the fact that $l \cap \nu$ is
nontrivial. Thus we can prove the following lemma:

\begin{lem}\label{minimal_splits}
Let $\til{G}$ be a complementary region to $\til{\Lambda}^+_m$, and suppose there is a leaf
$l$ of $\til{\Lambda}^+\backslash \til{\Lambda}^+_m$ contained in $\til{G}$.
Then the core circle of $G$ is not regulating in the negative direction, and for all
leaves $\lambda$ of $\til{\F}$ intersecting $\til{G}$ in polygons, the
intersection $l \cap \lambda$ is a leaf of $\Lambda^+(\lambda,\mu)$ for some fixed
$\mu$ which is a lowermost limit of $\til{\fG}$.
\end{lem}
\begin{pf}
This is basically established in the discussion above, except for the last claim.
Of course, $\mu$ is below every leaf intersecting $\til{\fG}$.
On the other hand, an isolated leaf of $\Lambda^+(\lambda,\mu)$ can be continued
in the negative direction until $\mu$. Since the leaf $l$ is isolated and intersects
$\til{\fG}$ leafwise, it follows that $\til{\fG}$ can be continued in the negative
direction until $\mu$, as claimed.
\end{pf}
Of course, it is possible that some interstitial region of
$\til{G}\backslash\til{\fG}$ could be continued below $\mu$, so it is only the
gut which stops at $\mu$. Notice too that if $l_i$ are a (necessarily finite)
collection of leaves of $\til{\Lambda}^+$ in $\til{G}$, the set of $\mu_i$
that they limit to are all limits of the ordered decreasing sequence of leaves that
$\til{\fG}$ intersects.

Now, the lamination $\Lambda^+_m$ is minimal, and therefore every leaf of
$(\Lambda^+_m)_\u$ is a limit of a sequence of translates of any other leaf.
In particular, no leaf is isolated. Notice that every non--isolated leaf of
$\Lambda^+_\u$ is a leaf of $\Lambda^+_m$.

\subsection{Combinatorics of the action on $S^1_\u$}

In order to pass flexibly between laminar relations and geodesic laminations,
we will refer to a pair of points in a circle related by some laminar relation as
a {\em leaf} of that laminar relation. These are the endpoints of the corresponding
leaf of the corresponding geodesic lamination of $\H^2$ bounded by that circle.
It should be clear from context which meaning is appropriate in each case.

The following technical lemma shows that for a convergent
sequence of leaves $l_i \to l$ of $\Lambda^+_\u$, the endpoints of the $l_i$ are
distinct from the endpoints of the $l$. The proof of this lemma is annoyingly fiddly.
Perhaps there is a broader context for this lemma allowing a more robust proof.

\begin{lem}\label{limiting_leaves_disjoint}
Let $l$ be a leaf of $(\Lambda^+_m)_\u$, and let $l_i$ be a sequence of leaves
limiting to $l$. Then the endpoints of the $l_i$ are disjoint from the endpoints of $l$.
\end{lem}
\begin{pf}
Suppose there are $l_i \to l$ with the $l_i$ asymptotic to an endpoint $p$
of $l$ for sufficiently large $i$.
If there is {\em some} leaf $l'$ which can be approximated by $l_i'$ whose
endpoints are disjoint from $l'$, then if $\alpha(l')$ is sufficiently
close to $l$, $\alpha(l'_j)$ will intersect $l_i$ transversely, for sufficiently large
$j$. It follows that every leaf can be approximated by leaves asymptotic to
exactly one end, in particular the same is true of the $l_i$. Since the $l_i$ could be
taken to be non--isolated on either side, the $l_i$ can be approximated from
either side. It follows that the approximating leaves on either side are
asymptotic to the same endpoint $p$ of $l_i$. Moreover, if $l$ is a boundary leaf of some
polygon $P$, and $l'$ is the other boundary leaf of $P$ asymptotic to $l$ at $p$, then
$l'$ is approximated on the other side by leaves asymptotic to $l'$ at $p$.
For, a sequence $\alpha_i(l)$ of translates of $l$ converges to $l$, and the
sequence $\alpha_i(l')$ necessarily converges to $l$ too from the same side.
So the $\alpha_i$ all fix $p$, and it follows that the leaves approximating
$l'$ are asymptotic to $l'$ at $p$. For each leaf $l$ of $(\Lambda^+_m)_\u$ it
follows that there is one end, the {\em sticky end}, such that approximating
sequences $l_i \to l$ share a sticky end with $l$. Moreover, for each polygon $P$,
the vertices alternate between sticky ends and {\em free ends} (i.e. the non--sticky
ends of leaves of $\partial P$). We will assume in the sequel that $l$ is
a boundary leaf of $P$ such that the other vertices of $P$ are clockwise from $p$
before the free end of $l$.

Let $\L$ denote the set of leaves of $\Lambda^+_\u$ which are asymptotic to $p$.
A dense subset of these, $\L_0$, project by $\phi_\lambda$ to leaves of
$\Lambda^+(\lambda)$ which are of the form $\Lambda^+(\cdot,\mu)$ for
various $\mu$. Here we are just using the fact that $\Lambda^+(\lambda)$ is
the closure of the union of $\Lambda^+(\lambda,\mu)$ as $\mu$ varies.
If $l_i$ denotes a leaf of $\L$ for which
$\phi_\lambda(l_i) \subset \Lambda^+(\lambda,\mu_i)$ for some $\mu_i$,
then by the definition of $\Lambda^+(\lambda,\mu_i)$,
the intersections with $S^1_\infty(\lambda)$ of the
markers from $\mu_i$ to $\lambda$ converge to $\phi_\lambda(p)$. Let $\L_0^+$ be
the subset of $\L_0$ for which the markers converge to $p$ from the anticlockwise side, say.
That is, the intersection of the markers with $S^1_\infty(\lambda)$
converge to $\phi_\lambda(p)$ from the anticlockwise side.
 
Then the markers from $\mu_i$ and $\mu_j$ to $\lambda$ are interspersed with each other
for any two $i,j$ with $m_i \in \L_0^+$, since they both contain infinite subsets
whose intersections with $S^1_\infty(\lambda)$ converge to $\phi_\lambda(p)$ from the
anticlockwise side, so the leaves $\mu_i,\mu_j$ are
comparable.
If there is a $m_k$ such that the markers from $\mu_k$ to $\lambda$
converge to $\phi_\lambda(p)$ from the clockwise side, including the degenerate case that
the markers from $\mu_k$ to $\lambda$ near $\phi_\lambda(p)$ are actually
equal to $\phi_\lambda(p)$, there is a translate of $m_k$ between two
elements $m_i,m_j$ of $\L_0^+$. Since the markers from $m_i,m_j$ and $m_k$
converge to $p$ from opposite sides, the leaf $\mu_k$ is incomparable with
$\mu_i < \mu_j$, say. On the other hand, markers from $\mu_k$ converge to
the free end of $m_k$ which is contained between the free ends of $m_i$ and $m_j$, and
therefore markers from $\mu_k$ to $\lambda$ link markers from $\mu_j$ to $\lambda$,
which implies they are comparable. This is a contradiction, and therefore one can
conclude that either $\L_0^+ = \L_0$, or for {\em every} $m_k$, the set of
markers from $\mu_k$ to $\lambda$ near $\phi_\lambda(p)$ are actually equal to
$\phi_\lambda(p)$. In the second case, then for any $i$,
there is an actual marker from $\mu_i$ to $\lambda$ which intersects $S^1_\infty(\lambda)$
at $\phi_\lambda(p)$. In particular, any two $\mu_i,\mu_j$ are asymptotic to $\lambda$
at the same point at infinity, and therefore $\mu_i$ and $\mu_j$ come arbitrarily
close. But this implies they are comparable. It follows in either case
that all the leaves in $\L_0$
are of the form $\Lambda^+(\cdot,\mu_i)$ where the $\mu_i$ are all {\em comparable}.

We claim that the free ends of $\L$ converge to $p$ from either side. Since otherwise,
by the fact $\Lambda^+_\u$ is closed, there is an ``outermost'' $l \in \L$.
But this is either a boundary leaf of some
complementary region $P$, in which case another side of $P$ is further out than
$l$, or it can be approximated on the outermost side by leaves, whose sticky
end must necessarily be at $p$; such approximating leaves would again be further out than
$l$. In either case, no such outermost $l$ can exist, and the claim is proved.

Fix $P$, a complementary region asymptotic to $p$ along some boundary leaves
$l,l'$ with the free end of $l'$ anticlockwise from the free end of $l$ in the
complement of $p$. Since $\L_0$ is dense in $\L$, the free ends of $\L_0$ converge to
$p$ from either side. Let $\mu_i$ for $i \in \Z$ be a sequence of leaves associated to the 
$m_i \in \L_0$ whose free ends converge to $p$ from the clockwise and the
anticlockwise direction as $i \to -\infty$ and $i \to +\infty$.
The leaves $\mu_i$ are all comparable, and we might as well take $\mu_i < \mu_j$
for $i<j$. Moreover, as $i \to -\infty$, {\em all} endpoints of markers
$M(\lambda,\mu_i)$ for any fixed $\lambda$ converge to $\phi_\lambda(p)$.
In particular, the laminations $\Lambda^+(\lambda,\mu_i)$, being the boundary of the
convex hull of these endpoints, converges in the Hausdorff topology in
$\lambda \cup S^1_\infty(\lambda)$ to the single point $\phi_\lambda(p)$.
It follows that the sequence $\mu_i$ decreases without bound as $i \to -\infty$; for,
if the $\mu_i$ were all bounded below by some leaf $\kappa$ of $\til{\F}$, the
convex hull of each $\Lambda^+(\lambda,\mu_i)$ would contain the convex hull
of $\Lambda^+(\lambda,\kappa)$ which is nonempty, since there are infinitely many
markers from $\kappa$ to $\lambda$.

If $\gamma$ is a power of the core geodesic of the complementary region $G$
corresponding to $P$, the lift $\til{\gamma}$ is regulating in the positive direction.
Moreover, for all $m_i \in \L_0$ whose free end is anticlockwise from $P$,
the complementary region $\til{G}$ is contained in the region bounded by
$\Lambda^+(\lambda,\mu_i)$, and therefore $\phi_{\mu_i}(P)$ is a nondegenerate
polygon in $\mu_i$. In particular, $\mu_i$ intersects $\til{G}$, and the
sequence of leaves $\gamma^n(\mu_i)$ increases unboundedly as $n \to \infty$.

It follows that for every $i$ sufficiently large, and for every $j>i$, there is
some positive $n$ such that $\gamma^n(\mu_i)>\mu_j$. Consequently, the free ends
of $m_j,\gamma^n(m_i)$ and the point $p$ are anticlockwise ordered, for sufficiently
large $n$. That is, the images of the free ends of the $m_i$ converge under the action
of positive powers of $\gamma$ to $p$, from the anticlockwise side.
Thus the dynamics of $\gamma$ on the clockwise interval from $p$ to
the free end of $l'$ has no fixed point, and the images of every point
in this interval converge to $p$ under increasing powers of $\gamma$.
Conversely, since $l$ is fixed by $\gamma$, it follows that $\gamma$ fixes some
leaf $\mu$ greater than all $\mu_i$ corresponding to $m_i \in \L_0$
whose free end is clockwise from the free end
of $l$. For all $m \in \L_0$ whose free end is anticlockwise from the free end of $l$,
the leaves $m$ are quotiented to a point by $\phi_\mu$, but no leaf $k \in \L_0$
whose free end is clockwise from the free end of $l$ is quotiented to
a point. Since $\gamma$ fixes $\mu$, it acts on $\mu$ by a hyperbolic isometry, and
therefore has exactly two fixed points in $S^1_\infty(\mu)$. One of these must be
$\phi_\mu(p)$; the other is some point $q$. Under the action of $\gamma$ on
$S^1_\infty(\mu)$, the leaves $\gamma^i(k)$ converge to the geodesic from 
$\phi_\mu(p)$ to $q$
(here $i \to \pm \infty$ as appropriate). In particular, there is a sequence $n_i$
in $\L_0$ whose projection to $\mu$ converges to the geodesic $\phi_\mu(p)q$,
such that the free ends are moving anticlockwise.
Since the free ends are moving anticlockwise, this sequence $n_i$
corresponds to leaves of $\Lambda^+(\mu,\nu_i)$ for some increasing sequence of
leaves $\nu_i$ whose upper bound is $\nu$. We must have
$\nu < \mu$, since if $\nu_i \to \mu$, then $n_i$ would converge to the point 
$\phi_\mu(p)$, whereas in fact they converge to the geodesic $\phi_\mu(p)q$.
The action of $\gamma$ preserves
the set of all $\nu_i$ for which $n_i$ has free end clockwise from $q$,
and therefore it fixes their upper bound $\nu$.
So there is some $\nu < \mu$ which is fixed
by $\gamma$. As before, we can consider the action of $\gamma$ on $S^1_\infty(\nu)$,
which has exactly two fixed points $\phi_\nu(p),r$ and
translates points on either side of the fixed
points in different directions. But the action of $\gamma$ on $S^1_\infty(\mu)$ and
$S^1_\infty(\nu)$ commute with the monotone map $\pi^\mu_\nu$. So $\gamma$ fixes the
extremal points in $(\pi^\mu_\nu)^{-1}(r)$ which are disjoint from $q$, since
$\pi^\mu_\nu(q) = \phi_\nu(p)$. In particular,
$\gamma$ has at least $3$ fixed points on $S^1_\infty(\mu)$, which is absurd.
See figure 5.

\begin{figure}[ht]
\centerline{\relabelbox\small \epsfxsize 2.0truein
\epsfbox{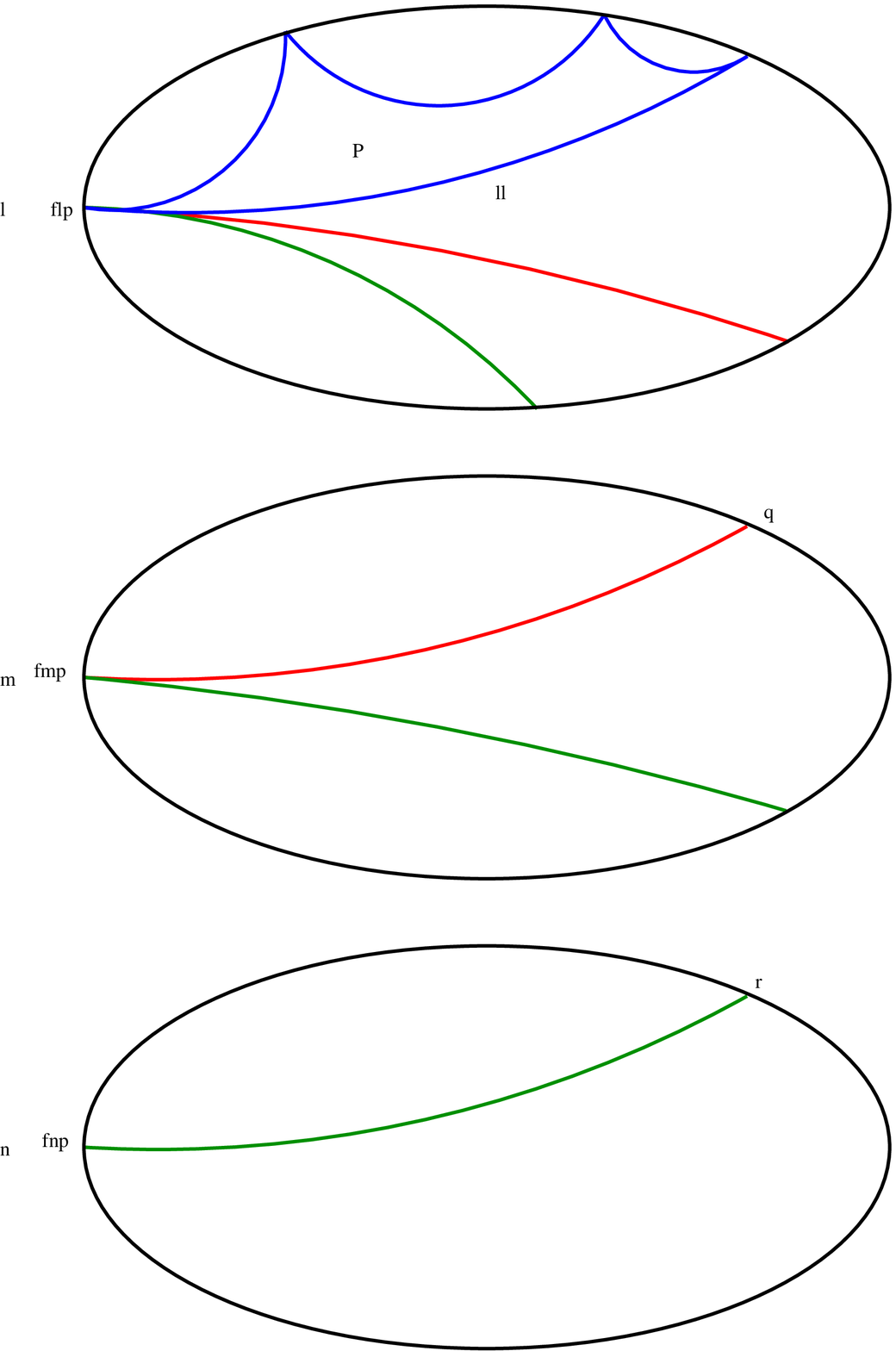}
\adjustrelabel <-3pt,0pt> {P}{$\phi_\lambda(P)$}
\adjustrelabel <0pt,-3pt> {ll}{$\phi_\lambda(l)$}
\adjustrelabel <0pt,20pt> {l}{$\lambda$}
\adjustrelabel <0pt,20pt> {m}{$\mu$}
\adjustrelabel <0pt,20pt> {n}{$\nu$}
\adjustrelabel <-20pt,0pt> {flp}{$\phi_\lambda(p)$}
\adjustrelabel <-50pt,0pt> {fmp}{$\phi_\mu(p) = \phi_\mu(l)$}
\adjustrelabel <-55pt,0pt> {fnp}{$\phi_\nu(p) = \pi^\mu_\nu(q)$}
\adjustrelabel <3pt,3pt> {q}{$q$}
\adjustrelabel <3pt,3pt> {r}{$r$}
\endrelabelbox}
\caption{$\gamma$ fixes $\mu$, and therefore fixes only the points $\phi_\mu(p)$ and $q$
in $S^1_\infty(\mu)$. But $\gamma$ leaves invariant an increasing
sequence of leaves $\nu_i$ such that $\Lambda^+(\mu,\nu_i)$
contains leaves converging to the geodesic $\phi_\mu(p)q$; in particular, $\gamma$
also fixes $\nu < \mu$. But then $\gamma$ fixes $r \in S^1_\infty(\nu)$, and therefore
fixes at least $3$ points in $S^1_\infty(\mu)$.}
\end{figure}

This contradiction proves the lemma.
\end{pf}

Part of the difficulty with proving this lemma is that the $3$--dimensional
origin of the lamination $\Lambda^+_\u$ is essential. We make the following definition:
\begin{defn}
A group $\Gamma$ of homeomorphisms of $S^1$ is
{\em sticky} if it satisfies the following properties:
\begin{enumerate}
\item{$\Gamma$ leaves invariant a lamination $\Lambda$ whose complementary regions are
finite sided ideal polygons which fall into finitely many orbit classes.}
\item{$\Lambda$ is minimal for the action of $\Gamma$
(i.e. the orbit of every leaf is dense).}
\item{Some (and therefore every) leaf of $\Lambda$ has a sticky end --- i.e. it is a limit
of a sequence of leaves which are asymptotic to one of its endpoints.}
\end{enumerate}
\end{defn}

The problem is that {\em sticky groups exist}, and even occur quite naturally. 

\begin{figure}[ht]
\scalebox{.6}{\includegraphics{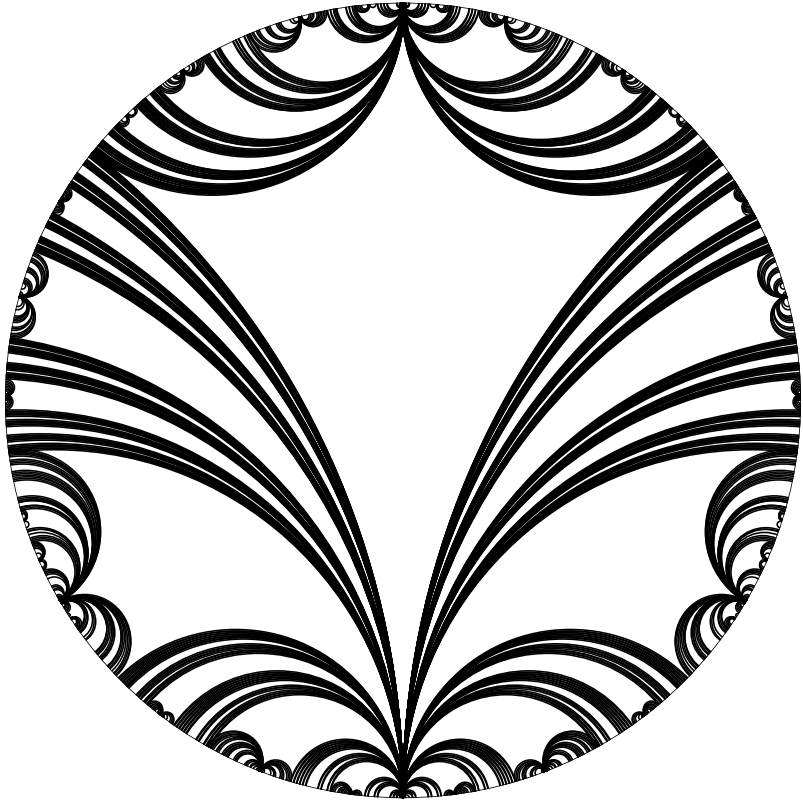}}
\caption{A lamination $\Lambda$ invariant under a sticky subgroup of $\hom(S^1)$}
\end{figure}

\begin{exa}
Let $T_\infty$ be the infinite $\infty$--valent tree.
For each vertex $v$, choose an identification
of the edges incident to $v$ with the dyadic rationals in $(-1,1)$. Thicken every edge
to the diagonal of a quadrilateral,
so that $T_\infty$ is embedded in an infinite $2$--complex
$\Sigma_\infty$. We embed $\Sigma_\infty \to \overline{\H^2}$ in the following way: the
vertices are mapped to ideal points, the squares map to an ideal quadrilateral, and the
ordering of the opposite vertices of the squares incident to a given vertex $v$ inherited
from $S^1_\infty - v$ should agree with the ordering
given by the identification with the dyadic
rationals. Such an embedding is easy to construct;
one particular choice is illustrated in figure 6.

For each ideal quadrilateral, the sticky ends correspond to vertices of $T_\infty$,
whereas the free edges correspond to nothing. If we build an abstract topological
space from the disjoint union of closures of complementary polygons, gluing two
such polygons along sticky vertices which map to the same point of $S^1$, the resulting
space deformation retracts to $T_\infty$.

The group of automorphisms of this lamination contains many interesting subgroups, such
as copies of Thompson's group, and many countable
finitely generated subgroups which act minimally.
For instance, the stabilizer of every vertex in $\Gamma$ admits a homomorphism to the
group of homeomorphisms of the dyadic rationals in $(-1,1)$; one can take the image to be
Thompson's group, and for appropriate $\Gamma$, one can find sections of this
homomorphism.
\end{exa}

In \cite{dC277} we actually prove that if $M$ is a $3$--manifold, and
$\pi_1(M)$ acts as a sticky group of homeomorphisms of $S^1$, then $M$ is Haken.

\vskip 12pt

Recall in lemma~\ref{produces_genuine_lamination} that we produced the genuine
lamination $\Lambda^+$ of $M$ from an {\it a priori} branched lamination, which
we denote here by $\Lambda^+_b$, whose branch locus consists of a union of $1$--manifolds
tangent to $\F$. These components are either circles or lines; a branch circle $\gamma$ is
split open to a tangential annulus $\gamma_s$ which can be taken to be a boundary annulus 
of a gut region. A branch line $\gamma$ is split open
to a tangential rectangle $\gamma_s$, which can be taken to be contained in an interstitial region.
In this case, the interstitial annulus bounding $\gamma_s$ away from the gut can be taken to
be transverse to $\F$. See figure 7.

\begin{figure}[ht]
\centerline{\relabelbox\small \epsfxsize 4.0truein
\epsfbox{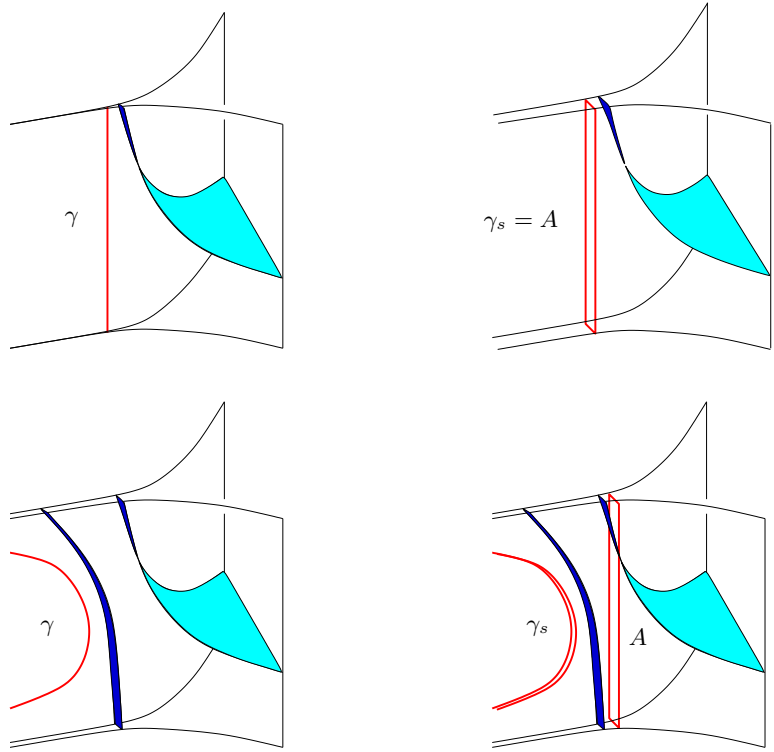}
\adjustrelabel <-5pt,0pt> {g1}{$\gamma$}
\adjustrelabel <-5pt,0pt> {g2}{$\gamma$}
\adjustrelabel <-5pt,0pt> {g1=a1}{$\gamma_s = A$}
\relabel {g2s}{$\gamma_s$}
\relabel {a2}{$A$}
\endrelabelbox}
\caption{If $\gamma$ --- a branch component of $\Lambda^+_b$ --- is a circle, it gets split open
to a tangential annulus $\gamma_s = A$. If $\gamma$ is a line, it gets split open to a tangential
rectangle $\gamma_s$, bounded from $\fG$ by a transverse annulus $A$. Polygons $P(\lambda)$ spiral
around $\gamma_s$.}
\end{figure}

We will show in fact that the branch locus is empty --- that is, that $\Lambda^+_b = \Lambda^+$.
We do this in two steps: firstly, by showing that each component of the branch locus is a circle,
and secondly by showing that there are no circle components. The first step does not use
the results we have proved about sticky laminations, and comes down to the fact that
if a branch locus is a line, then it is covered by infinitely many lines which all bound the
same interstitial region. Since the interstitial regions correspond to vertices of an
ideal polygon $P_\u$, we can construct a collection of leaves of $\Lambda^+_\u$ which
cross boundary leaves of $P_\u$, giving us a contradiction. The second step uses the
fact that $\Lambda^+_\u$ is not sticky, by considering the action of the element of the
fundamental group corresponding to a branch circle $\gamma$ on the circle at infinity
of a leaf $\lambda$ containing a cover $\til{\gamma}$. 

\begin{thm}\label{circle_branch_locus}
The branched lamination $\Lambda^+_b$ constructed in lemma~\ref{produces_genuine_lamination},
which is split open to produce $\Lambda^+$ has branch locus a (possibly empty) finite union
of circles.
\end{thm}
\begin{pf}
Since $\Lambda^+$ is genuine, and no complementary region is a product, there are only
finitely many complementary regions. It is easy to see that the same is therefore
true for $\Lambda^+_b$. Note that a branched lamination can
be naturally partitioned into guts and interstitial regions just as an ordinary
lamination can, where for any $\epsilon$, the interstitial region can be constrained
to be contained in the set of points where the
injectivity radius of the complementary region is at most $\epsilon$, and the
guts and interstices meet along annuli. In theorem~\ref{solid_torus_guts} we showed
that interstitial annuli can be straightened to be transverse or tangential to $\F$;
the same proof obviously applies to interstitial annuli in complementary
regions to branched laminations where the branch locus is tangent to $\F$, which is
exactly our situation. Each component $\gamma$ of the branch locus is in the
boundary of exactly one interstitial region $I$, meeting a gut $\fG$ along an annulus $A$,
where $A$ is transverse or tangent to $\F$. Assume $\gamma$ is noncompact.

$A$ runs between two leaves $l,m$ of $\Lambda^+_b$.
Let $\til{A}$ be a cover of $A$ in $\til{M}$, and let $\til{l},\til{m}$ be the leaves
of $\til{\Lambda}^+_b$ it runs between. Let $\alpha \in \pi_1(M)$ be the stabilizer of
$\til{A}$, corresponding to the loop which is the core of $A$.
The leaves $\til{l},\til{m}$ correspond to leaves $l_\u,m_\u$ of $\Lambda^+_\u$.
If $A$ were tangent to $\F$, then we could assume $\til{A}$ is contained in some leaf $\lambda$
of $\til{\F}$ stabilized by $\alpha$; but then $\alpha$ would preserve the distinct
geodesics $l(\lambda),m(\lambda)$ of $\lambda$,
which is absurd since $\alpha$ acts on $\lambda$ by an isometry.
So $\til{A}$ is transverse to $\F$, and $\alpha$ (say) moves points on $\til{A}$ in the
negative direction in $\til{M}$.

Let $\til{\gamma}$ denote a lift of $\gamma$ which
is a boundary component of $\til{l} \cap \til{m}$, and
suppose $\til{\gamma} \subset \mu$ for $\mu$ a leaf of $\til{\F}$. Choose a $\lambda$
with $\mu < \lambda$ and $\lambda \cap A$ nonempty.
The translates $\alpha^n(\til{\gamma})$ are distinct boundary components of
$\til{l} \cap \til{m}$, since if $\alpha$ stabilizes $\til{\gamma}$, it covers a
compact $\gamma$, contrary to assumption. Let $\mu_n = \alpha^n(\mu)$, and consider the
set of markers from $\mu_n$ to $\lambda$ for some $\lambda$ intersecting $\til{\fG}$.
By the discussion preceding corollary~\ref{finite_sided},
$l_\u,m_\u$ are adjacent sides of a finite sided
polygon $P_\u$ determining a finite sided ideal polygon $P(\lambda) \subset \lambda$.
The markers from $\mu_n$ to $\lambda$ are all contained in the arcs of circle $s_l(\lambda),
s_m(\lambda)$ bounded by $l(\lambda),m(\lambda)$, since $l(\mu_n)$ and $m(\mu_n)$ are the
same geodesic in $\mu_n$.
Moreover, there are infinitely many 
markers which limit to points in $s_l(\lambda)$ and $s_m(\lambda)$
respectively; it follows that there is a leaf $k_n(\lambda)$
of $\Lambda^+(\lambda,\mu_n)$ with endpoints
in $s_l(\lambda),s_m(\lambda)$ respectively, and which does not intersect the endpoint
$l(\lambda) \cap m(\lambda)$. Since $k_n(\lambda)$ cannot cross $l(\lambda),m(\lambda)$
transversely, it follows that the endpoints of $k_n(\lambda)$ must be the two
(nonintersecting) endpoints of $l(\lambda),m(\lambda)$, so all the
$k_n(\lambda)$ are equal. But the markers from distinct
$\mu_n$ to $\lambda$ are distinct; it follows that $k_n(\lambda)$ separates
$k_{n+1}(\lambda)$ from $k_{n-1}(\lambda)$; this contradiction implies there is only
one $\mu_n$ --- that is, $\til{\gamma}$ is stabilized by $\alpha$, and therefore
$\gamma$ is a circle. It is clear that there is at most one circle component of the
branch locus for each $A$, so the theorem follows.
\end{pf}

Lemma~\ref{limiting_leaves_disjoint} puts nontrivial constraints on ${\Lambda}^+_m$
and the action of $\pi_1(M)$. In particular, we can prove the following corollary:

\begin{cor}\label{no_tangential_annulus}
No gut region $A$ has a tangential interstitial annulus.
\end{cor}
\begin{pf}
If $G$ has such an annulus $A$, lift to $\til{G}$ with a tangential interstitial
rectangle $\til{A} \subset \lambda$ for some leaf $\lambda$ of $\til{\F}$. Then
the core $\alpha$ of $A$ fixes $\lambda$, stabilizes the two ideal points of
$\til{A} \subset \lambda$, and acts as a hyperbolic isometry.

But leaves of $\partial P$ are not isolated on the outside. It follows that the
axis of $\alpha$ is a limit of leaves $l_i$ of $\Lambda^+(\lambda)$ on one side (actually
on either side). By lemma~\ref{limiting_leaves_disjoint}, these leaves are
not asymptotic to the endpoints of the axis of $\alpha$. It follows that $\alpha^n(l_i)$
intersects $l_j$ transversely for some suitable $i,j,n$, which is absurd. The corollary
follows.
\end{pf}

Recall that hypothetical tangential interstitial annuli only occurred by splitting
open a branched $\Lambda^+_b$ along a branch circle. So corollary~\ref{no_tangential_annulus}
implies that there are no branch circles in $\Lambda^+_b$. Together with
theorem~\ref{circle_branch_locus},
this completes the proof that in fact $\Lambda^+_b = \Lambda^+$. That is, no
splitting was necessary to obtain a lamination.

\begin{cor}\label{no_splitting_necessary}
The branched lamination $\Lambda^+_b$ is actually a genuine lamination; i.e. its
branch locus is empty, and $\Lambda^+_b = \Lambda^+$.
\end{cor}

\begin{defn}
An element $\phi \in \hom(S^1)$ is {\em topologically pseudo--Anosov}
if it has at least $3$ and at most finitely many fixed points 
and translates points alternately clockwise and anticlockwise in the 
intervals complementary to these points.
\end{defn}

Since the action of a topologically pseudo--Anosov element on $S^1$
moves points in the complement of the fixed point set alternately clockwise
and anticlockwise, it must fix an {\em even} number of points.

\begin{lem}
Let $\alpha$ be the core curve of a complementary region $G$ corresponding to the
complementary polygon $P$ to $\Lambda^+_\u$. The action of some power of
$\alpha$ on $S^1_\u$ is topologically pseudo--Anosov, where
the vertices of $P$ are the repelling points and there is exactly one
attracting point in each complementary interval.
\end{lem}
\begin{pf}
Let $\lambda$ be some leaf of $\til{\F}$ intersecting $\til{G}$ in a finite
sided ideal polygon $P(\lambda)$, combinatorially equivalent of course to $P$.
For each $n,k \in \Z$ there is a commutative diagram

\begin{diagram}
S^1_\u                          & \rTo^{\alpha^k}         & S^1_\u \\
\dTo^{\phi_{\alpha^n(\lambda)}} &                         & \dTo^{\phi_{\alpha^{n+k}(\lambda)}} \\
S^1_\infty(\alpha^n(\lambda))   & \rTo^{\alpha^k}         & S^1_\infty(\alpha^{n+k}(\lambda)) \\
\end{diagram}

Now, the monotone maps 
$$\phi_{\alpha^n(\lambda)}:S^1_\u \to S^1_\infty(\alpha^n(\lambda))$$
quotient out less and less of the circle $S^1_\u$ as $n \to + \infty$.

Let $\lambda_t$ be the set of leaves of $\til{\F}$ which intersect the lifted gut
$\til{\fG}$, and
parameterize them by $t \in \R$. We know by theorem~\ref{dense_markers}
that as $t \to +\infty$, the leaves $\lambda_t$
escape to the positive end of $L$. On the other hand, as $t \to -\infty$ we don't
know whether the leaves escape to a negative end of $L$, or have a lower bound.
A lower bound here would consist of a union of incomparable uppermost leaves which
are all limits of $\lambda_t$ as $t \to -\infty$.

Fix a side $l$ of $P$ with vertices $p,q$. We will study the action of the power
of $\alpha$ fixing the vertices of $P$
on the complementary arc to $p,q$ outside $P$. For each leaf $\lambda_t$,
let $l(t)$ denote $\phi_{\lambda_t}(l)$.

For each $\mu$ an uppermost leaf below the $\lambda_t$, we can consider the
set $\Lambda^+(\lambda_t,\mu)$ for each $t$. By hypothesis we can choose $l$
so that every marker from
$\mu$ to $\lambda_t$ lies on the complementary arc to $l(t)$ outside $P(t)$, so
this lamination contains a closest
leaf $l_\mu(t)$ to $l(t)$. Notice that {\em every} leaf of $\Lambda^+(\lambda_t)$
can be approximated by leaves of $\Lambda^+(\lambda_t,\kappa)$ for some 
$\kappa < \lambda_t$; if $l_\kappa(t)$ denotes the closest leaf of
$\Lambda^+(\lambda_t,\kappa)$ to $l(t)$ for each $\kappa$, then as $\kappa_i$ increases
monotonely while staying below every $\lambda_t$, the endpoints of the
leaves $l_{\kappa_i}(t)$ move monotonely towards the endpoints of $l(t)$. It follows
that either $l_\mu(t) = l(t)$ for some uppermost $\mu$, or there
is some $\kappa$ with $l_\kappa(t)$ separating $l_\mu(t)$ from $l(t)$. Then
if $\nu$ is an uppermost leaf below the $\lambda_t$ satisfying
$\kappa \le \nu < \lambda_t$ for all $t$, then there is a leaf
$l_\nu(t)$ of $\Lambda^+(\lambda_t,\nu)$ which separates $l_\mu(t)$ from $l(t)$
(since by assumption, in this case there is no $\mu$ below every $\lambda_t$ with
$l_\mu(t) = l(t)$.)

We consider the size of $l_\mu(t)$ in the visual metric as seen from the
point $p_t = \til{\alpha} \cap \lambda_t$, where now $\til{\alpha}$ denotes the
relevant lift of the core of $\fG$. If $l_\mu(t) \ne l(t)$, there is a $\nu$
so that any leafwise geodesic $\delta_t$ from $p_t$ to a point in the convex hull of
$\Lambda^+(\lambda_t,\mu)$ must pass through the convex hull of
$\Lambda^+(\lambda_t,\nu)$. In particular, the sequence $\delta_t$ is eventually
arbitrarily close to some point in $\nu$, and therefore must be close to
$\nu$ for an arbitrarily long segment before getting to $l_\mu(t)$.
But this means that the visual angle of $l_\mu(t)$ converges to $0$ as
$t \to -\infty$. This means that under iterates of $\alpha$, the leaf 
$l_\mu$ of $\Lambda^+_\u$ converges to a point. It follows that for each edge
$l$ of $P$ and for every uppermost leaf $\mu$ below the $\lambda_t$, if
$l_\mu(t) \ne l(t)$ the leaves of $\Lambda^+_\u(\cdot,\nu)$ converge to a single
point of $S^1_\u$ under iteration of $\alpha$. 
(Here this notation is meant to indicate the leaf of $\Lambda^+_\u$
determined by the appropriate extremal pair of long markers from
$\nu$.) For $\mu,\nu$ distinct uppermost leaves less than every $\lambda_t$,
$\mu$ and $\nu$ are incomparable, so $l_\mu(t)$ and $l_\nu(t)$ can't both be equal
to $l(t)$; that is, there is {\em at most one} uppermost leaf
$\mu$ below the $\lambda_t$ whose markers limit to $l(t)$.

Now, a dense set of leaves of $\Lambda^+_\u$ outside $l$ are elements of
$\Lambda^+_\u(\cdot,\kappa)$ for some $\kappa$. If $\kappa$ is incomparable to
or equal to some $\lambda_t$, then every leaf of $\Lambda^+_\u(\cdot,\kappa)$ 
will project to a single point in $S^1_\infty(\lambda_t)$ for $t \to -\infty$.

Suppose there is {\em no} leaf $\mu$ with $l_\mu(t) = l(t)$,
by what we have just said, the visual measure of $\Lambda^+(\lambda_t,\mu)$ converges
to $0$ as seen from $p_t$ as $t \to -\infty$. It follows in this case
that for a dense set of leaves $l'$ of $\Lambda^+_\u$ outside $l$, the
images $\phi_{\lambda_t}(l') \in \lambda_t$ have visual measure going to $0$ as
$t \to -\infty$, as seen from $p_t$. Said another way, we have that for
any $l'$ outside $l$, and any $i$,
$$\lim_{n \to \infty}
\text{visual measure of }(\alpha^{n+i} \phi_{\alpha^{-n}(\lambda_0)}(l')) \text{ as seen
from } \til{\alpha} \cap \alpha^i(\lambda_0) \to 0$$
and therefore
$$\lim_{n \to \infty}
\text{visual measure of }\phi_{\alpha^i(\lambda_0)}(\alpha^{n+i}(l')) 
\text{ as seen from } \til{\alpha} \cap \alpha^i(\lambda_0) \to 0$$
But this is true for every positive $i$. Since $\alpha^i(\lambda_0)$
increases without bound as $i \to \infty$, it follows that the visual
measure of $\phi_{\mu}(\alpha^n(l'))$ converges to $0$ as seen from any point in
{\em any} $\mu$, as $n \to \infty$. In particular,
every leaf of $\Lambda^+_\u$ outside $l$
converges to a single point under iteration of $\alpha$. It follows that
$\alpha$ fixes at most one point in this complementary interval, and since
$\Lambda^+_\u$ is not sticky, $\alpha$ has a single attracting fixed point in
this interval.

Conversely, suppose that a unique uppermost leaf $\mu < \lambda_t$ for all $t$ exists
for which $l_\mu(t) = l(t)$. The uniqueness of $\mu$ implies that
$\mu$ is fixed by $\alpha$, and therefore $\alpha$
fixes exactly two points $m,n$ of $S^1_\infty(\mu) = \phi_\mu(S^1_\u)$.
The preimage of $m$, say, under $\phi_\mu^{-1}$ is the unique arc of $S^1_\u$ determined by
$l$ which collapses to $m$. We will show the
preimage of $n$ is a single point.

Since $\alpha$ fixes no leaf $\lambda_t$, we can look at the quotient
$\til{M}/\alpha$. The leaf $\mu$ covers a cylinder $\mu/\alpha$ in
this quotient, and each leaf $\lambda_t$ covers a plane (which we also refer to as
$\lambda_t$) which spirals around
$\mu/\alpha$. If $\gamma$ denotes the nontrivial geodesic loop in $\mu/\alpha$, and
$\tau$ is a transversal through some point of $\gamma$, then
$\tau$ intersects the image of each $\lambda_t$ in an infinite sequence of points
$\tau_i(t) \in \lambda_t$ going out to infinity. Fix $t$. For $i$ sufficiently
large, $\tau_i(t)$ is arbitrarily close to $\gamma$. It follows that
the leaf $\lambda_t$ stays arbitrarily close to $\mu/\alpha$ on an arbitrarily
large neighborhood of $\tau_i(t)$, and in particular, it stays close along a
geodesic arc $\gamma_i'(t) \subset \lambda_t$ from $\tau_i(t)$ to $\tau_{i+1}(t)$.
But this implies that for sufficiently large $i$, $\lambda_t$ stays arbitrarily
close to $\mu/\alpha$ along a quasigeodesic ray (which eventually is arbitrarily
close to a geodesic ray) $\gamma'(t)$. Let $n'(t) \in S^1_\infty(\lambda_t)$
denote the endpoint of this geodesic ray. Then it is clear we have actually
constructed a marker in $\til{M}$ from $\mu$ to $\lambda_t$ whose endpoints run
from $n$ to $n'(t)$. We show now that this
implies that the preimage of $n$ in $S^1_\u$ under $\phi^{-1}_\mu$
is a single point $n'$. Suppose that
the preimage of $n$ in $S^1_\u$ under $\phi^{-1}_\mu$ were an arc $k$. Then the union 
$$K = \bigcup_\nu \phi_{\nu}(k)$$
is a leaf of $\til{\Lambda}^+$ which intersects $\lambda_t$ for sufficiently
large $t$, but which does not intersect $\mu$ or leaves below $\mu$. If $n'(t)$ is
an endpoint of $\phi_{\lambda_t}(k)$, then $\gamma'(t)$ is eventually asymptotic
to $\phi_{\lambda_t}(k)$; otherwise, $\gamma'(t)$ crosses $\phi_{\lambda_t}(k)$.
Let $q \in \mu$ be a point on $\gamma(t)$ very close to
$q' \in \gamma'(t)$ for some $t$, where $q'$ is very close to $\phi_{\lambda_t}(k)$ in the first
case, or $q'$ is on the ray from $\phi_{\lambda_t}(k) \cap \gamma'(t)$ to $n'(t)$ in the second
case. There are lots of markers from
$\mu$ to $\lambda_t$; since $q$ can be chosen arbitrarily
close to $q'$, the set of endpoints of
markers from $\mu$ to $\lambda_t$ are an $\epsilon$--net in $S^1_\infty(\lambda_t)$. But
by choice of $q'$, this means that there are many distinct markers which end in the
arc of $S^1_\u$ bounded by $\phi_{\lambda_t}(k)$. On the other hand, by the
definition of $k$, there are no markers from $\mu$ to $\lambda_t$ in this arc. This
contradiction implies that there is no such $k$, and therefore the preimage of
$n$ in $S^1_\u$ under $\phi^{-1}_\mu$ is a single point $n'$, as claimed.

Now, $n'$ is a fixed point of $\alpha$ in $S^1_\u$. Let $l'$ be a leaf of
$\Lambda^+_\u$ arbitrarily close to $l$, and look at its iterates under positive
powers of $\alpha$. If there are any fixed points of
$\alpha$ between $n'$ and $l$, then $l'$ limits to a leaf $l''$ between $n'$ and $l$
fixed by $\alpha$. Then $\phi_\mu(l'')$ is fixed by $\alpha$, which is absurd, since
$\alpha$ fixes only $m$ and $n$ in $S^1_\infty(\mu)$.

It follows that $n'$ is the unique fixed point of $\alpha$ in the
complementary interval to $l$. Since leaves $\lambda_t$ spiral around $\mu/\alpha$ in
the quotient $M/\alpha$ as they wind around $\gamma$, one sees that $n'$ is an attracting
fixed point of $\alpha$ in this case too.

Consequently in either case $\alpha$ has the required
dynamics, and is topologically pseudo--Anosov.

See figure 8 for an illustration.
\end{pf}

\begin{rmk}\label{marker_structure}
Notice that the argument that $\phi_\mu^{-1}(n)$ is a single point $n' \in S^1_\u$
uses the geometric fact that for $K$ a leaf of $\til{\Lambda}^+$ which
doesn't intersect a leaf $\mu$ of $\til{\F}$, there is a {\em lower bound}
on the distance from $K$ to $\mu$. In turn, this relies on the compactness of $M$
and the fact that for $q \in \mu$ and $q' \in \lambda$ with $\mu < \lambda$ and
$q,q'$ sufficiently close, endpoints of markers from $\mu$ to $\lambda$ are
arbitrarily dense in $S^1_\infty(\lambda)$.

Notice that this implies that
if there exists a marker which is {\em asymptotically thin} --- i.e.
the vertical thickness of the marker goes to zero as we go to infinity, then
the following is true: for any pair of leaves $\mu < \lambda$ the set of
endpoints of markers ({\em not} long markers) from $\mu$ to $\lambda$, considered
as a subset of $S^1_\infty(\lambda)$, contains no points which are isolated on either
side (i.e. every point is a double--sided limit). Moreover, when we take the closure to
produce long markers, we observe that every long marker is a {\em nontrivial}
limit of markers on at least one side, and that anything isolated on exactly
one side is a long marker but {\em not} a marker. An asymptotically thin
marker exists if for some closed geodesic loop in a leaf of $\F$, there is
{\em strictly contracting} holonomy on at least one side. 
There are certain points in this paper where the logic would be
streamlined if we assumed the existence of an asymptotically thin marker.
\end{rmk}

For a topologically pseudo--Anosov $\alpha$ associated to a complementary
polygon $P$, there is a ``dual'' polygon $P'$ with vertices in $S^1_\u$
the attracting fixed points of $\alpha$. Let $\gamma$ be one of the
boundary edges of $P'$.

\begin{figure}[ht]
\scalebox{.6}{\includegraphics{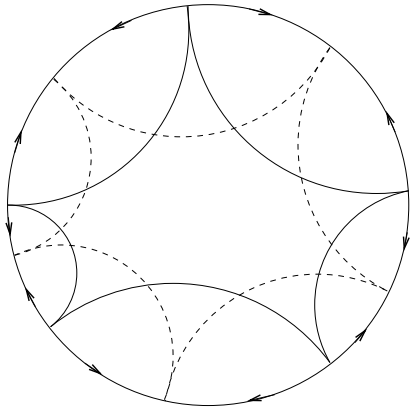}}
\caption{The action of the core of a gut region on $S^1_\u$ is topologically
pseudo--Anosov. The dual polygon $P'$ is illustrated in dashed lines.}
\end{figure}

\begin{lem}\label{weak_pinching_lamination}
Let $M$ be atoroidal. For every $\beta \in \pi_1(M)$,
$\beta(\gamma)$ is either equal to or disjoint from
$\gamma$.
\end{lem}
\begin{pf}
Associated to $\gamma$ there is a plane $\pi$ properly embedded in $\til{M}$
which is just the union of $\phi_\mu(\gamma) \subset \mu$ as we vary over
all leaves $\mu$ of $\til{\F}$. This is a single connected plane, since for every
leaf $\lambda$ of $\til{\F}$, $\pi$ is either a single geodesic or empty,
and for $\lambda > \mu$ a point in $\mu \cap \pi$ is joined to a point in
$\lambda \cap \pi$. The foliation of $\pi$ is certainly not $\R$--covered,
since it intersects $\til{\Lambda}^+$ transversely, and therefore intersects
incomparable leaves of $\til{\F}$.

An appropriate power of $\alpha$ stabilizes $\pi$ and
quotients it out to give a cylinder $C$ which projects to $M$. We claim this
projection is a covering map. For, if $\beta$ acting on $S^1_\u$ translates
$\gamma$ to intersect itself, this intersection determines a line of 
intersection of $\pi$ with $\beta(\pi)$. As we move along this line of
intersection in the positive direction and translate back to a fixed region 
in $\til{M}$ by iterates of $\alpha^{-1}$, the action of $\alpha^{-1}$ 
on $S^1_\u$ keeps the endpoints of $\beta(\gamma)$ away from the 
endpoints of $\gamma$. It follows
that the line of intersection of $\pi$ with $\beta(\pi)$ must project to
a compact portion of the cylinder $C$. The projection of $C$ to $M$ is
an immersion, so it makes sense to talk about a path in $M$ contained in
a sheet of the projection. The projection of $C$ to $M$ is ``locally proper'',
in the sense that a path in $M$ contained in a sheet which exits to an end of $C$
must be infinitely long. It follows that a self--intersection of $C$ in $M$
is either compact, or lifts to a properly embedded curve in $C$. Since one end
of this line of intersection lifts to a bounded region of $C$, the entire line does,
and therefore it is periodic under the action of $\alpha$. In particular, $\alpha, \beta$
generate the image of a Baumslag--Solitar group. By considering the action of
$\alpha$ and $\beta$ on $S^1_\u$ one sees by a standard argument
that this is absurd. For more details one can consult the
analogous lemma 5.3.7 in \cite{dC99b}. Or, one can appeal directly to a theorem
of Shalen (\cite{Shalen_BS})
on the possible homomorphisms of Baumslag--Solitar groups into 
$3$--manifold groups; such images are all virtually abelian. Since $M$ is
atoroidal, one can conclude that the image is actually cyclic.

It follows that no such $\beta$ can exist.
\end{pf}

\begin{rmk}
The technique of proof of lemma~\ref{weak_pinching_lamination} is based on
a sketch of an argument by Thurston proving similar facts about $\R$--covered foliations.
\end{rmk}

\begin{defn}
We make the definition 
$$\Lambda^-_\u = \overline{\bigcup_\alpha \alpha(\lambda)}$$
It is clear from lemma~\ref{weak_pinching_lamination} that $\Lambda^-_\u$
is a $\pi_1(M)$--equivariant
lamination of $S^1_\u$. It therefore determines a possibly branched lamination
of $\til{M}$ which can be split open to
a $\pi_1(M)$--equivariant lamination $\til{\Lambda}^-$ of $\til{M}$ covering
an essential lamination of $M$ which we refer to as $\Lambda^-$.
\end{defn}

\begin{thm}
Let $M$ be atoroidal. Then $\Lambda^-$ has solid torus guts. 
Moreover, a lift of any one of these gut 
regions determines a finite sided ideal polygon in $S^1_\u$ which is 
exactly a ``dual'' polygon to a polygon associated to a gut region of 
$\Lambda^+$, and all such dual polygons arise in this way. The pair of
laminations $\Lambda^\pm$ bind every leaf of $\til{\F}$ so that complementary
regions are finite sided compact polygons.
\end{thm}
\begin{pf}
We inherit the notation $P,P'$ from above, and let $v$ denote the vertex
of $P$ isolated from its neighbors by $\gamma$. Let $e_1,e_2$ be the edges of
$P$ adjacent to $v$, so that $\gamma$ intersects both $e_1$ and $e_2$. Let
$w_1,w_2$ be the other vertices of $e_1,e_2$ so that $v,w_1,w_2$ are in
clockwise order on $S^1_\u$.

Clearly, the proof of theorem~\ref{solid_torus_guts} applies to $\Lambda^-$,
since it too is transverse to $\til{\F}$. So $\Lambda^-$ has solid torus guts.
The first claim is that
$\gamma$ is a ``boundary'' leaf of $\Lambda^-$. Firstly, $\gamma$ is isolated on
the side moving away from $v$. For if not, there is
some $\beta \in \pi_1(M)$ which moves an endpoint of $\gamma$ very slightly
off itself in the direction moving away from $v$.

Define
$$\gamma' = \lim_{n \to -\infty} \alpha^n\beta(\gamma)$$

If $\Lambda^-_\u$ is sticky, $\gamma'$ joins a vertex of $P'$ to a vertex, say
$w_2$, of $P$. From the vantage point of a point on a leaf $\lambda$ of $\til{\F}$
far out along the geodesic from $v$ to $w_2$, the geodesics from $w_2$ to $v$
and from $w_2$ to $w_1$ are almost asymptotic, and parallel to the image of $\gamma'$.
By translating back to a compact region of $\til{M}$ and
extracting a convergent subsequence, we see that $\Lambda^-_\u$ contains
a leaf of $\Lambda^+_\u$. 

If $\Lambda^-_\u$ is not sticky, $\gamma'$ is a diagonal of $P$,
joining $w_1$ to $w_2$. But all of $\Lambda^+_m$ is
in the closure of the translates of any diagonal of $P$, so 
$\Lambda^+_m \subset \Lambda^-_\u$.

In either case, $\Lambda^+_m \subset \Lambda^-_\u$ which is absurd since $\gamma$ is
transverse to $\Lambda^+_m$. 

It follows that $\gamma$ is isolated on the side away from $v$. If $\gamma$ is
isolated on the side towards $v$, then it is an interior edge of a polygon $P''$.
Since the endpoints of $\gamma$ are fixed by $\alpha$, every vertex of $P''$
is fixed by $\alpha$. It follows that there is only one vertex of $P''$ on the
side of $\gamma$ towards $v$, which is $v$ itself. But then there are edges
of $P''$ joining a vertex of $\gamma$ to $v$, which as before implies 
$\Lambda^+_m \subset \Lambda^-_\u$, which is absurd. 

It follows that there is a polygonal region
$P''$ in $S^1_\u$ of which $\gamma$ is a boundary leaf, corresponding to a
gut region of $\Lambda^-$. Since every leaf of $\Lambda^-$ is a limit of translates of
$\gamma$, there are no isolated leaves in $\Lambda^-$, and therefore $\Lambda^-$ is
actually minimal. Since $\alpha(P'')$ cannot be transverse to
$P''$, the dynamics of $\alpha$ implies the vertices of $P''$ are a
subset of the vertices of $P'$. If the vertices of $P''$ are a {\em strict} subset of
the vertices of $P'$, it follows that there is an edge of $P''$ which intersects
two non--adjacent sides of $P$. Since $\Lambda^-$ is minimal, this edge of
$P''$ is a limit of nearby leaves of $\Lambda^-$; translating these by powers of
$\alpha$, one sees that either $\Lambda^-$ is sticky and there is some
leaf of $\Lambda^-$ from a vertex of $P''$ to a vertex of $P$, or there is
some leaf of $\Lambda^-$ which is a diagonal of $P$. As before, this
implies $\Lambda^+_m \subset \Lambda^-_\u$, which gives a contradiction. In
particular, $P'' = P'$.

Finally, suppose there is some polygon $Q$ associated to a gut region of
$\Lambda^+$ with stabilizer $\delta$ also stabilizing the dual polygon $Q'$,
where $Q'$ is {\em not} a polygon associated to a gut region of
$\Lambda^-$. Then no leaf of $\Lambda^-_\u$ can link a leaf of $Q$ or $Q'$.
For, if it does, then under iteration of $\delta$ or $\delta^{-1}$ it would either
converge to a leaf of $Q$ or $Q'$ contrary to assumption, or else it would
converge to a diagonal of $Q$ or $Q'$. But the images of any diagonal of $Q$
contain a subsequence converging to any leaf of $\Lambda^+_m$, and similarly
for the images of any diagonal of $Q'$, again contrary to assumption.

It follows that the complementary region to $\Lambda^-$
intersecting $Q$ has infinite hyperbolic area. But this contradicts
theorem~\ref{solid_torus_guts}. This contradiction proves the theorem.
\end{pf}

\begin{rmk}
Since the complementary regions of $\Lambda^-_\u$ correspond to polygons
$P'$ dual to polygons $P$ associated to $\Lambda^+_\u$, it follows immediately
that $\Lambda^-$ is genuine.
\end{rmk}

\begin{defn}
Let $X$ transverse to $\F$ be a vector field, where $\F$ branches only in
the negative direction. Say that $X$ is {\em regulating in the positive
direction} if each integral curve of $\til{X}$ (the lift of $X$ to $\til{M}$)
contains a ray --- i.e. a subset homeomorphic to $\R^+$ --- which projects
properly to the leaf space of $\til{\F}$.
\end{defn}

\begin{cor}
If $\F$ is a taut foliation with one--sided branching in the negative 
direction, there is a transverse vector field $X$ which is regulating in the
positive direction.
\end{cor}
\begin{pf}
The construction of $X$ from $\Lambda^\pm$ is standard. One method is to
use the leafwise hyperbolic metric to canonically identify a point on one
leaf with a point on any higher leaf by the stratification of $\mu$ into
leaves and complementary regions of $\Lambda^\pm(\mu)$. Another method
is to ``quotient out'' the complementary regions to the laminations to give
a pair of branched foliations of $M$ which intersect each other transversely
in one--dimensional leaves whose tangent vectors give $X$. See for instance
\cite{dC99b} or \cite{wT97}.

To see the flow is regulating in the positive direction, suppose that
$\gamma^+$ is a positive integral curve of $\til{X}$. Then for each
leaf $\lambda$ of $\til{\F}$, $\gamma^+\cap \lambda$ is bounded by 
a finite collection of leaves of $\Lambda^\pm(\lambda)$. The corresponding
leaves are nondegenerate in $\Lambda^\pm(\mu)$ for all $\mu>\lambda$,
and therefore $\gamma^+$ cannot ``escape to infinity'' and must intersect
every leaf above $\lambda$. That is, $\gamma^+$ is regulating in the
positive direction.
\end{pf}

\begin{cor}
$X$ as constructed above gives a pseudo--Anosov flow transverse to $\F$.
\end{cor}
\begin{pf}
This is also standard from the construction, and from the dynamics of
the flow in a neighborhood of the closed orbits. The singular closed orbits
correspond exactly to the gut regions of $\Lambda^\pm$, and therefore
have pseudo--Anosov dynamics. To understand the dynamics near the nonsingular
closed orbits, we appeal to general facts. Theorem 5.3.15 of
\cite{dC99b} applies in our situation, and implies that every element
$\alpha$ of $\pi_1(M)$ acts on $S^1_\u$ in a manner either conjugate to an
element of $PSL(2,\R)$, or some finite power is topologically pseudo--Anosov. 
We give a sketch of the argument here: since the laminations
$\Lambda^\pm_\u$ are transverse, every point $p \in S^1_\u$ is a limit of a nested
sequence of leaves either of $\Lambda^+_\u$ or of $\Lambda^-_\u$. It follows
that for any $\alpha \in \pi_1(M)$, the set of fixed points of $\alpha$ is either
all of $S^1_\u$ (in which case $\alpha$ is the identity), or is finite.
Moreover, since each fixed point is the limit of a nested sequence of leaves of
one lamination, the dynamics of $\alpha$ at that fixed point is either attracting
or repelling (i.e. it can't look like a parabolic fixed point). In particular,
$\alpha$ is either conjugate to an element of $PSL(2,\R)$, or is
topologically pseudo--Anosov. In particular, the flow in the neighborhood of the
closed orbits has the appropriate dynamics, and $X$ can be blown down to
be pseudo--Anosov.
\end{pf}

We remark that this implies that the lamination $\Lambda^-_\u$ is not sticky
either.

\subsection{Deforming foliations}

\begin{thm}\label{onesided_family}
Suppose $\F_t$ is a one--parameter family of taut foliations with
one--sided branching of an
atoroidal $3$--manifold $M$. Then the representation
$$\pi_1(M) \to \hom(S^1_\u(\F_\epsilon))$$
is constant. Moreover, the laminations $\Lambda^\pm_\epsilon$ are
all isotopic.
\end{thm}
\begin{pf}
A continuous family $\F_t$ of foliations of $M$ may be thought of as a foliation
of $M \times I$ by Riemann surface leaves. Therefore Candel's theorem applies with
parameters, and shows that the hyperbolic structures on $\F_t$ may be chosen to vary
continuously. For small $\epsilon$, $\Lambda^\pm_0$ is transverse to
$\F_t$ and leaves of $\til{\Lambda}^\pm_0$ intersect leaves of
$\til{\F}_\epsilon$ quasigeodesically. 

After straightening these leafwise,
we can assume $\til{\Lambda}^\pm_0$ intersects leaves of $\til{\F}_\epsilon$
geodesically. It is not hard to see that this straightening is
continuous. For, if $\gamma$ is a $\delta$--quasigeodesic in $\H^2$ for some
fixed $\delta$ and $\gamma_0$ is a subarc of length $t$,
the geodesic representatives of $\gamma$ and of $\gamma_0$ agree
to within Hausdorff distance $o(e^{-t})$ on most of the length of
$\gamma_0$. If $\lambda_\epsilon^t$ is a family of leaves of
$\til{\F}_\epsilon$ and $\mu$ is a leaf of $\til{\Lambda}^+_0$, say,
then the leafwise intersections $\gamma_t$ vary continuously in the
geometric topology on arbitrarily large compact subsets. Since the
straightenings on compact subsets are continuous, the straightenings are
continuous on all of $\mu$. Notice that straightening
preserves the cyclic order on the ends of quasigeodesic rays.

For each $\epsilon$, let $L_\epsilon$ be the leaf space of $\til{\F}_\epsilon$.
Let $\alpha_i$ be a sequence
of elements in $\pi_1(M)$ which blow up some $I \subset L_\epsilon$ to
a bi--infinite properly embedded ray $r \subset L_\epsilon$. Let
$E_\infty(\epsilon)$ denote the circle bundle at infinity of
$\til{\F}_\epsilon$, thought of as a union $\bigcup_\tau UT\til{\F}_\epsilon|_\tau$
for transversals $\tau$ projecting to arcs in $L_\epsilon$.
Recall $\mu$ was a straightened leaf of $\til{\Lambda}^+_0$, transverse to
$\til{\F}_\epsilon$. Since $\mu$ is transverse to $\til{\F}_\epsilon$, it intersects
some line worth of leaves of $L_\epsilon$. Let $\tau$ be a (possibly infinite)
transversal to $\til{\F}_\epsilon$, intersecting exactly this set of leaves, contained
in $\mu$. For each point $p \in \tau \cap \lambda$ where $\lambda$ is a leaf of
$\til{\F}_\epsilon$, the geodesic $\mu \cap \lambda$ is split into two rays
emanating from $p$, which determine a pair of endpoints in $S^1_\infty(\lambda)$.
As we vary over points in $\mu$, this determines a pair
of transversals to $E_\infty(\epsilon)$, contained in the image of
$UT\til{\F}_{\epsilon}|_\tau \subset E_\infty(\epsilon)$. It is clear that
these transversals do not depend on the choice of $\tau$. In particular, for
each endpoint $p$ of a leaf $m$ of $\Lambda^+_0$ associated to a leaf $\mu$ of
$\til{\Lambda}^+$, and each transversal $\tau$ to $\til{\F}_\epsilon$
which intersects exactly the leaves of $\til{\F}_\epsilon$,
which $\mu$ intersects, there is an associated transversal
$$s(p) \subset E_\infty(\epsilon)$$
which intersects exactly the circles at infinity of the leaves intersecting $\tau$.
As we vary over all such $p$ we get a system of transversal to $E_\infty(\epsilon)$.
Since straightening preserves the cyclic order of ends of rays, these transversals
cannot cross each other (though {\it a priori} they might run into each other).

Pick some such $s = s(p)$.
For some subsequence of $\alpha_i$, the transversals
$\alpha_i(s)$ converge in $E_\infty(\epsilon)$ to a bi--infinite properly embedded
line $s' = \subset E_\infty(\epsilon)$. By passing to a further subsequence, we
can assume the limit of $\alpha_i(p)$ converges in $(S^1_\u)_0$ to some $q$. We
call $s' = s'(q)$. Then for each point $r \in (S^1_\u)_0$ in the orbit of $q$
we have a corresponding bi--infinite properly embedded line $s'(r)$ in 
$E_\infty(\epsilon)$. As before, if $q_1,q_2,q_3$ is a positively cyclically
ordered triple in the orbit of $q$ in $(S^1_\u)_0$, the cyclic ordering on 
$s'(q_1),s'(q_2),s'(q_3)$ is positive or degenerate, as measured in sufficiently
positive leaves of $\til{\F}_\epsilon$ which these three lines intersect.

If we knew these lines were ``vertical'' in $E_\infty(\epsilon)$,
in the sense that they did not cross
endpoint transversals of long markers for $\til{\F}_\epsilon$, we could identify
them with points in $(S^1_\u)_\epsilon$. This idea suggests the following strategy.

For each leaf $\lambda$ of $\til{\F}_\epsilon$ in the support of $s(r)$,
$r_\lambda = s(r) \cap S^1_\infty(\lambda)$
determines a (possibly degenerate)
interval 
$$\phi_\lambda^{-1}(r_\lambda) = I_r(\lambda) \subset (S^1_\u)_\epsilon$$
As we vary $\lambda$ in $\til{\F}_\epsilon$ in the positive direction,
these intervals or points do not necessarily vary continuously (they might
shrink discontinuously as we move in the positive direction or grow discontinuously
as we move in the negative direction), but their union as we vary over the support
of $s(r)$ is easily seen to be a {\em connected}
set $I_r \subset (S^1_\u)_\epsilon$. We claim $I_r$ is not all of $(S^1_\u)_\epsilon$.
For if $s_0(r)$ denotes the restriction of $s(r)$
to some interval $I \subset L_\epsilon$ with 
$$\bigcup_{\nu \in I} I_r(\nu) = (S^1_\u)_\epsilon$$
then we can compress $I$ by the action of some $\beta \in \pi_1$ to an arbitrarily
small interval $J \subset L_\epsilon$, by corollary~\ref{compressible_action},
and still have 
$$\bigcup_{\nu \in J} I_{\beta(r)}(\nu) = (S^1_\u)_\epsilon$$
 If $\beta_i$ is a sequence
with $\beta_i(J) \to \lambda$ for $\lambda$ a single leaf of $\til{\F}_\epsilon$,
then the sequence of transversals $\beta_i(s_0)$ converge geometrically to
$S^1_\infty(\lambda) \subset E_\infty(\epsilon)$. In particular, if
$s'$ arising from some $\mu'$ crosses $S^1_\infty(\lambda)$ transversely,
then $s'$ intersects $\beta_i(s_0)$ transversely for some $i$,
which is impossible.

Since $I_r$ is not all of $(S^1_\u)_\epsilon$, it makes sense to define
$\rho(s)$ to be the anticlockwisemost $\limsup$ of the $I_r(\lambda)$ as $\lambda$
goes to infinity in the positive direction.

A varying sequence of points in $(S^1_\u)_0$ can be approximated by a varying
sequence of endpoints of leaves of $(\Lambda^+_\u)_0$ which determine a
varying sequence of transverse sections of $(E_\infty)_\epsilon$. These sections
intersect sufficiently high leaves $\lambda \in \til{\F}_\epsilon$. A triple of points
in $(S^1_\u)_0$ corresponding to endpoints
of leaves of $(\Lambda^+_\u)_0$ determines a distinct 
triple of points in $S^1_\infty(\lambda)$ for sufficiently high leaves $\lambda$ of
$\til{\F}_\epsilon$ with the same cyclic order. We now have a
$\pi_1(M)$--equivariant
map $\rho$ from the orbit of $q$ to $(S^1_\u)_\epsilon$ which does not reverse,
but might degenerate, the cyclic order of triples of points. Since the orbit of
$q$ is dense, we can extend this map by taking anticlockwise limits to
all of $(S^1_\u)_0$, which is not {\it a priori} continuous, but is anticlockwise
semicontinuous (i.e continuous from the anticlockwise side). 

This determines a monotone relation between $(S^1_\u)_0$ to
$(S^1_\u)_\epsilon$ --- i.e. a single $S^1$ mapping monotonely to both. The master
$S^1$ is obtained by blowing up the points of $(S^1_\u)_0$ where $\rho$ is discontinuous,
and inserting the image in $(S^1_\u)_\epsilon$ there.
The whole picture admits a $\pi_1(M)$-equivariant
action, restricting to the standard action on either $S^1_\u$.
But the action of $\pi_1(M)$ on either circle is
minimal. It follows that the relation is actually an {\em isomorphism}, since if there are
nontrivial intervals in $(S^1_\u)_0$ whose preimages map to a
point in $(S^1_\u)_\epsilon$ by this relation,
the closure of the set of endpoints of such intervals is a proper invariant set,
contrary to minimality. It follows that this map is an equivariant isomorphism,
and therefore these two actions are conjugate.

The laminations $(\Lambda^\pm_\u)_\epsilon$ are determined by the action of
$\pi_1(M)$, and these determine the topology of the $(\Lambda^\pm)_\epsilon$.
We show that this implies these two laminations are isotopic. We produce a
{\em new} $\pi_1(M)$--equivariant homeomorphism $\phi_\epsilon$
from $\til{M}$ to itself, covering a
self map of $M$ homotopic to the identity,
taking $\til{\Lambda}^\pm_0$ to $\til{\Lambda}^\pm_\epsilon$. Since $M$ contains
very full laminations, any self--homeomorphism homotopic to the identity is actually
isotopic to the identity, by the main result of \cite{dGwK93}.
Since both $\til{\Lambda}^\pm_0$ and $\til{\Lambda}^\pm_\epsilon$ can be canonically
collapsed to singular foliations and recovered from these foliations, we assume
that both collapsings have been done. Then a point $p \in \til{M}$ is determined
by which leaf $\lambda$
of $\til{\F}_0$ it lies on, and which (possibly singular)
leaves $m,l$ of $\phi_\lambda(\Lambda_\u^\pm)_0 \cap \lambda$
it lies on. There is an identification of $m,l \in (\Lambda_\u^\pm)_0$ with
corresponding leaves $m',l' \in (\Lambda_\u^\pm)_\epsilon$.
Moreover, $p$ lies on a leaf $\lambda_\epsilon$ of $\til{\F}_\epsilon$.
So define
$$\phi_\epsilon(\lambda,m,l) = (\lambda_\epsilon,m',l')$$
It is clear that this satisfies the required properties.
\end{pf}

\begin{thm}
Let $\F_t$ be a one--parameter family of taut foliations of an
atoroidal $3$--manifold $M$, where
$\F_0$ is $\R$--covered, and $\F_t$ has one--sided branching for
$t>0$. Then the representation
$$\pi_1(M) \to \hom(S^1_\u(\F_\epsilon))$$
is constant for $t\ge 0$. Moreover, the laminations $\Lambda^\pm_\epsilon$
are all isotopic.
\end{thm}
\begin{pf}
The notation and construction of the laminations $\Lambda^\pm$ and
circle $S^1_\u$ arising from an $\R$--covered foliation are consistent
with the notation used throughout this paper. The existence of such a
circle and such laminations for an $\R$--covered foliation
is established in \cite{dC99b}. In any case, by exactly repeating the
argument of theorem~\ref{onesided_family} we obtain a
$\pi_1(M)$--equivariant map $(S^1_\u)_0 \to (S^1_\u)_\epsilon$.
Minimality again implies that this map is an isomorphism, and therefore
that the representations are conjugate and the laminations are isotopic.
\end{pf}

\section{Examples}

\begin{exa}[Meigniez]
The first examples of taut foliations with one--sided branching of hyperbolic
$3$--manifolds are constructed by Meigniez in \cite{gM91}. These
foliations are deformations of surface bundles over circles with
fiber $\Sigma$ and monodromy some pseudo--Anosov
$\phi:\Sigma \to \Sigma$, where the
leaves of the foliation are bent in one direction along one of the
transverse laminations obtained by suspending the stable/unstable
laminations of $\Sigma$ fixed by $\phi$.
\end{exa}

\begin{exa}
If $M$ is a graph manifold obtained by plumbing together a
collection of Seifert fibered manifolds along tori, and $\F$ is a foliation
transverse to the circle fibers of each piece and transverse to the separating
tori, $\F$ is $\R$--covered. The restriction of $\F$ to each Seifert--fibered
piece arises from a slithering of that piece over $S^1$ (see \cite{wT97}).
Moreover, $\pi_1(M)$ acts on $L$, the
leaf space of $\til{\F}$ by {\em coarsely projective transformations}.

Here we say a homeomorphism $\phi:\R \to \R$ is a coarse projective transformation
if there is some constant $k>0$ and a scale factor $\alpha >0$ such that
for all $p,q \in \R$ we have
$$\alpha(\phi(p) - \phi(q)) - k \le p - q \le \alpha(\phi(p) - \phi(q)) + k$$

In such a situation, it is relatively easy to find embedded loops $\gamma$
transverse to $\F$ whose lifts $\til{\gamma}$ to $\til{M}$ are regulating
in the positive direction but {\em not} the negative direction, and whose
complements are hyperbolic. It follows that for $f:N \to M$ a finite cover branched
over $\gamma$ with sufficiently large branching degree, $N$ is a hyperbolic manifold,
and the pulled back foliation $\G = f^{-1}(\F)$ is taut and has one--sided branching
in the negative direction. This example is an adaptation of a construction of
nonuniform $\R$--covered foliations of hyperbolic $3$--manifolds in \cite{dC99a} and
is explained in more detail in \cite{CalPhD}.

Instead of looking at a branched cover, if one does not mind treating manifolds
with cusps, one can just let $N$ be the complement $M - \gamma$.
\end{exa}

\begin{exa}
The knot complement $M = S^3 - K$ where $K$ is the knot $5_2$ 
from Rolfsen's table is an example of the previous kind.
The manifold $M$ is a $2$--fold branched cover of a graph manifold.     

There is a {\em taut ideal triangulation} (see \cite{Lack}) of $M$ carrying
this foliation with a projectively invariant transverse measure. The limit set
of a leaf of $\til{\F}$ is a dendrite, as illustrated in the figure below. Thanks to
Nathan Dunfield for producing this picture.
\end{exa}

\begin{figure}[ht]
\scalebox{.6}{\includegraphics{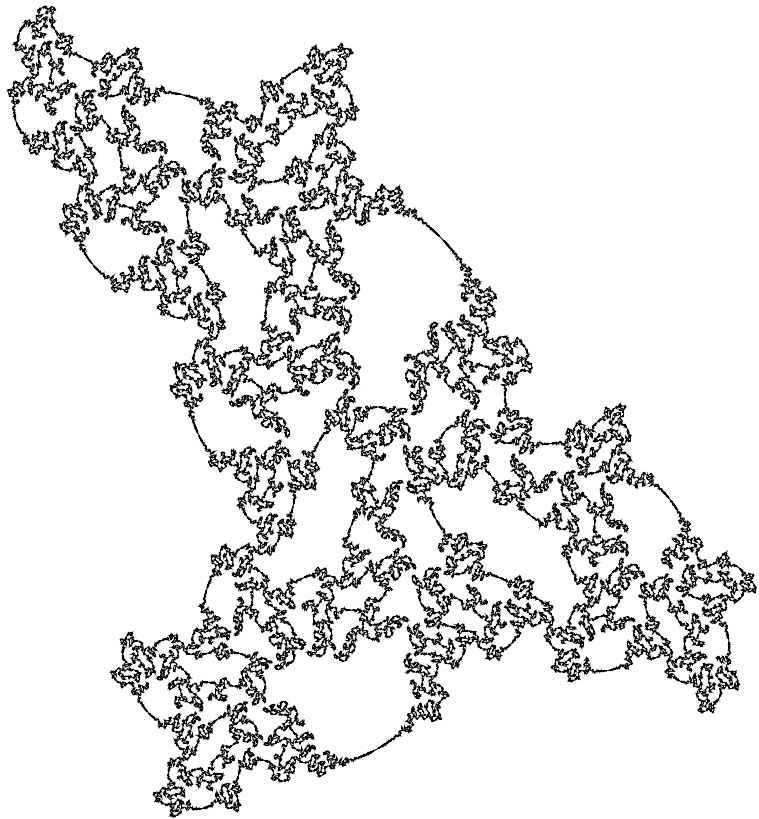}}
\caption{The limit set of a leaf of $\til{\F}$ is a {\em dendrite}.}
\end{figure}

\end{document}